\documentclass{gtpart}

\usepackage{graphicx}
\usepackage[latin1]{inputenc}
\usepackage{enumerate}
\usepackage[all]{xy}
\usepackage{multirow}
\usepackage{color}
\usepackage[T1]{fontenc}
\usepackage{hyperref}

\title{Finite approximations of $p$-local compact groups}
\author{A. Gonzalez}
\givenname{Alex}
\surname{Gonzalez}
\address{Department of Mathematics\\
Kansas State University\\
66506 Manhattan KS\\
United States of America}
\email{agondem@math.ksu.edu}
\urladdr{http://www.math.ksu.edu/~agondem}


\arxivreference{arXiv:1503.08856}
\arxivpassword{hnf9m}

\volumenumber{}
\issuenumber{}
\publicationyear{}
\papernumber{}
\startpage{}
\endpage{}
\doi{}
\MR{}
\Zbl{}
\received{}
\revised{}
\accepted{}
\published{}
\publishedonline{}
\proposed{}
\seconded{}
\corresponding{}
\editor{}
\version{}

\newtheoremstyle{slant}{}{}{\slshape}{}{\bfseries}{.}{.5em}{}%
\newtheoremstyle{special}{}{}{\slshape}{}{\bfseries}{.}{.5em}{\thmnote{#3}}

\theoremstyle{plain}

\newtheorem{Th}{Theorem}
\newtheorem{thm}{Theorem}[section]
\newtheorem{prop}[thm]{Proposition}
\newtheorem{lmm}[thm]{Lemma}
\newtheorem{cor}[thm]{Corollary}

\theoremstyle{definition}
\newtheorem{rmk}[thm]{Remark}{\rm}
\newtheorem{defi}[thm]{Definition}
\newtheorem{expl}[thm]{Example}{\rm}
{\rm}
\newtheorem{hyp}[thm]{Hypothesis}{\rm}
{\rm}

\DeclareMathAlphabet\EuR{U}{eur}{m}{n}
\SetMathAlphabet\EuR{bold}{U}{eur}{b}{n}

\newcommand{\3}[1]{\widebar{#1}}
\newcommand{\4}[1]{\widetilde{#1}}
\newcommand{\5}[1]{\widehat{#1}}

\newcommand{\res}{\mathrm{res}}
\newcommand{\Inj}{\operatorname{Inj}\nolimits}

\newcommand{\incl}{\mathrm{incl}}
\newcommand{\Hom}{\mathrm{Hom}}
\newcommand{\Mor}{\mathrm{Mor}}
\newcommand{\Map}{\mathrm{Map}}

\newcommand{\Rep}{\mathrm{Rep}}
\newcommand{\Aut}{\mathrm{Aut}}
\newcommand{\Inn}{\mathrm{Inn}}
\newcommand{\Out}{\mathrm{Out}}
\newcommand{\Iso}{\mathrm{Iso}}

\newcommand{\Ob}{\mathrm{Ob}}
\newcommand{\Ker}{\mathrm{Ker}}
\newcommand{\Id}{\mathrm{Id}}
\newcommand{\Syl}{\operatorname{Syl}\nolimits}

\newcommand{\hocolim}{\mathrm{hocolim \,}}

\newcommand{\ind}{\mathrm{ind}}

\newcommand{\F}{\mathbb{F}}
\newcommand{\N}{\mathbb{N}}

\newcommand{\bb}{\mathcal{B}}
\newcommand{\cat}{\mathcal{C}}
\newcommand{\ee}{\mathcal{E}}

\newcommand{\FF}{\mathcal{F}}
\newcommand{\g}{\mathcal{G}}
\newcommand{\hh}{\mathcal{H}}

\newcommand{\kk}{\mathcal{K}}
\newcommand{\LL}{\mathcal{L}}
\newcommand{\mm}{\mathcal{M}}
\newcommand{\nn}{\mathcal{N}}
\newcommand{\oo}{\mathcal{O}}
\newcommand{\pp}{\mathcal{P}}

\newcommand{\TT}{\mathcal{T}}

\newcommand{\frakx}{\mathfrak{X}}

\newcommand{\colim}{\operatornamewithlimits{colim}}
\newcommand{\holim}{\operatornamewithlimits{holim}}

\newcommand{\curs}{\EuR}
\newcommand{\widebar}[1]{\overset{\mskip3mu\hrulefill\mskip3mu}{#1} \vphantom{#1}}
\newcommand{\proj}{\operatorname{proj}\nolimits}

\newcommand{\gen}[1]{\langle #1 \rangle}
\newcommand{\defin}{\stackrel{def} =}

\newcommand{\padic}{\Z^\wedge_p}
\newcommand{\prufferp}{\Z/p^{\infty}}

\newcommand{\typ}{\mathrm{typ}}
\newcommand{\ploc}{(S, \FF, \LL)}
\newcommand{\bullf}{(-)^{\bullet}}
\newcommand{\wL}{\widetilde{\LL}}
\newcommand{\Cat}{\curs{\mathbf{Cat}}}

\newcommand{\map}{\operatorname{Map}\nolimits}

\renewcommand{\Im}{\operatorname{Im}\nolimits}

\newcommand{\longleft}[1]{\;{\leftarrow%
\count255=0 \loop \mathrel{\mkern-6mu}%
    \relbar\advance\count255 by1\ifnum\count255<#1\repeat}\;}
\newcommand{\longright}[1]{\;{\count255=0 \loop \relbar\mathrel{\mkern-6mu}%
    \advance\count255 by1\ifnum\count255<#1\repeat\rightarrow}\;}
\newcommand{\Right}[2]{\overset{#2}{\longright#1}}
\newcommand{\RIGHT}[3]{\mathrel{\mathop{\kern0pt\longright#1}
        \limits^{#2}_{#3}}}
\newcommand{\Left}[2]{{\buildrel #2 \over {\longleft#1}}}
\newcommand{\LEFT}[3]{\mathrel{\mathop{\kern0pt\longleft#1}\limits^{#2}_{#3}}
}
\newcommand{\dRIGHT}[3]{\mathrel{%
   \mathop{\vcenter{\baselineskip=0pt\hbox{$\kern0pt\longright#1$}%
   \hbox{$\kern0pt\longright#1$}}}\limits^{#2}_{#3}}}
\newcommand{\LRIGHT}[3]{\mathrel{%
   \mathop{\vcenter{\baselineskip=0pt\hbox{$\kern0pt\longleft#1$}%
   \hbox{$\kern0pt\longright#1$}}}\limits^{#2}_{#3}}}
\newcommand{\RLEFT}[3]{\mathrel{%
   \mathop{\vcenter{\baselineskip=0pt\hbox{$\kern0pt\longright#1$}%
   \hbox{$\kern0pt\longleft#1$}}}\limits^{#2}_{#3}}}
\newcommand{\onto}[1]{\;{\count255=0 \loop \relbar\joinrel
    \advance\count255 by1
    \ifnum\count255<#1 \repeat \twoheadrightarrow}\;}

\numberwithin{equation}{section}

\begin{document}

\begin{abstract}    

In this paper we show how every $p$-local compact group can be described as a telescope of $p$-local finite groups. As a consequence, we deduce several corollaries, such as a Stable Elements Theorem for the mod $p$ cohomology of their classifying spaces, and a generalized Dwyer-Zabrodsky description of certain related mapping spaces.

\end{abstract}

\maketitle

For a fixed prime $p$, the concept of $p$-local compact group was introduced in the last decade by C. Broto, R. Levi and B. Oliver in \cite{BLO3}. It provides a unifying, categorical language to study several classes of groups from a $p$-local point of view, including finite groups, compact Lie groups, $p$-compact groups and algebraic groups over algebraic closures of fields of characteristic other than $p$.

In view of the many classes of groups and other mathematical objects that $p$-local compact groups model, it is no surprise that this theory is not very well understood yet. Nevertheless several structural results have already been conjectured since its introduction in \cite{BLO3}, of which two are of special relevance in this work. The first conjecture asks about the existence of a version of the Stable Elements Theorem by H. Cartan and S. Eilenberg \cite[Theorem 10.1]{CE} for $p$-local compact groups, while the second conjecture suggests the extension of the classic result of W. Dwyer and A. Zabrodsky \cite[Theorem 1.1]{DZ} on mapping spaces to $p$-local compact groups.

A reader familiar with the theory of $p$-local finite groups will know that the proofs of these results in the finite case, as well as many others, depend ultimately on a certain biset, which is at the core this theory. And it is there where the main obstacle to study $p$-local compact groups resides, since there is no such biset in the context of $p$-local compact groups.

An alternative to the biset is necessary, and evidence that such alternative may exist is found in a result by E. M. Friedlander and G. Mislin, \cite{FM1, FM2}, that relates compact Lie groups to towers of finite groups. In this paper we prove a version of the result of Friedlander and Mislin for $p$-local compact groups.

Many statements for $p$-local finite groups extend to $p$-local compact groups as a consequence of this approximation result that replaces the use of the biset by arguments on telescopes of classifying spaces of $p$-local finite groups. In fact, in view of the results deduced in this work, the existence of finite approximations in all cases seems to be in itself one of the most essential properties of $p$-local compact groups.

In order to state our main results, let us be more precise. Let $\Z/p^{\infty}$ be the union of all the cyclic groups $\Z/p^n$, with $n \geq 1$, via the obvious inclusions. A \emph{discrete $p$-toral group} is a group $S$ that contains a normal subgroup $T$ isomorphic to a finite direct product of copies of $\Z/p^{\infty}$, and such that $T$ has finite index in $S$.

A \emph{saturated fusion system} over $S$ is a category $\FF$ whose object set is the collection of subgroups of $S$, and all of whose morphisms are actual group monomorphisms subject to a certain set of conditions. A \emph{centric linking system} associated to $\FF$ is a category $\LL$ whose object set is a certain subset of subgroups of $S$, and whose morphisms are again subject to another set of conditions. Roughly speaking, a (saturated) fusion system models the conjugation action of a group over a fixed Sylow $p$-subgroup (without the assumption that such group exists), while a centric linking system supplies the necessary data to produce a ``classifying space'' associated to the corresponding fusion system.

A \emph{$p$-local compact group} is a triple $\g = \ploc$, where $S$ is a discrete $p$-toral group, $\FF$ is a saturated fusion system over $S$, and $\LL$ is a centric linking system associated to $\FF$. The \emph{classifying space} of $\g$ is the $p$-completed nerve of the category $\LL$, $B\g = |\LL|^{\wedge}_p$. The exact definitions are included in Section \ref{Background}. Here and for the rest of the paper, $(-)^{\wedge}_p$ denotes the $p$-completion functor developed by A. K. Bousfield and D. M. Kan in \cite{BK}.

By an \emph{approximation of $\g$ by $p$-local finite groups} we mean a family of $p$-local finite groups $\{\g_i = (S_i, \FF_i, \LL_i)\}_{i \geq 0}$, such that $\LL_i \subseteq \LL_{i+1} \subseteq \LL$ for all $i$, and such that the homotopy colimit of the spaces $B\g_i$ is equivalent to $B\g$ after $p$-completion (see Section \ref{Sfam} for the precise definition).

Unstable Adams operations for $p$-local compact groups were defined and constructed by F. Junod, R. Levi and A. Libman in \cite{JLL}. Given a $p$-local compact group $\g = \ploc$ and an unstable Adams operation $\Psi$ on $\g$, we construct an approximation of $\g$ by $p$-local finite groups by taking fixed points of $\g$ by $\Psi$. This process essentially mimics the construction of fixed points of an algebraic group by successive iterations of a given Frobenius automorphism, and in this sense relates to the work of Friendlander and Mislin \cite{FM1, FM2}.

\begin{Th}\label{thmA}

Every $p$-local compact group admits an approximation by $p$-local finite groups.

\end{Th}

With this, one can deduce the Stable Elements Theorem \cite[Theorem 5.8]{BLO2} for $p$-local compact groups, which computes the cohomology of the classifying spaces of a $p$-local compact group as a subalgebra of the cohomology of its Sylow $p$-subgroup.

\begin{Th}\label{thmB}

Let $\g = \ploc$ be a $p$-local compact group, and let $M$ be a finite $\Z_{(p)}$-module with trivial $S$-action. Then, the natural map
$$
H^{\ast}(B\g; M) \Right3{\cong} H^{\ast}(\FF;M) \defin \varprojlim_{\FF} H^{\ast}(-; M) \subseteq H^{\ast}(BS; M)
$$
is an isomorphism.

\end{Th}

A version of the spectral sequence of A. D\'iaz \cite{Diaz} for $p$-local compact groups is also deduced easily. The notation $(-)^{\FF}$ below denotes stable elements in $\FF$, see Section  \ref{Sstable} for an explicit definition.

\begin{Th}\label{thmD}

Let $\FF$ be a saturated fusion system over a discrete $p$-toral group $S$. Let also $R\leq S$ be a strongly $\FF$-closed subgroup of $S$, and let $M$ be a finite $\Z_{(p)}$-module with trivial $S$-action. Then there is a first quadrant cohomological spectral sequence with second page
$$
E_2^{n,m} = H^n(S/R; H^m(R;M))^{\FF}
$$
and converging to $H^{n+m}(B\g;M)$.

\end{Th}

The above spectral sequence has a rather interesting feature: it uses only a strongly closed subgroup of a fusion system, instead of a \emph{normal} fusion subsystem, and thus is much more versatile. There are many interesting corollaries in \cite[Sections 5 and 6]{Diaz}, which we suspect to admit generalizations to the compact case. We leave this out of the scope of this paper for the sake of brevity.

Another relevant consequence of Theorem \ref{thmA} is the following result describing the homotopy type of certain mapping spaces, which generalizes \cite[Theorem 6.3]{BLO2}, and which was also conjectured in \cite{BLO3}.

\begin{Th}\label{thmC}

Let $\g = \ploc$ be a $p$-local compact group, and let $P$ be a discrete $p$-toral group. Then, each component of the mapping space $\map(BP, B\g)$ has the homotopy type of the classifying space of a $p$-local compact group, which corresponds to the ``centralizer'' of the image of the corresponding homomorphism $P \to S$.

\end{Th}\

So far, Theorem \ref{thmA} provides the only known way to circumvent the absence of a biset in the compact case. As shown in work by several authors, see for example \cite{DL1}, \cite{DL2}, \cite{RS}, \cite{GeRe} and \cite{GeReYa}, this biset is closely related to the Burnside ring of the corresponding $p$-local finite group. In this sense, it would be interesting to analyze the Burnside ring of a $p$-local compact group $\g = \ploc$, and compare it with the Burnside rings of a given approximation of $\g$ by $p$-local finite groups. Representations of $p$-local compact groups are another aspect of the theory that may benefit from Theorem \ref{thmA}. Very little is known in the compact case, see \cite{CanCas}, in comparison with what is known in the finite case, see \cite{CL}. Of particular interest for the author would be a generalization to $p$-local compact groups of the Peter-Weyl Theorem on existence of faithful complex representations for $p$-compact groups, see \cite[Theorem 1.6]{AGMV} (for $p$ odd) and \cite[Remark 7.3]{Andersen-Grodal} (for $p = 2$).


\emph{Organization of the paper.} In Section \ref{Background} we review the basic definitions and results about $p$-local compact groups that we need in later sections. It includes two subsections. The first one, \ref{Ssisotyp}, introduces isotypical equivalences, and reviews the main properties of unstable Adams operations. The second subsection, \ref{Squotient}, presents the constructions of normalizers, centralizers, and quotients in a $p$-local compact groups. This section does not contain any original result, but it is necessary in order to make this work self-contained.

Section \ref{Sbig} describes a general procedure to add new objects to a given transporter system. When applied to a centric linking system, this process generates a \emph{quasicentric} linking system. These enlarged transporter systems play a crucial role in the proof of Theorem \ref{thmA} in the next section. The construction introduced in this section is not as general as that in \cite{BCGLO1} for finite centric linking systems, but it turns out to be much simpler.

In Section \ref{Sfam} we construct finite approximations of $p$-local compact groups. This section is rather long, and thus we opted for splitting it in shorter subsections. Subsection \ref{Sprelim} includes the definition of an approximation of a $p$-local compact group by $p$-local finite groups and its main properties. In preparation for the following subsection, this part also contains the definition of the fixed points of a centric linking system by an unstable Adams operation. Subsection \ref{Ssapprox} is the core of this section, it essentially includes the proof of Theorem \ref{thmA}. The last part of this section, Subsection \ref{Ssexample}, contains an example of the construction of approximations of a $p$-local compact group by $p$-local finite groups.

Section \ref{Sstable} contains the proofs of Theorem \ref{thmB} and \ref{thmD}, while Section \ref{Smap} contains the proof of Theorem \ref{thmC}. The proofs for these results follow closely the corresponding proofs for their finite counterparts. We include a table of contents for the reader's convenience.

\tableofcontents


\textbf{Acknowledgements.} The author is specially grateful to C. Broto, R. Levi and B. Oliver for numerous inspirational conversations and their support. The author would like to thank also N. Castellana, A. Libman and R. Stancu for their help in various stages of this work, and A. Chermak and R. Molinier for their support during the final stages of this work. Finally, the author is truly indebted to the referee, without whose assistance and thorough job the reader's experience with this paper would be greatly diminished.

\textbf{Acknowledgements of financial support.} The author was partially supported by FEDER/MEC grants BES-2005-11029 and MTM2010-20692, and by SGR2009-1092.


\section{Background on \texorpdfstring{$p$}{p}-local compact groups}\label{Background}

Let $p$ a prime, to remain fixed for the rest of the paper unless otherwise stated. In this section we review all the definitions and results about $p$-local compact groups that we need in this paper. The main references for this section are the works of C. Broto, R. Levi and B. Oliver \cite{BLO2, BLO3, BLO6}.

Roughly speaking, $p$-local compact groups are abstractions of the fusion data obtained from finite and compact Lie groups. This idea already implies the existence of some sort of \emph{Sylow $p$-subgroup}. In the finite case this role was played by finite $p$-groups, while in this more general setup we use discrete $p$-toral groups. Let $\prufferp = \bigcup_{n \geq 1} \Z/p^n$ be the union of the cyclic $p$-groups $\Z/p^n$ under the obvious inclusions.

\begin{defi}\label{defidiscptor}

A \emph{discrete $p$-toral group} $P$ is a group that contains a normal subgroup $P_0$ of finite index and which is isomorphic to $(\prufferp)^{\times r}$ for some finite $r \geq 0$.

\end{defi}

In other words, a discrete $p$-toral group is a group $P$ fitting in an exact sequence
$$
\{1\} \to P_0 \Right4{} P \Right4{} \pi \to \{1\},
$$
where $\pi$ is a finite $p$-group and $P_0 \cong (\prufferp)^{\times r}$. The \emph{rank} of $P$, denoted by $\mathrm{rk}(P)$, is $r$, and the \emph{order} of $P$ is then defined as the pair $|P| \defin (\mathrm{rk}(P), |\pi(P)|)$, considered as an element of $\N^2$. This way we can compare the order of two discrete $p$-toral groups, by writing $|Q|\leq |P|$ if either $\mathrm{rk}(Q) < \mathrm{rk}(P)$, or $\mathrm{rk}(Q) = \mathrm{rk}(P)$ and $|\pi(Q)|\leq |\pi(P)|$. Given a discrete $p$-toral group $S$ and subgroups $P, Q\leq S$, define
$$
\Hom_S(P, Q) = \{f = c_x \in \Hom(P, Q) \,\, | \,\, \exists x \in S \mbox{ such that } x P x^{-1}\leq Q\}.
$$

\begin{defi}

Given a discrete $p$-toral group $S$, a \emph{fusion system} over $S$ is a category $\FF$ with $\Ob(\FF) = \{P\leq S\}$, and whose morphisms are actual homomorphisms satisfying the following:
\begin{enumerate}[(i)]

\item $\Hom_S(P,Q) \subseteq \Hom_{\FF}(P,Q) \subseteq \Inj(P,Q)$  for all $P, Q \in \Ob(\FF)$; and

\item every morphism in $\FF$ is the composition of an isomorphism in $\FF$, followed by an inclusion.

\end{enumerate}
The \emph{rank of $\FF$} is the rank of $S$.

\end{defi}

The following notation will be used tacitly throughout the rest of the paper. Let $\FF$ be a fusion system over $S$, and let $P, Q, X\leq S$. As objects in $\FF$, we say that $P$ and $Q$ are \emph{$\FF$-conjugate} if they are isomorphic as objects in $\FF$. The notation $P^X$ and $P^{\FF}$, respectively, stands for the $X$-conjugacy and $\FF$-conjugacy classes of $P$. Note also that $\Aut_{\FF}(P)$ is a group, by definition of fusion system, and that $\Inn(P)\leq \Aut_{\FF}(P)$. Thus, it is reasonable to define
$$
\Out_S(P) \defin \Aut_S(P)/\Inn(P) \qquad \mbox{and} \qquad \Out_{\FF}(P) \defin \Aut_{\FF}(P)/\Inn(P).
$$
Finally, we say that $P$ is \emph{fully $\FF$-centralized}, respectively \emph{fully $\FF$-normalized}, if $|C_S(P)| \geq |C_S(Q)|$ for all $Q \in P^{\FF}$, respectively  if $|N_S(P)| \geq |N_S(Q)|$ for all $Q \in P^{\FF}$.

\begin{defi}\label{defisat}

Let $S$ be a discrete $p$-toral group, and let $\FF$ be a fusion system over $S$. We say that $\FF$ is a \emph{saturated fusion system} if the following conditions are satisfied.
\begin{enumerate}[(I)]

\item If $P\leq S$ is a fully $\FF$-normalized subgroup, then it is also fully $\FF$-centralized. Moreover, in this case $\Out_{\FF}(P)$ is a finite group, and $\Out_S(P) \in \Syl_p(\Out_{\FF}(P))$.

\item Suppose $P\leq S$ and $f \in \Hom_{\FF}(P,S)$ are such that $f(P)$ is fully $\FF$-centralized, and set
$$
N_f = \{g \in N_S(P) \,\, | \,\, f \circ c_g \circ f^{-1} \in \Aut_S(f(P))\}.
$$
Then, there exists $\4{f} \in \Hom_{\FF}(N_f, S)$ such that $\4{f}|_P = f$.

\item Let $P_1\leq P_2\leq P_3\leq \ldots$ be a sequence of subgroups of $S$, and set $P = \bigcup_{n = 1}^{\infty} P_n$. If $f \in \Hom(P,S)$ is a homomorphism such that $f|_{P_n} \in \Hom_{\FF}(P_n,S)$ for all $n$, then $f \in \Hom_{\FF}(P,S)$.

\end{enumerate}

\end{defi}

We also recall the definition of centric and radical subgroups, which are crucial concepts in the $p$-local group theory.

\begin{defi}

Let $\FF$ be a saturated fusion system over a discrete $p$-toral group $S$.
\begin{itemize}

\item A subgroup $P\leq S$ is \emph{$\FF$-centric} if $C_S(Q) = Z(Q)$ for all $Q \in P^{\FF}$.

\item A subgroup $P\leq S$ is \emph{$\FF$-radical} if $\Out_{\FF}(P)$ contains no nontrivial normal $p$-subgroup.

\end{itemize}

\end{defi}

Given a saturated fusion system $\FF$ over a discrete $p$-toral group $S$, we denote by $\FF^c$ and $\FF^r$ the full subcategories of $\FF$ with object sets the collections of $\FF$-centric and $\FF$-radical subgroups, respectively. We also write $\FF^{cr} \subseteq \FF$ for the full subcategory of $\FF$-centric $\FF$-radical subgroups.

Proving that a given fusion system is saturated is a rather difficult task, even when the fusion system is finite, but there are some techniques that may be helpful. One of these techniques, which we will use in later sections, is \cite[Theorem A]{BCGLO1}, restated as Theorem \ref{5A} below.

\begin{defi}

Let $\FF$ be a fusion system over a finite $p$-group $S$, and let $\hh \subseteq \Ob(\FF)$ be a subset of objects.
\begin{itemize}

\item $\FF$ is \emph{$\hh$-generated} if every morphism in $\FF$ can be described as a composite of restrictions of morphisms in $\FF$ between subgroups in $\hh$.

\item $\FF$ is \emph{$\hh$-saturated} if the saturation axioms hold for all subgroups in the set $\hh$.

\end{itemize}

\end{defi}

\begin{thm}\label{5A}

Let $\FF$ be a fusion system over a finite $p$-group $S$, and let $\hh$ be a subset of objects of $\FF$ closed under $\FF$-conjugacy and such that $\FF$ is $\hh$-generated and $\hh$-saturated. Suppose further that, for each $\FF$-centric subgroup $P \notin \hh$, $P$ is $\FF$-conjugate to some $Q$ such that
$$
\Out_S(Q) \cap O_p(\Out_{\FF}(Q)) \neq \{1\}.
$$
Then $\FF$ is saturated.

\end{thm}

The following result will be useful in later sections when checking the condition displayed in the previous theorem. We state it in full generality since in fact it will apply in different situations throughout this paper.

\begin{lmm}\label{Kpgp}

Let $\FF$ be a fusion system over a discrete $p$-toral group $S$. Let also $P\leq S$ be a subgroup, and let $P_0 \lhd P$ be a normal subgroup such that $f|_{P_0} \in \Aut_{\FF}(P_0)$ for all $f \in \Aut_{\FF}(P)$. Set
$$
K_P \defin \Ker(\Aut_{\FF}(P) \Right2{} \Aut_{\FF}(P_0) \times \Aut(P/P_0)).
$$
Then, $K_P\leq O_p(\Aut_{\FF}(P))$.

\end{lmm}

\begin{proof}

By definition $K_P$ is normal in $\Aut_{\FF}(P)$, and thus we only have to show that $K_P$ is a discrete $p$-toral subgroup of $\Aut_{\FF}(P)$ (that is, every element of $K_P$ has order a power of $p$). This in turn follows by \cite[Theorem 3.2]{Gor}: although the result in \cite{Gor} is stated for finite groups, the arguments in its proof apply here without modification, since $\Aut_{\FF}(P)$ is a locally finite group.
\end{proof}

The concept of transporter system associated to a fusion system was first introduced in \cite{OV} for fusion systems over finite $p$-groups, and then extended to discrete $p$-toral groups in \cite{BLO6}, with centric linking systems as a particular case. We refer the reader to the aforementioned sources for further details.

Let $G$ be a group and let $\hh$ be a set of subgroups of $G$ that is closed by overgroups, i.e. if $H \in \hh$ and $K \geq H$, then $K \in \hh$, and closed by conjugation in $G$, i.e. if $H \in \hh$ and $g \in G$, then $gHg^{-1} \in \hh$. The transporter category of $G$ with respect to $\hh$ is the category $\TT_{\hh}(G)$ whose object set is $\hh$, and with morphism sets
$$
\Mor_{\TT_{\hh}(G)}(P,Q) = \{x \in G \,\, | \,\, x \cdot P \cdot x^{-1}\leq Q\}
$$
for each $P,Q \in \hh$.

\begin{defi}\label{defitransporter}

Let $S$ be a discrete $p$-toral group, and let $\FF$ be a fusion system over $S$. A \emph{transporter system} associated to $\FF$ is a nonempty category $\TT$ whose object set $\Ob(\TT)$ is a subset of $\Ob(\FF)$ that is closed by overgroups and conjugation in $\FF$, together with functors
$$
\TT_{\Ob(\TT)}(S) \Right4{\varepsilon} \TT \qquad \mbox{and} \qquad \TT \Right4{\rho} \FF
$$
satisfying the following conditions.
\begin{itemize}

\item[(A1)] The functor $\varepsilon$ is the identity on objects and an inclusion on morphism sets, and the functor $\rho$ is the inclusion on objects and a surjection on morphism sets.

\item[(A2)] For each $P, Q \in \Ob(\TT)$, the set $\Mor_{\TT}(P,Q)$ has a free action of
$$
E(P) \defin \Ker \big[\rho_P \colon \Aut_{\TT}(P) \Right2{} \Aut_{\FF}(P) \big]
$$
by right composition, and $\rho_{P,Q}$ is the orbit map of this action. Also, $E(Q)$ acts freely on $\Mor_{\TT}(P,Q)$ by left composition.

\item[(B)] Let $P,Q \in \Ob(\TT)$. Then, the map $\varepsilon_{P,Q} \colon N_S(P,Q) \to \Mor_{\TT}(P,Q)$ is injective, and
$$
(\rho_{P,Q} \circ \varepsilon_{P,Q})(g) = c_g \in \Hom_{\FF}(P,Q)
$$
for all $g \in \Mor_{\TT_{\Ob(\TT)}(S)}(P, Q) = N_S(P, Q)$.

\item[(C)] For all $P, Q \in \Ob(\TT)$, for all $\varphi \in \Mor_{\TT}(P,Q)$, and for all $g \in P$, the following is a commutative diagram in $\TT$.
$$
\xymatrix{
P \ar[r]^{\varphi} \ar[d]_{\varepsilon_P(g)} & Q \ar[d]^{\varepsilon_Q(\rho(\varphi)(g))} \\
P \ar[r]_{\varphi} & Q
}
$$

\item[(I)] Each isomorphism class of objects in $\Ob(\TT)$ contains an element $P$ such that
$$
\varepsilon_P(N_S(P)) \in \Syl_p(\Aut_{\TT}(P));
$$
or, in other words, such that $\varepsilon(N_S(P))$ has finite index prime to $p$ in $\Aut_{\TT}(P)$.

\item[(II)] Let $P, Q \in \Ob(\TT)$ be isomorphic objects, and let $\varphi \in \Iso_{\TT}(P,Q)$. Let also $\4{P}\leq N_S(P)$ and $\4{Q}\leq N_S(Q)$ be such that $\varphi \circ \varepsilon_P(\4{P}) \circ \varphi^{-1}\leq \varepsilon_Q(\4{Q})$. Then there is some morphism $\4{\varphi} \in \Mor_{\TT}(\4{P}, \4{Q})$ such that
$$
\4{\varphi} \circ \varepsilon_{P, \4{P}}(1) = \varepsilon_{Q, \4{Q}}(1) \circ \varphi.
$$

\item[(III)] Let $P_1\leq P_2\leq P_3\leq \ldots$ be a sequence in $\Ob(\TT)$, and let $\varphi_n \in \Mor_{\TT}(P_n,S)$ be such that $\varphi_n = \varphi_{n+1} \circ \varepsilon_{P_n, P_{n+1}}(1)$ for all $n \geq 1$. Then, upon setting $P = \bigcup_{n \geq 1} P_n$, there is a morphism $\varphi \in \Mor_{\TT}(P,S)$ such that $\varphi_n = \varphi \circ \varepsilon_{P_n,P}(1)$ for all $n \geq 1$.

\end{itemize}
The \emph{rank of $\TT$} is the rank of $S$. A \emph{centric linking system} associated to a saturated fusion system $\FF$ is a transporter system $\LL$ such that $\Ob(\LL)$ is the collection of all $\FF$-centric subgroups of $S$ and $E(P) = \varepsilon(Z(P))$ for all $P \in \Ob(\LL)$.

\end{defi}

\begin{rmk}\label{rmktransp}

The above definition of centric linking system is taken from \cite{BLO6}, and it is seen in \cite[Corollary A.5]{BLO6} to coincide with the original \cite[Definition 4.1]{BLO3}. Notice also that axiom (I) above differs from the corresponding axiom for the finite case, see \cite[Definition 3.1]{OV}, in that condition (I) above seems to be more restrictive than the corresponding condition in \cite{OV}:
\begin{itemize}

\item[(I')] $\varepsilon_{S,S}(S) \in \Syl_p(\Aut_{\TT}(S))$.

\end{itemize}
However, \cite[Proposition 3.4]{OV} implies that both definitions, \cite[Definition 3.1]{OV} and the above, agree in the finite case.

\end{rmk}

\begin{lmm}\label{epimono}

In a transporter system, all morphisms are monomorphisms and epimorphisms in the categorical sense.

\end{lmm}

\begin{proof}

This is \cite[Proposition A.2 (d)]{BLO6}.
\end{proof}

\begin{defi}

A \emph{$p$-local compact group} is a triple $\g = \ploc$ formed by a discrete $p$-toral group $S$, a saturated fusion system $\FF$ over $S$, and a centric linking system $\LL$ associated to $\FF$. The \emph{classifying space} of a $p$-local compact group $\g$ is the $p$-completed nerve of $\LL$, denoted by $B\g = |\LL|^{\wedge}_p$. The \emph{rank of $\g$} is the rank of $S$.

\end{defi}

Generalizing work of \cite{Chermak} and \cite{Oliver}, it is proved in \cite{Levi-Libman} that every saturated fusion system over a discrete $p$-toral group has an associated centric linking system which is unique up to isomorphism. Thus, from now on we speak of \emph{the} associated centric linking system for a given saturated fusion system. 

Finally, we recall the ``bullet construction'' on a $p$-local compact group.

\begin{defi}\label{defibullet}

Let $\FF$ be a saturated fusion system over a discrete $p$-toral group $S$. Let also $T\leq S$ be the maximal torus, and let $W = \Aut_{\FF}(T)$. Set the following

\begin{enumerate}[(i)]

\item The exponent of $S/T$ is $e = \exp(S/T) = \min\{k \in \N \, | \, x^{p^k} \in T \mbox{ for all } x \in S\}$.

\item For each $P\leq T$, let $I(P) = \{t \in T \, | \, \omega(t) = t \mbox{ for all } \omega \in W \mbox{ such that } \omega|_P = \Id_P\}$, and let $I(P)_0$ denote its maximal torus.

\item For each $P\leq S$, set $P^{[e]} = \{x^{p^e} \, | \, x \in P\}\leq T$, and set
$$
P^{\bullet} = P \cdot I(P^{[e]})_0 = \{xt \, | \, x \in P, \, t \in I(P^{[e]})_0\}.
$$

\item Let $\FF^{\bullet}$ be the full subcategory of $\FF$ with object set $\Ob(\FF^{\bullet}) = \{P^{\bullet} \, | \, P\leq S\}$.

\end{enumerate}

\end{defi}

The following summarizes the main properties of the ``bullet construction''.

\begin{prop}\label{3.2BLO3}

Let $\g = \ploc$ be a $p$-local compact group. Then, for each $P, Q \in \Ob(\FF)$ and each $f \in \Hom_{\FF}(P,Q)$ there is a unique $f^{\bullet} \in \Hom_{\FF}(P^{\bullet}, Q^{\bullet})$ whose restriction to $P$ is $f$. This way, the ``bullet construction'' makes $P \mapsto P^{\bullet}$ into a functor $\bullf \colon \FF \to \FF$ that satisfies the following properties.
\begin{enumerate}[(i)]

\item The set $\Ob(\FF^{\bullet}) = \{P^{\bullet} \,  | \, P\leq S\}$ contains finitely many $S$-conjugacy classes of subgroups of $S$.

\item For all $P\leq S$, $(P^{\bullet})^{\bullet} = P^{\bullet}$.

\item If $P\leq Q\leq S$, then $P^{\bullet}\leq Q^{\bullet}$.

\item For all $P, Q\leq S$, $N_S(P, Q) \subseteq N_S(P^{\bullet}, Q^{\bullet})$.

\item For all $P\leq S$, $C_S(P) = C_S(P^{\bullet})$.

\item The functor $\bullf$ is a left adjoint to the inclusion of $\FF^{\bullet}$ as a full subcategory of $\FF$.

\item All $\FF$-centric $\FF$-radical subgroups of $S$ are in $\FF^{\bullet}$. In particular, there are only finitely many $\FF$-conjugacy classes of such subgroups.

\end{enumerate}
Moreover, if we denote by $\LL^{\bullet} \subseteq \LL$ the full subcategory with $\Ob(\LL^{\bullet}) = \{P^{\bullet} \, | \, P \in \Ob(\LL)\}$, then there is a unique functor $\bullf \colon \LL \to \LL^{\bullet}$ such that the following holds.
\begin{enumerate}[(a)]

\item $\bullf \circ \rho = \rho \circ \bullf \colon \LL \to \FF$.

\item For all $P, Q \in \Ob(\LL)$ and all $\varphi \in \Mor_{\LL}(P,Q)$, we have $\varepsilon_{Q, Q^{\bullet}}(1) \circ \varphi = \varphi^{\bullet} \circ \varepsilon_{P, P^{\bullet}}(1)$.

\item For all $P, Q \in \Ob(\LL)$ and all $g \in N_S(P,Q)$, we have $\varepsilon_{P,Q}(g)^{\bullet} = \varepsilon_{P^{\bullet}, Q^{\bullet}}(g)$.

\item The functor $(-)^{\bullet} \colon \LL \to \LL$ is left adjoint to the inclusion of $\LL^{\bullet}$ as a full subcategory of $\LL$. In particular, the inclusion $\LL^{\bullet} \subseteq \LL$ induces an equivalence $|\LL^{\bullet}|\simeq |\LL|$.

\end{enumerate}

\end{prop}

\begin{proof}

The first part of the statement corresponds to \cite[Proposition 3.3]{BLO3}. Parts (i), (ii) and (iii)  correspond to \cite[Lemma 3.2 (a), (b) and (c)]{BLO3} respectively. Part (iv) is an easy variation of \cite[Lemma 3.2 (b)]{BLO3} (details are left to the reader). For part (v), let $P\leq S$. Since $P\leq P^{\bullet}$, we have $C_S(P) \geq C_S(P^{\bullet})$. Let $x \in C_S(P)$. By (iv), we have $x \in N_S(P^{\bullet})$. Since $c_x = \Id \in \Aut_{\FF}(P)$ extends uniquely to $c_x = \Id \in \Aut_{\FF}(P^{\bullet})$, it follows that $x \in C_S(P^{\bullet})$. Part (vi) corresponds to \cite[Corollary 3.4]{BLO3}, and part (vii) corresponds to \cite[Corollary 3.5]{BLO3}. The last part of the statement, including parts (a), (b) and (c) corresponds to \cite[Proposition 1.12]{JLL}. Part (d) corresponds to \cite[Proposition 4.5 (a)]{BLO3}.
\end{proof}


\subsection{Isotypical equivalences and unstable Adams operations}\label{Ssisotyp}

In this subsection we review the concept of isotypical equivalence, with particular interest on the unstable Adams operations for $p$-local compact groups originally introduced in \cite{JLL}.

\begin{defi}

Let $(\TT, \varepsilon, \rho)$ be a transporter system associated to a fusion system $\FF$. An automorphism $\Psi \colon \TT \to \TT$ is \emph{isotypical} if $\Psi(\varepsilon_P(P)) = \varepsilon_{\Psi(P)}(\Psi(P))$ for each $P \in \Ob(\TT)$.

\end{defi}

We denote by $\Aut_{\typ}^I(\TT)$ the group of isotypical automorphisms $\Psi$ of $\TT$ which in addition satisfy $\Psi(\varepsilon_{P,Q}(1)) = \varepsilon_{\Psi(P), \Psi(Q)}(1)$ whenever $P\leq Q$. Notice that if $\Psi \in \Aut_{\typ}^{I}(\TT)$, then $\Psi$ induces an automorphism of $S$ by restricting to the object $S \in \Ob(\TT)$. By abuse of notation, we will denote the induced automorphism by $\Psi \in \Aut(S)$.

Next we review the concept of unstable Adams operations for $p$-local compact groups. Our definition corresponds to the definition of \emph{normal Adams operation} in \cite[Definition 3.3]{JLL}, conveniently adapted to our notation. By $(\padic)^{\times}$ we denote the subgroup of multiplicative units in the ring of $p$-adic integers $\padic$.

\begin{defi}\label{uAo}

Let $\g = \ploc$ be a $p$-local compact group. An \emph{unstable Adams operation of degree $\zeta \in (\padic)^{\times}$} on $\g$ is an isotypical equivalence $\Psi \in \Aut_{\typ}^{I}(\LL)$ such that the induced automorphism $\Psi \in \Aut(S)$ satisfies
\begin{enumerate}[(i)]

\item the restriction of $\Psi$ to the maximal torus $T\leq S$ is the $\zeta$-power automorphism; and

\item $\Psi$ induces the identity on $S/T$.

\end{enumerate}
An unstable Adams operation is \emph{fine} if its degree is $\zeta \neq 1$, with $\zeta$ congruent to $1$ modulo $p$.

\end{defi}

As proved in \cite[Theorem 4.1]{JLL}, unstable Adams operations exist for all $p$-local compact groups, and in particular this applies to the existence of fine unstable Adams operations.

\begin{thm}

Let $\g = \ploc$ be a $p$-local compact group. Then, for some large enough $m \in \N$, $\g$ has unstable Adams operations of degree $\zeta$, for each $\zeta \in 1 + p^m\padic$.

\end{thm}

\begin{rmk}\label{uAo1}

Roughly speaking, the construction of unstable Adams operations in \cite{JLL} is done by defining $\Psi$ to fix enough objects and morphisms in $\LL$. More specifically, $\Psi$ fixes
\begin{enumerate}[(a)]

\item a set $\hh$ of representatives of the $S$-conjugacy classes $\Ob(\LL^{\bullet})$; and

\item for each $P \in \hh$, a set of representatives $\mm_P \subseteq \Aut_{\LL}(P)$ of the classes in $\Aut_{\LL}(P)/P \cong \Out_{\FF}(P)$.

\end{enumerate}
This properties will be crucial in our constructions in Section \ref{Sfam}.

\end{rmk}

Let $S$ be a discrete $p$-toral group, let $\FF$ be a fusion system over $S$ (not necessarily saturated), and let $\TT$ be a transporter system associated to $\FF$. Let also $\Psi \in \Aut_{\typ}^{I}(\TT)$ be an isotypical automorphism. Set also
\begin{equation}\label{fixS}
C_S(\Psi) = \{g \in S \, | \, \Psi(\varepsilon_S(g)) = \varepsilon_S(g)\}\leq S,
\end{equation}
the subgroup of fixed points of $S$ by $\Psi$. The following result is the main tool in detecting objects and morphisms in $\TT$ that are invariant by $\Psi$.

\begin{lmm}\label{invar1}

The following holds.
\begin{enumerate}[(i)]

\item Let $P\leq C_S(\Psi)$, and let $Q \in P^S$. Then, $Q\leq C_S(\Psi)$ if and only if, for some $x \in N_S(Q,P)$,
$$
x^{-1} \cdot \Psi(x) \in C_S(Q).
$$

\item Let $P, P' ,Q, Q'\leq C_S(\Psi)$ be such that $P' \in P^S$ and $Q' \in Q^S$, and suppose $P, P', Q, Q' \in \Ob(\LL)$. Let also $x \in N_S(P',P)$ and $y \in N_S(Q', Q)$, and let $\varphi \in \Mor_{\TT}(P,Q)$ be such that $\Psi(\varphi) = \varphi$. Set $\varphi' = \varepsilon(y^{-1}) \circ \varphi \circ \varepsilon(x) \in \Mor_{\wL}(P', Q')$. Then, $\Psi(\varphi') = \varphi'$ if and only if
$$
\varepsilon(y^{-1} \cdot \Psi(y)) \circ \varphi' = \varphi' \circ \varepsilon(x^{-1} \cdot \Psi(x)).
$$

\end{enumerate}

\end{lmm}

\begin{proof}

For part (i), let $P\leq C_S(\Psi)$, and let $Q \in P^S$. Let also $g \in Q$ and $x \in N_S(Q,P)$, and set $h = x \cdot g \cdot x^{-1} \in P$. Since $P\leq C_S(\Psi)$, we get
$$
x \cdot g \cdot x^{-1} = h = \Psi(h) = \Psi(x \cdot g \cdot x^{-1}) = \Psi(x) \cdot \Psi(g) \cdot \Psi(x)^{-1}.
$$
Thus, if $Q\leq C_S(\Psi)$ then clearly $x^{-1} \cdot \Psi(x) \in C_S(Q)$, and conversely if $x^{-1} \cdot \Psi(x) \in C_S(Q)$ then $h \in C_S(\Psi)$. Since the argument works for any $g \in Q$ and any $x \in N_S(Q, P)$, part (i) follows.

For part (ii), let $P, P' ,Q, Q'\leq C_S(\Psi)$, with $P' \in P^S$ and $Q' \in Q^S$. Let also $x \in N_S(P',P)$ and $y \in N_S(Q', Q)$, and let $\varphi \in \Mor_{\TT}(P,Q)$ be such that $\Psi(\varphi) = \varphi$, with $\varphi' = \varepsilon(y^{-1}) \circ \varphi \circ \varepsilon(x) \in \Mor_{\wL}(P', Q')$. We have
$$
\varepsilon(y) \circ \varphi' \circ \varepsilon(x^{-1}) = \varphi = \Psi(\varphi) = \Psi(\varepsilon(y) \circ \varphi' \circ \varepsilon(x^{-1})) = \varepsilon(\Psi(y)) \circ \Psi(\varphi') \circ \varepsilon(\Psi(x)^{-1}),
$$
and (ii) follows easily.
\end{proof}


\subsection{Normalizers, centralizers, and related constructions}\label{Squotient}

In this subsection we review the construction of the centralizer and normalizer $p$-local compact subgroups for a given $p$-local compact group. The main references here are \cite[Appendix A]{BLO2} and \cite[Section 2]{BLO6}. For the rest of this subsection, fix a $p$-local compact group $\g = \ploc$, a subgroup $A\leq S$, and a subgroup $K\leq \Aut(A)$, and define the following:
\begin{itemize}

\item $\Aut_{\FF}^K(A) = K \cap \Aut_{\FF}(A)$;

\item $\Aut_S^K(A) = K \cap \Aut_S(A)$; and

\item $N_S^K(A) = \{x \in N_S(A) \, | \, c_x \in K\}$.

\end{itemize}
The subgroup $A$ is \emph{fully $K$-normalized in $\FF$} if we have $|N_S^K(A)| \geq |N_S^{^{f}K}(f(A))|$ for each $f \in \Hom_{\FF}(A, S)$, where $^{f}K = \{f \gamma f^{-1} \, | \, \gamma \in K\}\leq \Aut(f(A))$.

\begin{defi}\label{definorm}

The \emph{$K$-normalizer fusion system of $A$ in $\FF$}, is the fusion system $N_{\FF}^K(A)$ over $N_S^K(A)$ with morphism sets
$$
\begin{aligned}
\Hom_{N_{\FF}^K(A)}&(P,Q) = \\
 & = \{f \in \Hom_{\FF}(P,Q) \,\, | \,\, \exists \4{f} \in \Hom_{\FF}(PA, QA) \mbox{ with } \4{f}|_P = f \mbox{ and } \4{f}|_A \in K\}
\end{aligned}
$$
for each $P, Q\leq N_S^K(A)$.

\end{defi}

By \cite[Theorem 2.3]{BLO6} we know that $N_{\FF}^K(A)$ is a saturated fusion system whenever $A$ is fully $K$-normalized in $\FF$. For this reason, for the rest of this subsection we assume that $A$ satisfies this property.

\begin{lmm}\label{centricNFKA}

If $P\leq N_S^K(A)$ is $N_{\FF}^K(A)$-centric, then $P \cdot A$ is $\FF$-centric.

\end{lmm}

\begin{proof}

Let $P\leq N_S^K(A)$ be $N_{\FF}^K(A)$-centric. We have to check that, for each $\gamma \in \Hom_{\FF}(P \cdot A, S)$, there is an inclusion $C_S(\gamma(P \cdot A))\leq \gamma(P \cdot A)$. We can apply \cite[Proposition A.2]{BLO2}, since the proof in \cite{BLO2} works without modifications in the compact setup, and it follows that the subgroup $A$ is fully centralized in $\FF$, and there is some $f \in \Hom_{\FF}(N_S^{^{\gamma}K}(\gamma(A)) \cdot \gamma(A),S)$ such that $(f \circ \gamma)|_A \in K$. Thus $f \circ \gamma$ is a morphism in $N_{\FF}^K(A)$.

Note that $C_S(\gamma(P \cdot A))\leq C_S(\gamma(A))\leq N_S^{^{\gamma} K}(\gamma(A))$, and we have inclusions
$$
\begin{aligned}
f(C_S(\gamma(P \cdot A))) &\leq C_S((f \circ \gamma)(P \cdot A)) = C_S((f \circ \gamma)(P) \cdot A)\leq \\
 &\leq C_S((f \circ \gamma)(P)) \cap C_S(A)\leq C_S((f \circ \gamma)(P)) \cap N_S^K(A)\leq (f \circ \gamma)(P),
\end{aligned}
$$
where the last inequality holds since $P \in N_{\FF}^K(A)^c$. Thus,
$$
C_S(\gamma(P \cdot A))\leq \gamma(P)\leq \gamma(P) \cdot \gamma(A) = \gamma(P \cdot A),
$$
and this proves that $P \cdot A \in \FF^c$.
\end{proof}

In view of the above, we can now define $N_{\LL}^K(A)$ as the category with objects the set of $N_{\FF}^K(A)$-centric subgroups of $N_S^K(A)$ and with morphism sets
$$
\Mor_{N_{\LL}^K(A)}(P,Q) = \{\varphi \in \Mor_{\LL}(PA, QA) \,\, | \,\, \rho(\varphi)|_P \in \Hom_{N_{\FF}^K(A)}(P,Q) \mbox{ and } \rho(\varphi)|_A \in K\}.
$$
In general, $N_{\LL}^K(A)$ need not be a transporter system associated to $N_{\FF}^K(A)$, but there are two particular situations where this is indeed the case.

\begin{lmm}

If either $K = \{\Id\}$ or $K = \Aut(A)$, then the category $N_{\LL}^K(A)$ is a centric linking system associated to $N_{\FF}^K(A)$.

\end{lmm}

\begin{proof}

The case $K = \{\Id\}$ corresponds to \cite[Proposition 2.5]{BLO2} in the finite case, while the case $K = \Aut(A)$ corresponds to \cite[Lemma 6.2]{BLO2} for $p$-local finite groups. In both situations, the proof for $p$-local finite groups applies here without modification to show that $N_{\LL}^{\{\Id\}}(A)$ satisfies all the condition of a centric linking system, except perhaps axiom (III), which is easily checked.
\end{proof}

\begin{defi}\label{rmknorm}

Let $\g = \ploc$ be a $p$-local compact group, and let $A\leq S$.
\begin{enumerate}[(a)]

\item If $A$ is fully $\FF$-centralized, the \emph{centralizer $p$-local compact group of $A$ in $\g$} is the triple
$$
C_{\g}(A) = (C_S(A), C_{\FF}(A), C_{\LL}(A)) \defin (N_S^{\{\Id\}}(A), N_{\FF}^{\{\Id\}}(A), N_{\LL}^{\{\Id\}}(A)).
$$

\item If $A$ is fully $\FF$-normalized, the \emph{normalizer $p$-local compact group of $A$ in $\g$} is the triple
$$
N_{\g}(A) = (N_S(A), N_{\FF}(A), N_{\LL}(A)) \defin (N_S^{\Aut(A)}(A), N_{\FF}^{\Aut(A)}(A), N_{\LL}^{\Aut(A)}(A)).
$$

\end{enumerate}
A subgroup $A\leq S$ is called \emph{central in $\FF$} if $C_{\g}(A) = \g$. Similarly, $A\leq S$ is called \emph{normal in $\FF$} if $N_{\g}(A) = \g$. Clearly, if $A\leq S$ is central in $\FF$ then in particular it is normal in $\FF$.

\end{defi}

\begin{lmm}\label{central1}

Let $\g = \ploc$ be a $p$-local compact group, and let $P\leq S$. Then $P$ is fully $\FF$-centralized if and only if $P^{\bullet}$ is fully $\FF$-centralized. Furthermore, if this is the case then $C_{\g}(P) = C_{\g}(P^{\bullet})$.

\end{lmm}

\begin{proof}

Suppose first that $P^{\bullet}$ is fully $\FF$-centralized. By Proposition \ref{3.2BLO3} (v), we have $C_S(P) = C_S(P^{\bullet})$. If $Q \in P^{\FF}$, then $Q^{\bullet} \in (P^{\bullet})^{\FF}$ by Proposition \ref{3.2BLO3}, and we have
$$
|C_S(Q)| = |C_S(Q^{\bullet})|\leq |C_S(P^{\bullet})| = |C_S(P)|,
$$
which implies that $P$ is fully $\FF$-centralized.

Conversely, suppose that $P$ is fully $\FF$-centralized, and let $R \in (P^{\bullet})^{\FF}$ be fully $\FF$-centralized. Choose some $\gamma \in \Hom_{\FF}(P^{\bullet} C_S(P^{\bullet}), S)$ such that $\gamma(P^{\bullet}) = R$, and set $Q = \gamma(P)$, with $Q^{\bullet} = \gamma(P^{\bullet}) = R$. By Proposition \ref{3.2BLO3} (v), we have $C_S(P) = C_S(P^{\bullet})$, and thus
$$
\gamma(C_S(P^{\bullet})) = \gamma(C_S(P)) = C_S(Q) = C_S(Q^{\bullet}) = C_S(R),
$$
where the leftmost and rightmost equalities hold by Proposition \ref{3.2BLO3} (v), and the equality in the middle holds since $P$ is fully $\FF$-centralized. It follows that $P^{\bullet}$ is fully $\FF$-centralized.

To finish the proof, suppose that $P$ and $P^{\bullet}$ are fully $\FF$-centralized, and consider $C_{\g}(P)$ and $C_{\g}(P^{\bullet})$, which are $p$-local compact groups with Sylow $C_S(P) = C_S(P^{\bullet})$. By definition, it is enough to show that $C_{\FF}(P) = C_{\FF}(P^{\bullet})$. Notice that there is an obvious inclusion $C_{\FF}(P^{\bullet}) \subseteq C_{\FF}(P)$. Let $Q, R\leq C_S(P)$, and let $f \in \Hom_{C_{\FF}(P)}(Q,R)$. By definition of $C_{\FF}(P)$, there is some $\4{f} \in \Hom_{\FF}(QP, RP)$ such that $\4{f}|_Q = f$ and $\4{f}|_P = \Id$. Let $\gamma = \4{f}$, and consider $\gamma^{\bullet} \in \Hom_{\FF}((QP)^{\bullet}, (RP)^{\bullet})$. Then $\gamma^{\bullet}$ restricts to a morphism $\omega \in \Hom_{\FF}(QP^{\bullet}, RP^{\bullet})$. Furthermore, by definition of $\omega$, we have
$$
\omega|_Q = \gamma^{\bullet}|Q = \gamma|Q = f \qquad \mbox{and} \qquad \omega|_{P^{\bullet}} = \gamma^{\bullet}|_{P^{\bullet}} = (f|_P)^{\bullet} = \Id.
$$
Thus $f$ is a morphism in $C_{\FF}(P^{\bullet})$, and $C_{\FF}(P) = C_{\FF}(P^{\bullet})$.
\end{proof}

\begin{cor}\label{central2}

Let $\g = \ploc$ be a $p$-local compact group. Then, for each $P\leq S$ which is fully $\FF$-centralized, there is a sequence of finite subgroups $P_0\leq P_1\leq \ldots$ such that $P = \bigcup_{n \geq 0} P_n$ and such that the following conditions hold for all $n \geq 0$.
\begin{enumerate}[(i)]

\item $P_n$ is fully $\FF$-centralized and $P_n^{\bullet} = P^{\bullet}$.

\item $C_{\g}(P_n) = C_{\g}(P)$.

\end{enumerate}

\end{cor}

\begin{proof}

Let $P\leq S$. Since $S$ is locally finite, so is $P$, and we can find some sequence of finite subgroups $P_0\leq P_1\leq \ldots$ such that $P = \bigcup_{n \geq 0} P_n$. Furthermore, there is some $M \in \N$ such that $P_n^{\bullet} = P^{\bullet}$ for all $n \geq M$, and we may assume for simplicity that $M = 0$. By Lemma \ref{central1}, $P^{\bullet}$ is fully $\FF$-centralized, and then so is $P_n$, for all $n \geq 0$. Furthermore,
$$
C_{\FF}(P_n) = C_{\FF}(P_n^{\bullet}) = C_{\FF}(P^{\bullet}) = C_{\FF}(P),
$$
and this finishes the proof.
\end{proof}

To finish this section, we recall the construction of the quotient of a transporter system by a $p$-group. This quotient was already explored in \cite[Appendix A]{Gonza2}, and here we only recall the necessary definitions. Let $(\TT, \varepsilon, \rho)$ be a transporter system associated to a fusion system $\FF$, and let $A\leq S$ be a normal subgroup in $\FF$. If $P, Q\leq S$ are such that $A\leq P,Q$, then each morphism $f \in \Hom_{\FF}(P,Q)$ restricts to an automorphism of $A$, and hence it also induces a homomorphism $\ind(f) \colon P/A \to Q/A$. For a subgroup $P/A\leq S/A$, we will denote by $P\leq S$ the unique subgroup of $S$ that contains $A$ with image $P/A$ through the projection $S \to S/A$.

\begin{defi}\label{quotient1}

Let $A$ is a normal subgroup in $\FF$. The \emph{quotient} of $\TT$ by $A$ is the transporter system $(\TT/A, \3{\varepsilon}, \3{\rho})$ associated to the fusion system $\FF/A$, where
\begin{itemize}

\item $\FF/A$ is the fusion system over $S/A$ with morphism sets
$$
\begin{aligned}
\Hom_{\FF/A}&(P/A, Q/A) = \\
 & = \{\overline{f} \in \Hom(P/A,Q/A) \,\, | \,\, \exists f \in \Hom_{\FF}(P,Q) \mbox{ such that } \overline{f} = \ind(f)\}.
\end{aligned}
$$

\item $\TT/A$ is the category with object set $\{P/A\leq S/A \, | \, A\leq P \in \Ob(\TT)\}$ and morphism sets
$$
\Mor_{\TT/A}(P/A,Q/A) = \Mor_{\TT}(P,Q)/\varepsilon_P(A).
$$
The structural functors $\3{\varepsilon}$ and $\3{\rho}$ are induced, respectively, by the structural functors $\varepsilon$ and $\rho$ of $\TT$.

\end{itemize}

\end{defi}

\begin{rmk}\label{quotient21}

By \cite[Proposition A.2]{Gonza2}, $(\TT/A, \3{\varepsilon}, \3{\rho})$ is a transporter system associated to $\FF/A$. 

\end{rmk}

\begin{lmm}\label{quotient22}

Let $S$ be a discrete $p$-toral group, let $\FF$ be a saturated fusion system over $S$, and let $\TT$ be a transporter system associated to $\FF$, such that $\Ob(\TT)$ contains all the centric subgroups of $\FF$. Let also $A\leq S$ be normal in $\FF$. Then the following holds:
\begin{enumerate}[(i)]

\item the fusion system $\FF/A$ is saturated; and

\item the transporter system $\TT/A$ contains all the $\FF/A$-centric subgroups of $S/A$.

\end{enumerate}

\end{lmm}

\begin{proof}

Part (i) corresponds to \cite[Proposition A.3]{Gonza2}, and part (ii) is easily checked: let $P/A\leq S/A$ be $\FF/A$-centric, and let $Q/A$ be $\FF/A$-conjugate to $P/A$. If $P, Q\leq S$ denote the preimages of $P/A$ and $Q/A$ in $S$ respectively, then $Q$ is $\FF$-conjugate to $P$ by definition of $\FF/A$. Moreover,
$$
C_S(Q)A/A\leq C_{S/A}(Q/A)\leq Q/A,
$$
and thus $C_S(Q)\leq Q$ (since $A\leq Q$). It follows that $P$ is $\FF$-centric, and hence an object in $\TT$. This implies that $P/A \in \Ob(\TT/A)$.
\end{proof}


\section{Telescopic transporter systems}\label{Sbig}

In this section we describe a general procedure to add new objects to a given transporter system. The constructions in this section play a crucial role in the next section.

\begin{defi}

Let $S$ be a discrete $p$-toral group, let $\FF$ be a fusion system over $S$ (not necessarily saturated), and let $\TT$ be a transporter system associated to $\FF$. The transporter system $\TT$ is \emph{telescopic} if it satisfies the following condition.
\begin{itemize}

\item[(T)] For each $P \in \Ob(\TT)$ there is a sequence $P_0\leq P_1\leq \ldots$ of objects in $\TT$ such that $P = \bigcup_{i \geq 0} P_i$, and such that $P_i$ is a finite subgroup of $S$ for all $i \geq 0$.

\end{itemize}

\end{defi}

\begin{lmm}\label{equinerv}

Let $\g = \ploc$ be a $p$-local compact group, and let $\wL$ be a telescopic transporter system associated to $\FF$. Suppose in addition that $\wL$ contains $\LL$ as a full subcategory. Then, the inclusion $\LL \subseteq \wL$ induces an equivalence between the corresponding nerves.

\end{lmm}

\begin{proof}

This is an immediate consequence of \cite[Proposition A.9]{BLO6}.
\end{proof}

In terms of the above definition, in this section we study some situations where we can add objects to a given transporter system to produce a telescopic transporter system, without changing the homotopy type of the nerve of the original transporter system.

\begin{defi}\label{compsyst}

Let $S$ be a discrete $p$-toral group, let $\FF$ be a saturated fusion system over $S$, and let $\TT$ be a transporter system associated to $\FF$. Let also $(-)^{\star}_{\FF} \colon \FF \to \FF$ and $(-)^{\star}_{\TT} \colon \TT \to \TT$ be a pair of idempotent functors, and let $\cat^{\star} \subseteq \cat$, with $\cat = \FF, \TT$, be the full subcategory with $\Ob(\cat^{\star}) = \{P^{\star} \, | \, P \in \Ob(\cat)\}$. The pair $((-)^{\star}_{\FF}, (-)^{\star}_{\TT})$ is a \emph{finite retraction pair} if the following conditions are satisfied.
\begin{itemize}

\item[(1)] For each $P\leq S$, $P\leq (P)^{\star}_{\FF}$. Moreover, if $P \in \Ob(\TT)$, then $(P)^{\star}_{\FF} = (P)^{\star}_{\TT} = P^{\star}$.

\item[(2)] For each $P, Q\leq S$ and each $f \in \Hom_{\FF}(P,Q)$, $(f)^{\star}_{\FF} \in \Hom_{\FF}(P^{\star}, Q^{\star})$ extends $f$, and it is the unique extension.

\item[(i)] $\Ob(\FF^{\star})$ contains finitely many $S$-conjugacy classes of subgroups of $S$.

\item[(ii)] For all $P\leq S$, $P^{\star} = (P^{\star})^{\star}$.

\item[(iii)] If $P\leq Q\leq S$, then $P^{\star}\leq Q^{\star}$.

\item[(iv)] For all $P, Q\leq S$, $N_S(P,Q) \subseteq N_S(P^{\star}, Q^{\star})$.

\item[(v)] For all $P\leq S$, $C_S(P) = C_S(P^{\star})$.

\item[(a)] $(-)^{\star}_{\FF} \circ \rho = \rho \circ (-)^{\star}_{\TT}$.

\item[(b)] For all $P, Q \in \Ob(\TT)$ and all $\varphi \in \Mor_{\TT}(P,Q)$, we have $\varepsilon_{Q, Q^{\star}}(1) \circ \varphi = (\varphi)^{\star}_{\TT} \circ \varepsilon_{P, P^{\star}}(1)$.

\end{itemize}

\end{defi}

The above definition is inspired in the pair of ``bullet'' functors described in Proposition \ref{3.2BLO3}. In particular, conditions (i)-(v) and (a)-(b) are labelled to emphasize the relation with the motivating example. Properties (vi), (vii), (c) and (d) in \ref{3.2BLO3} are actually consequences of the definition, as we prove below.

\begin{lmm}\label{extraprop}

Let $\FF$ be a saturated fusion system over a discrete $p$-toral group $S$, let $\TT$ be a transporter system associated to $\FF$, and let $((-)^{\star}_{\FF}, (-)^{\star}_{\TT})$ be a finite retraction pair. Then the following properties hold (where we label the properties according to \ref{3.2BLO3} to emphasize the correspondence).
\begin{itemize}

\item[(vi)] The functor $(-)^{\star}_{\FF}$ is left adjoint to the inclusion of $\FF^{\star}$ as a full subcategory of $\FF$.

\item[(vii)] All $\FF$-centric $\FF$-radical subgroups of $S$ are in $\FF^{\star}$.

\item[(c)] For all $P, Q \in \Ob(\TT)$ and all $g \in N_S(P,Q)$, we have $(\varepsilon_{P,Q}(g))^{\star}_{\TT} = \varepsilon_{P^{\star}, Q^{\star}}(g)$.

\item[(d)] The functor $(-)^{\star}_{\TT}$ is left adjoint to the inclusion of $\TT^{\star}$ as a full subcategory of $\TT$. In particular, the inclusion $\TT^{\star} \subseteq \TT$ induces an equivalence $|\TT^{\star}| \simeq |\TT|$.

\end{itemize}

\end{lmm}

\begin{proof}

Property (vi) follows from condition (2) \ref{compsyst}, since it implies that, for all $P, Q\leq S$, the restriction map $\Hom_{\FF}(P^{\star}, Q^{\star}) \to \Hom_{\FF}(P, Q^{\star})$ is a bijection. To prove property (vii), let $P\leq S$ be $\FF$-centric, and suppose that $P \notin \Ob(\FF^{\star})$. Then, $P \lneqq P^{\star}$, which implies that $P \lneqq N_{P^{\star}}(P)$. Since every element of $\Aut_{\FF}(P)$ extends uniquely to an element of $\Aut_{\FF}(P^{\star})$, it is not hard to see that $1 \neq N_{P^{\star}}(P)/\Inn(P)$ is normalized by $\Aut_{\FF}(P)$, and hence $P$ cannot be $\FF$-radical. Property (c) is an immediate consequence of property (b) in \ref{compsyst}, applied to $\varphi = \varepsilon_{P,Q}(g)$, together with Lemma \ref{epimono}.

Finally, we prove property (d). For each $P \in \Ob(\TT)$, set as usual
$$
E(P) = \Ker(\Aut_{\TT}(P) \Right2{}\Aut_{\FF}(P)).
$$
We claim that $E(P) = E(P^{\star})$ for all $P \in \Ob(\TT)$. Clearly, restriction from $P^{\star}$ to $P$ maps $E(P^{\star})$ to $E(P)$, and this restriction is injective by \ref{epimono}. Let now $\varphi \in E(P)$, and consider $(\varphi)^{\star}_{\TT} \in \Aut_{\TT}(P^{\star})$. By assumption, $f = \rho(\varphi) = \Id$, and thus $(f)^{\star}_{\FF} = \Id \in \Aut_{\FF}(P^{\star})$ by property (2) in \ref{compsyst}. By property (a) in \ref{compsyst} we get
$$
\rho((\varphi)^{\star}_{\TT}) = (\rho(\varphi))^{\star}_{\FF} = (\Id)^{\star}_{\FF} = \Id,
$$
and $(\varphi)^{\star}_{\TT} \in E(P^{\star})$. Using axiom (A2) of transporter systems, together with property (vi) above, it is easy to deduce now that the functor $(-)^{\star}_{\TT}$ is left adjoint to the inclusion of $\TT^{\star}$ as a full subcategory of $\TT$, and property (d) follows.
\end{proof}

\begin{defi}\label{wT}

Let $((-)^{\star}_{\FF}, (-)^{\star}_{\TT})$ be a finite retraction pair. Define $\4{\TT}$ to be the category with object set $\Ob(\4{\TT}) = \{P\leq S \, | \, P^{\star} \in \Ob(\TT)\}$, and with morphism sets
$$
\Mor_{\4{\TT}}(P,Q) = \{\varphi \in \Mor_{\TT}(P^{\star}, Q^{\star}) \, | \, \varphi \circ \varepsilon_{P^{\star}}(g) \circ \varphi^{-1} \in \varepsilon_{Q^{\star}}(Q) \mbox{, for all } g \in P\},
$$
for all $P, Q \in \Ob(\4{\TT})$. Composition in $\4{\TT}$ is given by composition in $\TT$. Define also functors
$$
\TT_{\Ob(\4{\TT})}(S) \Right3{\4{\varepsilon}} \4{\TT} \qquad \mbox{and} \qquad \4{\TT} \Right3{\4{\rho}} \FF
$$
as follows. The functor $\4{\varepsilon}$ is the identity on objects, and the functor $\4{\rho}$ is injective on objects. For all $P, Q \in \Ob(\4{\TT})$, all $g \in N_S(P,Q)$, and all $\varphi \in \Mor_{\4{\TT}}(P,Q)$, define
$$
\4{\varepsilon}_{P,Q}(g) = \varepsilon_{P^{\bullet}, Q^{\bullet}}(g) \in \Mor_{\wL}(P,Q) \qquad \mbox{and} \qquad \4{\rho}(\varphi \colon P \to Q) = \rho(\varphi \colon P^{\bullet} \to Q^{\bullet})|_P.
$$
The properties of the functors $(-)^{\star}_{\FF}$ and $(-)^{\star}_{\TT}$ imply that both $\4{\TT}$ and the above functors are well-defined, and that $\TT$ is a full subcategory of $\4{\TT}$.

\end{defi}

\begin{prop}\label{extendL}

For each finite retraction pair $((-)^{\star}_{\FF}, (-)^{\star}_{\TT})$, the category $\4{\TT}$, with the functors $\4{\varepsilon}$ and $\4{\rho}$ defined above, is a telescopic transporter system associated to $\FF$. Furthermore, the functor $(-)^{\star}_{\TT} \colon \TT \to \TT$ extends to a functor $(-)^{\star}_{\4{\TT}} \colon \4{\TT} \to \TT$, which is unique satisfying the following properties
\begin{enumerate}[(a)]

\item There is an equality $(-)^{\star}_{\FF} \circ \4{\rho} = \rho \circ (-)^{\star}_{\4{\TT}} \colon \4{\TT} \to \FF$.

\item For all $P, Q \in \Ob(\4{\TT})$ and all $\varphi \in \Mor_{\4{\TT}}(P,Q)$, we have $\4{\varepsilon}_{Q,Q^{\star}}(1) \circ \varphi = (\varphi)^{\star}_{\4{\TT}} \circ \4{\varepsilon}_{P,P^{\star}}(1)$.

\end{enumerate}
In particular, the inclusion of $\TT$ in $\4{\TT}$ as a full subcategory induces an equivalence $|\TT| \simeq |\4{\TT}|$.

\end{prop}

\begin{proof}

By definition, $\4{\TT}$ contains $\TT$ as a full subcategory, and the functor $(-)_{\TT}^{\star} \colon \TT \to \TT^{\star}$ can be extended to a functor $(-)_{\4{\TT}}^{\star} \colon \4{\TT} \to \TT^{\star}$ as follows. On objects, $(P)^{\star}_{\4{\TT}} = (P)^{\star}_{\FF} = P^{\star}$. On morphisms, $(-)^{\star}_{\4{\TT}}$ is defined by the inclusion
$$
\Mor_{\4{\TT}}(P,Q) \subseteq \Mor_{\TT}(P^{\star}, Q^{\star})
$$
given by definition of $\4{\TT}$. The proof of (a) and (b), as well as the uniqueness of $(-)^{\star}_{\4{\TT}}$ satisfying these conditions, is left to the reader as an easy exercise.

Next we show that $\4{\TT}$ is indeed a transporter system. Conditions (A1), (B) and (C) are clear. Condition (A2) follows from the properties of the functor $(-)^{\star}_{\FF} \colon \FF \to \FF^{\star}$. Indeed, for each $P \in \Ob(\4{\TT})$ set
$$
\begin{array}{c}
\4{E}(P) = \Ker(\4{\rho}_P \colon \Aut_{\4{\TT}}(P) \to \Aut_{\FF}(P)) \\
E(P^{\star}) = \Ker(\rho_{P^{\star}} \colon \Aut_{\TT}(P^{\star}) \to \Aut_{\FF}(P^{\star})).
\end{array}
$$
If $\varphi \in E(P^{\star})$, then, by definition of $\4{\TT}$, together with axiom (C) on $\TT$, it follows that $\varphi \in \4{E}(P)$. Conversely, if $\varphi \in \4{E}(P)$, then by definition $\varphi \in \Aut_{\TT}(P^{\star})$ is such that $\rho(\varphi)|_P = \Id$. By property (a) on $(-)^{\star}_{\4{\TT}}$, we have
$$
\rho(\varphi) = \rho((\varphi)^{\star}_{\4{\TT}}) = (\rho(\varphi)|_P)^{\star}_{\FF} = (\Id)^{\star}_{\FF} = \Id,
$$
and thus $\varphi \in E(P^{\star})$. Hence $\4{E}(P) = E(P^{\star})$, and the freeness of the action of $\4{E}(P)$ on $\Mor_{\4{\TT}}(P,Q)$ follows from property (A2) in $\TT$. That $\4{\rho}_{P,Q}$ is the orbit map of this action follows easily.

Next we check condition (I) for $\4{\TT}$. Fix $Q \in \Ob(\4{\TT})$. If $Q \in \Ob(\TT)$ then there is nothing to show, since $\TT$ is a full subcategory of $\4{\TT}$. Thus, assume that $Q \notin \Ob(\TT)$. We can choose $Q$ such that $Q^{\star}$ is fully $\FF$-normalized, so that
$$
\varepsilon_{Q^{\star}}(N_S(Q^{\star})) \in \Syl_p(\Aut_{\TT}(Q^{\star})).
$$
Set for short $G = \Aut_{\TT}(Q^{\star})$ and $K = \varepsilon_{Q^{\star}}(N_S(Q^{\star}))$.

We claim first that every subgroup $H$ of $G$ has Sylow $p$-subgroups. Notice that $G/Q^{\star} \cong \Out_{\FF}(Q^{\star})$, which is a finite group. Thus, $H/(H \cap Q^{\star}) \cong HQ^{\star}/Q^{\star}\leq G/Q^{\star}$, and thus $H \cap Q^{\star}$ is a discrete $p$-toral normal subgroup of $H$ with finite index. The claim follows by \cite[Lemma 8.1]{BLO3}.

Now, by definition we can consider $H = \Aut_{\4{\TT}}(Q)$ as a subgroup of $G$, and in particular the above discussion implies that $H$ has Sylow $p$-subgroups. Fix $R \in \Syl_p(H)$ such that $\varepsilon_Q(N_S(Q))\leq R$. Since $K \in \Syl_p(G)$, there is some $\varphi \in G$ such that $\varphi \circ R \circ \varphi^{-1}\leq K$. Set $P = \rho(\varphi)(Q)\leq Q^{\star}$. Note that $P^{\star}\leq Q^{\star}$ by definition of $P$, and $P$ is $\FF$-conjugate to $Q$, which implies that $P^{\star}$ is $\FF$-conjugate to $Q^{\star}$ This implies that $P^{\star} = Q^{\star}$. Thus,
$$
\Aut_{\4{\TT}}(P) = \varphi \circ H \circ \varphi^{-1} = \varphi \circ \Aut_{\4{\TT}}(Q) \circ \varphi^{-1}\leq G,
$$
and $\varepsilon_P(N_S(P)) \in \Syl_p(\Aut_{\4{\TT}}(P))$. Condition (I) follows.

Condition (II) for $\4{\TT}$ follows easily from condition (II) for $\TT$. Indeed, let $\varphi \in \Iso_{\4{\TT}}(P,Q)$, $P \lhd \4{P}\leq S$ and $Q \lhd \4{Q}\leq S$ be such that $\varphi \circ \4{\varepsilon}_P(\4{P}) \circ \varphi^{-1}\leq \4{\varepsilon}_Q(\4{Q})$. By applying the functor $(-)_{\4{\TT}}^{\star} \colon \4{\TT} \to \TT^{\star}$, we get $\varphi^{\star} \in \Iso_{\TT}(P^{\star}, Q^{\star})$, and
$$
P^{\star} \lhd \widehat{P} \defin N_{(\4{P})^{\star}}(P^{\star})\leq S \qquad \mbox{and} \qquad Q^{\star} \lhd \widehat{Q} \defin N_{(\4{Q})^{\star}}(Q^{\star})\leq S,
$$
such that $\varphi^{\star} \circ \varepsilon_{P^{\star}}(\widehat{P}) \circ (\varphi^{\star})^{-1}\leq \varepsilon_{Q^{\star}}(\widehat{Q})$. Axiom (II) in $\TT$ implies that there exists some $\widehat{\varphi} \in \Mor_{\TT}(\widehat{P}, \widehat{Q})$ such that
$$
\widehat{\varphi} \circ \varepsilon_{P^{\star}, \widehat{P}}(1) = \varepsilon_{Q^{\star}, \widehat{Q}}(1) \circ \varphi^{\star}.
$$
Note that $\4{P}\leq \widehat{P}$ and $\4{Q}\leq \widehat{Q}$ by property (iv) in \ref{compsyst}. Thus, we may restrict the morphism $\widehat{\varphi}$ to $\4{P}$, and condition (II) follows.

Condition (III) for $\4{\TT}$ follows easily by condition (III) for $\TT$, together with the properties of the functor $(-)^{\star}_{\4{\TT}} \colon \4{\TT} \to \TT^{\star}$, since $\TT^{\star}$ contains finitely many isomorphism classes of objects by property (i) in \ref{compsyst}.

Let us now prove that that $\4{\TT}$ is a telescopic transporter system. Let $P \in \Ob(\4{\TT})$. By \cite[Lemma 1.9]{BLO3}, there is a sequence $P_0\leq P_1\leq \ldots$ of finite subgroups of $P$ such that $P = \bigcup_{i \geq 0} P_i$. Since $\FF^{\star}$ contains finitely many $S$-conjugacy classes of subgroups by property (i) in \ref{compsyst}, it follows that there is some $M \in \N$ such that $(P_i)^{\star} = P^{\star}$ for all $i \geq M$, and we may assume that $M = 0$ for simplicity. This way, since $P \in \Ob(\4{\TT})$, it follows that $P_i \in \Ob(\4{\TT})$ for all $i \geq 0$. Thus $\4{\TT}$ is a telescopic transporter system.

Finally, we check that the inclusion of $\TT$ in $\4{\TT}$ as a full subcategory induces an equivalence between the corresponding nerves. Recall from property (d) in \ref{extraprop} that the inclusion $\TT^{\star} \subseteq \TT$ induces an equivalence $|\TT^{\star}| \simeq |\TT|$. Thus, we only need to show that $(-)^{\star}_{\4{\TT}} \colon \4{\TT} \to \TT^{\star}$ is (left) adjoint to the inclusion of $\TT^{\star}$ as a full subcategory of $\4{\TT}$. That is, given $P, Q \in \Ob(\4{\TT})$, we have to show that the restriction map
$$
\Mor_{\4{\TT}}(P^{\star}, Q^{\star}) \Right2{} \Mor_{\4{\TT}}(P, Q^{\star})
$$
is a bijection. Let $\4{E}(P) = E(P^{\star})$ as above, and recall that there is a bijection between the sets $\Hom_{\FF}(P^{\star}, Q^{\star})$ and $\Hom_{\FF}(P, Q^{\star})$, given by the restriction map, by (vi) in \ref{extraprop}. Thus, by axiom (A2) of transporter systems,
$$
\Mor_{\TT}(P^{\star}, Q^{\star})/E(P^{\star}) = \Hom_{\FF}(P^{\star}, Q^{\star}) \cong \Hom_{\FF}(P,Q) = \Mor_{\4{\TT}}(P, Q^{\star})/\4{E}(P),
$$
and the claim follows.
\end{proof}

We call $\4{\TT}$ the \emph{telescopic transporter system associated to $((-)^{\star}_{\FF}, (-)^{\star}_{\TT})$}, or simply the telescopic transporter system associated to $\TT$ if there is no need to specify $((-)^{\star}_{\FF}, (-)^{\star}_{\TT})$.

\begin{prop}\label{extendL2}

Each $\Psi \in \Aut_{\typ}^I(\TT)$ extends uniquely to some $\4{\Psi} \in \Aut_{\typ}^I(\4{\TT})$.

\end{prop}

\begin{proof}

Let $\Psi \in \Aut_{\typ}^I(\TT)$ and let $\psi \in \Aut(S)$ be the automorphism induced by $\Psi$. Then $\Psi$ extends to $\4{\TT}$ by the formulas
$$
\4{\Psi}(P) = \psi(P) \qquad \mbox{and} \qquad \4{\Psi}(\varphi \colon P \to Q) = \Psi(\varphi^{\star}_{\4{\TT}} \colon P^{\star} \to Q^{\star}).
$$
Clearly, this determines an isotypical equivalence $\4{\Psi}$ of $\4{\TT}$. Moreover, since $\4{\Psi}$  is isotypical, that is $\4{\Psi}(\4{\varepsilon}_{P,Q}(1)) = \4{\varepsilon}_{\4{\Psi}(P), \4{\Psi}(Q)}(1)$ for all $P, Q \in \Ob(\4{\TT})$ with $P\leq Q$, and since morphisms in $\4{\TT}$ are monomorphisms and epimorphisms in the categorical sense by Lemma \ref{epimono}, it follows that $\4{\Psi}$ is the unique extension of $\Psi$ to $\4{\TT}$.
\end{proof}

Below we analyze some examples which will be of interest in later sections.

\begin{expl}\label{expl1}

We start with the most obvious example. Let $\g = \ploc$ be a $p$-local compact group, and let $(-)^{\bullet}_{\FF} \colon \FF \to \FF$ and $(-)^{\bullet}_{\LL} \colon \LL \to \LL$ be the usual ``bullet'' functors. Then, clearly $((-)^{\bullet}_{\FF}, (-)^{\bullet}_{\LL})$ is a finite retraction pair (by \ref{3.2BLO3}), and this way we obtain a telescopic transporter system $\wL$ which in addition satisfies the following properties.
\begin{itemize}

\item[(1)] For each $P \in \Ob(\wL)$, we have 
$$
\Ker(\Aut_{\wL}(P) \to \Aut_{\FF}(P)) = \varepsilon_P(C_S(P)) = \varepsilon_P(Z(P^{\bullet})).
$$

\end{itemize}
Let $P \in \Ob(\wL)$, and set $E(P) = \Ker(\Aut_{\wL}(P) \to \Aut_{\FF}(P))$. If $P \in \Ob(\LL)$ then $\Aut_{\wL}(P) = \Aut_{\LL}(P)$, and there is nothing to prove. Suppose that $P \notin \Ob(\LL)$. By definition $P^{\bullet} \in \Ob(\LL)$, and there is a commutative diagram of group extensions
$$
\xymatrix{
E(P^{\bullet}) \ar[r] & \Aut_{\wL}(P^{\bullet}) \ar[r] & \Aut_{\FF}(P^{\bullet}) \\
E(P) \ar[r] \ar[u] & \Aut_{\LL}(P) \ar[r] \ar[u] & \Aut_{\FF}(P) \ar[u]
}
$$
where all the vertical arrows are inclusions. Thus, we have
$$
E(P)\leq E(P^{\bullet}) = \varepsilon_{P^{\bullet}}(Z(P^{\bullet})) = \varepsilon_{P^{\bullet}}(C_S(P^{\bullet})) = \varepsilon_P(C_S(P)),
$$
where $C_S(P^{\bullet}) = Z(P^{\bullet})$ since $P^{\bullet}$ is $\FF$-centric, and $C_S(P^{\bullet}) = C_S(P)$ by property (v) in \ref{compsyst}. The inclusion $\varepsilon_P(C_S(P))\leq E(P)$ is clear. This proves (1). In particular, every object in $\wL$ is \emph{quasicentric} (that is $C_{\FF}(P)$ is the fusion system of $C_S(P)$ for all $P \in \Ob(\wL)$), and in this sense $\wL$ is a \emph{quasicentric linking system}.
\begin{itemize}
 
\item[(2)] There is an isomorphism $\Aut_{\typ}^{I}(\LL) \cong \Aut_{\typ}^{I}(\wL)$.

\end{itemize}
This follows by Proposition \ref{extendL2}, together with the observation that every isotypical automorphism of $\wL$ must restrict to an isotypical automorphism of $\LL$, since $\Ob(\LL)$ is the set of all $\FF$-centric subgroups of $S$. Moreover, this restriction is injective as a consequence of Lemma \ref{epimono},  and since every morphism in $\wL$ is the restriction of some morphism in $\LL$.

\end{expl}

\begin{expl}\label{expl3}

The following is a less obvious example. Again, let $\g = \ploc$ be a $p$-local compact group, let $((-)^{\bullet}_{\FF}, (-)^{\bullet}_{\LL})$ be the finite retraction pair in \ref{expl1}, and let $\wL$ be the telescopic transporter system associated to $\LL$. Let also $A\leq S$ be a fully $\FF$-normalized subgroup such that $N_S(A)$ has finite index in $S$ (for example $A\leq T$ a subgroup of the maximal torus of $S$).

Let $N_{\g}(A) = (N_S(A), N_{\FF}(A), N_{\LL}(A))$ be the normalizer $p$-local compact group of $A$, as defined in \ref{rmknorm}. Again, in general $N_{\LL}(A)$ is not a telescopic linking system. In this case, since $N_{\g}(A)$ is a $p$-local compact group, one could apply Example \ref{expl1} to produce a telescopic transporter system associated to $N_{\LL}(A)$. However, this way there is no obvious correspondence between the telescopic transporter systems for $N_{\LL}(A)$ and $\LL$, mainly because usually the ``bullet'' functors in $N_{\FF}(A)$ and in $\FF$ do not agree, that is $(P)^{\bullet}_{\FF} \neq (P)^{\bullet}_{N_{\FF}(A)}$ in general.

Instead, we propose a different construction. Set for short $N = N_S(A)$, $\ee = N_{\FF}(A)$ and $\TT = N_{\LL}(A)$. We define a finite retraction pair $((-)^{\star}_{\ee}, (-)^{\star}_{\TT})$ as follows. For each $P\leq N$, notice that $(P)^{\bullet}_{\FF}\leq N_S(A)$, since $N_S(A)$ already contains the maximal torus of $S$. Thus, we can define
$$
P^{\star} = P^{\bullet}.
$$
On morphisms, let $P, Q\leq N$, and let $f \in \Hom_{\ee}(P,Q)$. By definition of $\ee$, the morphism $f$ extends to some $\gamma \in \Hom_{\FF}(PA, QA)$ such that $\gamma|_A \in \Aut_{\FF}(A)$. Applying the functor $(-)^{\bullet}_{\FF}$ to the commutative square
$$
\xymatrix{
PA \ar[r]^{\gamma} & QA\\
P \ar[u]^{\incl} \ar[r]_{f} & Q \ar[u]_{\incl}
}
$$
we see that $(f)^{\bullet}_{\FF}$ extends to $(\gamma)^{\bullet}_{\FF}$, and the latter restricts in turn to a morphism $\5{\gamma} \in \Hom_{\FF}((P)^{\bullet}_{\FF}A, (Q)^{\bullet}_{\FF}A)$ such that $\5{\gamma}|_A = \gamma|_A \in \Aut_{\FF}(A)$. We define
$$
(f)^{\star}_{\ee} = (f)^{\bullet}_{\FF}.
$$
Properties (1)-(2) and (i)-(v) in \ref{compsyst} for $(-)^{\bullet}_{\FF}$ imply that $(-)^{\star}_{\ee}$ also satisfies these conditions.

On $\TT$, define $(-)^{\star}_{\TT}$ as follows. Let $P, Q \in \Ob(\TT)$, and let $\varphi \in \Mor_{\TT}(P,Q)$. By definition, $\varphi$ is a morphism in $\Mor_{\LL}(PA, QA)$ such that $\rho(\varphi)|_P \in \Mor_{\ee}(P,Q)$, and $\rho(\varphi)|_A \in \Aut_{\FF}(A)$. Clearly,
$$
P^{\bullet}A\leq (PA)^{\bullet} \qquad \mbox{and} \qquad Q^{\bullet}A\leq (QA)^{\bullet},
$$
and thus $(\varphi)^{\bullet}_{\LL}$ restricts to a morphism $\4{\varphi} \in \Mor_{\LL}(P^{\bullet}A, Q^{\bullet}A)$ such that $\rho(\4{\varphi})|_{P^{\bullet}} \in \Hom_{\ee}(P^{\bullet}, Q^{\bullet})$ and $\rho(\4{\varphi})|_A \in \Aut_{\FF}(A)$. Define $(\varphi)^{\star}_{\TT} = \4{\varphi}$. It is not difficult to check that $(-)^{\bullet}_{\TT}$ satisfies properties (a)-(b) in \ref{compsyst}.

Let $\wL$ and $\4{\TT}$ be the associated telescopic transporter systems for $\LL$ and $\TT$, respectively. In general, $\TT$ is not a subcategory of $\LL$, and neither is $\4{\TT}$ a subcategory of $\wL$. Let $\TT_{\geq A} \subseteq \TT$ be the full subcategory of subgroups that contain $A$, and let $\4{\TT}_{\geq A} \subseteq \4{\TT}$ be the full subcategory of subgroups $P$ such that $A\leq P^{\bullet}$. Then the following is easily checked.
\begin{enumerate}[(a)]

\item $\TT_{\geq A}$ is a subcategory of $\LL$, and $\4{\TT}_{\geq A}$ is a subcategory of $\wL$.

\item $\TT_{\geq A}$ contains all the centric radical subgroups of $\ee$.

\item The functor $(-)^{\star}_{\4{\TT}}$ coincides with the functor $(-)^{\bullet}_{\wL}$ on the subcategory $\4{\TT}_{\geq A}$.

\end{enumerate}
This example, including the above remarks, will be very useful in the next section, when we have to compare certain constructions on a $p$-local compact group $\g = \ploc$ and on the normalizer $N_{\g}(A) = (N_S(A), N_{\FF}(A), N_{\LL}(A))$ of a certain subtorus $A\leq S$.

\end{expl}


\section{Families of unstable Adams operations}\label{Sfam}

In this section we prove Theorem \ref{thmA}, restated as Theorem \ref{fix6} below: every $p$-local compact group can be approximated by $p$-local finite groups. Roughly speaking, given a $p$-local compact group we produce an \emph{approximation of $\g$ by $p$-local finite groups} (defined below) by considering \emph{fixed point subcategories} of a telescopic transporter system associated to $\LL$ by iterations of a given unstable Adams operation on $\g$. 

Essentially, we follow the same lines as \cite{Gonza1}. However, the introduction of telescopic transporter systems means a great deal of simplification, and it is actually thanks to this that we are finally able to prove Proposition \ref{fix5}, basically the missing step in \cite{Gonza1} in proving Theorem \ref{thmA}. We have opted for reproving here every property that we need from \cite{Gonza1} for the sake of completeness as well as for correcting mistakes: while working on \ref{fix5} below, the author realized that the statement of \cite[Lemma 2.11]{Gonza1} is false. Nevertheless, this does not affect neither the main results of \cite{Gonza1} nor the results that we present in this paper, and this comment is just intended as a warning to the interested reader.

\begin{defi}\label{defiapprox}

Let $\g = \ploc$ be a $p$-local compact group, and let $\wL$ be a telescopic transporter system associated to $\FF$ and containing $\LL$ as a full subcategory. An \emph{approximation of $\g$ by $p$-local finite groups} is a family $\{(S_i, \FF_i, \LL_i)\}_{i \geq 0}$ satisfying the following conditions.
\begin{enumerate}[(i)]

\item $S = \bigcup_{i \geq 0} S_i$.

\item For each $i \geq 0$, $S_i$ is a finite $p$-group, $\FF_i$ is a saturated fusion system over $S_i$, and $\LL_i$ is a linking system associated to $\FF_i$. Furthermore, $\Ob(\FF_i^{cr}) \subseteq \Ob(\LL_i)$, and there are inclusions $\LL_i \subseteq \LL_{i+1}$ and $\LL_i \subseteq \wL$.

\item For each $P, Q \in \Ob(\wL)$ and each $\varphi \in \Mor_{\wL}(P,Q)$ there exists some $M \in \N$ such that, for all $i \geq M$, there are objects $P_i, Q_i \in \Ob(\LL_i)$ and morphisms $\varphi_i \in \Mor_{\LL_i}(P_i, Q_i)$, such that $P = \bigcup_{i \geq M} P_i$ and $Q = \bigcup_{i \geq M} Q_i$, and $\4{\varepsilon}_{Q_i, Q}(1) \circ \varphi_i = \varphi \circ \4{\varepsilon}_{P_i, P}(1)$.

\end{enumerate}

\end{defi}

Although condition (i), or at least a weaker version of it, can be deduced from condition (iii) applied to $P = Q = S$ and to any $\varphi \in \Aut_{\LL}(S)$, we prefer to include (i)  in the definition for the sake of clarity. We now show some basic properties of approximations of $p$-local compact groups by $p$-local finite groups.

\begin{lmm}\label{finmorph}

Let $\g = \ploc$ be a $p$-local compact group, and let $\{(S_i, \FF_i, \LL_i)\}_{i \geq 0}$ be a finite approximation of $\g$ by $p$-local finite groups with respect to some telescopic transporter system $\wL$ satisfying the conditions in \ref{defiapprox}. Then, for every finite subgroup $P\leq S$ and every $f \in \Hom_{\FF}(P,S)$, there exists some $M \in \N$ such that, for all $i \geq M$, $P\leq S_i$, and $f$ is the composition of a morphism $\gamma \in \Hom_{\FF_i}(P, S_i)$ with the inclusion $S_i\leq S$.

\end{lmm}

\begin{proof}

Let $Q = f(P)$, which is also a finite subgroup of $S$, and let $\gamma \in \Iso_{\FF}(P,Q)$ be the restriction of $f$ to its image. By property (i) in \ref{defiapprox}, it is clear that there is some $M_0 \in \N$ such that $P, Q\leq S_i$ for all $i \geq M_0$. By Alperin's Fusion Theorem \cite[Theorem 3.6]{BLO3}, there exist $W_0 = P, W_1, \ldots, W_n = Q\leq S$, $U_1, \ldots, U_n \in \Ob(\LL) \subseteq \Ob(\wL)$, and morphisms $\varphi_j \in \Aut_{\wL}(U_j)$, for $j = 1, \ldots, n$, such that, for each $j$
$$
W_{j-1}, W_j\leq U_j \qquad \mbox{and} \qquad \rho(\varphi_j)(W_{j-1}) = W_j,
$$
and $\gamma = \rho(\varphi_n) \circ \ldots \circ \rho(\varphi_1)$. Combining properties (i) and (iii) in \ref{defiapprox}, we see that for each $j = 1, \ldots, n$ there exists some $M_j \in \N$ such that, for all $i \geq M_j$, there exist $U_{j,i}, V_{j, i} \in \Ob(\LL_i)$, together with an isomorphism $\varphi_{j.i} \in \Iso_{\LL_i}(U_{j,i}, V_{j,i})$, such that
$$
U_j = \bigcup_{i \geq M_j} U_{j,i} = \bigcup_{i \geq M_j} V_{j,i} \qquad \mbox{and} \qquad \4{\varepsilon}_{V_{j,i}, U_j}(1) \circ \varphi_{j,i} = \varphi_j \circ \4{\varepsilon}_{U_{j,i}, U_j}(1).
$$
Moreover, since $W_{j-1}, W_j$ are finite subgroups, we may assume without loss of generality that $W_{j-1}\leq U_{j,i}$ and $W_j\leq V_{j, i}$ for all $i \geq M_j$. Let $M = \max\{M_0, \ldots, M_n\}$. Then, for all $i \geq 0$, it follows that
$$
\gamma = \rho_i(\varphi_{n,i}) \circ \ldots \circ \rho_i(\varphi_{1,i}) \in \Mor(\FF_i),
$$
and this finishes the proof.
\end{proof}

\begin{lmm}\label{Quill}

Let $\cat$ be a nonempty category that satisfies the following conditions:
\begin{enumerate}[(i)]

\item given objects $a_1, a_2$, there is an object $b$ and morphisms $a_1 \Right1{g_1} b \Left1{g_2} a_2$; and

\item given morphisms $g_1, g_2 \colon a \to b$, there is $h \colon b \to c$ such that $h \circ g_1 = h \circ g_2$.

\end{enumerate}
Then, the nerve of $\cat$ is contractible.

\end{lmm}

\begin{proof}

This is \cite[Corollary 2]{Quillen}.
\end{proof}

\begin{lmm}\label{approx-2}

Let $\cat$ be a small category all of whose morphisms are epimorphisms and monomorphisms in the categorical sense, and let $\cat_0 \subseteq \cat_1 \subseteq \ldots $ be a sequence of subcategories such that $\cat = \bigcup_{i \geq 0} \cat_i$. Then, $|\cat| \simeq \hocolim |\cat_i|$.

\end{lmm}

\begin{proof}

Let $I$ be the category of natural ordinals, with objects $\Ob(I) = \{i \in \N\}$, and where the morphism set $\Mor_I(i,j)$ is $\{\sigma_{i,j}\}$ if $i\leq j$, or empty otherwise. Define a functor
$$
\Theta \colon I \Right3{} \Cat
$$
by $\Theta(i) = \cat_i$ and $\Theta(\sigma_{i,j}) = \incl \colon \cat_i \to \cat_j$. The Grothendieck construction on $\Theta$, namely $G(\Theta)$, is the category with object set
$$
\{(i, X) \, | \, i \in \Ob(I) \mbox{ and } X \in \Ob(\cat_i)\}.
$$
The morphism sets $\Mor_{G(\Theta)}((i,X), (j, Y))$ are empty whenever $j < i$. Otherwise they consist of the pairs $(\sigma_{i,j}, \varphi)$, with $\varphi \in \Mor_{\cat_j}(X,Y)$. By \cite[Theorem 1.2]{Thomason}, we have an equivalence $\hocolim |\cat_i| \simeq |G(\Theta)|$.

Consider now the projection functor $\tau \colon G(\Theta) \to \cat$ that sends an object $(i, X)$ to $\tau(i, X) = X$, and a morphism $(\sigma_{i,j}, \varphi)$ to $\tau(\sigma_{i,j}, \varphi) = \varphi$. We claim that this functor induces an equivalence between the corresponding nerves. For each $X \in \cat$, let $\tau/X$ be the category with object set
$$
\Ob(\tau/X) = \{((j, Y), \varphi) \,\, | \,\, (j, Y) \in \Ob(G(\Theta)) \mbox{ and } \varphi \in \Mor_{\cat}(Y, X)\}.
$$
A morphism in $\tau/X$ from $((j,Y), \varphi)$ to $((k, Z), \psi)$ is $(\sigma_{i,j}, \gamma) \in \Mor_{\cat}((j,Y), (k, Z))$ such that $\varphi = \psi \circ \gamma$. By \cite[Theorem A and Corollary 2]{Quillen}, it is enough to check that $\tau/X$ satisfies the conditions of Lemma \ref{Quill}. Clearly, $\tau/X$ is nonempty, since $((i, X), 1_X) \in \tau/X$ for some $i \in \N$ big enough. Let $((j,Y), \varphi), ((k,Z), \psi) \in \Ob(\tau/X)$, and let $m = \max\{i,j,k\}$. Then condition (i) in Lemma \ref{Quill} holds with
$$
((j,Y), \varphi) \Right3{(\sigma_{j,m},\varphi)} ((m, X), 1_X) \Left3{(\sigma_{k,m}, \psi)} ((j,Y), \psi).
$$
Regarding condition (ii), notice that $\Mor_{\tau/X}(((j,Y), \varphi), ((k, Z), \psi))$ is either empty, or contains a single morphisms, since morphisms in $\cat$ are all epimorphisms and monomorphisms in the categorical sense. Thus, $|\tau/X|$ is contractible for all $X \in \cat$, and the claim follows.
\end{proof}

\begin{lmm}\label{approx0}

Let $\g = \ploc$ be a $p$-local compact group, and suppose $\g$ admits an approximation by $p$-local finite groups $\{(S_i, \FF_i, \LL_i)\}_{i \geq 0}$ with respect to some telescopic transporter system $\wL$ satisfying the conditions in \ref{defiapprox}. Then, there is an equivalence $B\g \simeq (\hocolim |\LL_i|)^{\wedge}_p$.

\end{lmm}

\begin{proof}

By \ref{defiapprox}, $\wL$ contains $\LL$ as a full subcategory, and thus by \ref{equinerv} the inclusion $\LL \subseteq \wL$ induces an equivalence $|\LL| \simeq |\wL|$. It is enough to show that $\hocolim |\LL_i| \simeq |\wL|$.

Set $\LL^{\circ} \defin \bigcup_{i \geq 0} \LL_i \subseteq \wL$. By Lemma \ref{epimono}, all morphisms in $\wL$ are epimorphisms and monomorphisms in the categorical sense, and thus the same applies to $\LL^{\circ}$. By Lemma \ref{approx-2} it follows that
$$
\hocolim |\LL_i| \simeq |\LL^{\circ}|.
$$
Thus, to finish the proof it is enough to show that the inclusion functor $\iota \colon \LL^{\circ} \to \wL$ induces an equivalence of nerves. For each $P \in \Ob(\wL)$, the undercategory $\iota/P$ has object set
$$
\Ob(\iota/P) = \{(Q, \varphi) \, | \, Q \in \Ob(\LL^{\circ}) \mbox{ and } \varphi \in \Mor_{\wL}(Q,P)\}.
$$
A morphism in $\iota/P$ from $(Q, \varphi)$ to $(R, \psi)$ is a morphism $\gamma \in \Mor_{\LL^{\circ}}(Q,R)$ such that $\varphi = \psi \circ \gamma$.

We show that $\iota/P$ satisfies the conditions of Lemma \ref{Quill}, which implies that $|\iota/P|$ is contractible. Clearly, $\iota/P$ is nonempty. Let $(Q, \varphi), (R, \psi) \in \Ob(\iota/P)$. By property (iii) in Definition \ref{defiapprox}, there is some $X \in \Ob(\LL^{\circ})$, with $X\leq P$, such that $\varphi$ and $\psi$ restrict to morphisms $\varphi \colon Q \to X$ and $\psi \colon R \to X$ in $\LL^{\circ}$. Thus, condition (i) of \ref{Quill} is satisfied with
$$
(Q, \varphi) \Right3{\varphi} (X, \varepsilon(1)) \Left3{\psi} (R, \psi).
$$
Regarding condition (ii) in \ref{Quill}, the set $\Mor_{\iota/P}((Q, \varphi), (R, \psi))$ is either empty or contains a single morphism. Since the argument works for all $P \in \Ob(\LL)$, it follows that $|\LL^{\circ}| \simeq |\wL|$.
\end{proof}

\begin{rmk}\label{approx-1}

Suppose the $p$-local compact group $\g = \ploc$ admits an approximation by $p$-local finite groups $\{(S_i, \FF_i, \LL_i)\}_{i \geq 0}$. For each $i$, the fusion system $\FF_i$ is saturated, and we may consider its associated centric linking system $\TT_i$. Let also $\hh_i = \Ob(\LL_i) \cap \Ob(\TT_i)$, and let $\LL_{\hh_i} \subseteq \LL_i$ be the full subcategory with object set $\hh_i$. Since both $\Ob(\LL_i)$ and $\Ob(\TT_i)$ contain $\Ob(\FF_i^{cr})$, it follows that $\Ob(\FF_i^{cr}) \subseteq \hh_i$. Moreover, there is a commutative diagram
$$
\xymatrix@R=1.2cm{
\ldots & \TT_i & \TT_{i+1} & \ldots \\
\ldots \ar[r] & \LL_{\hh_i} \ar[u]^{\iota_i} \ar[d]_{\incl} \ar[r]^{\incl} & \LL_{\hh_{i+1}} \ar[u]_{\iota_{i+1}} \ar[d]^{\incl} \ar[r] & \ldots \\
\ldots \ar[r] & \LL_i \ar[r]_{\incl} & \LL_{i+1} \ar[r] & \ldots
}
$$
where all the vertical arrows induce homotopy equivalences between the realizations of the corresponding nerves, by \cite[Theorem B]{BCGLO1}. Thus, if we denote $B\g_i = |\TT_i|^{\wedge}_p$, then Lemma \ref{approx0} implies that $B\g \simeq (\hocolim B\g_i)^{\wedge}_p$.

\end{rmk}


\subsection{Preliminary constructions}\label{Sprelim}

In this subsection we establish the notation and basic facts necessary for the proof that every $p$-local compact group has an approximation by $p$-local finite groups, in the next subsection. 

\begin{hyp}\label{hyp1}

For the rest of this subsection, let $\g = \ploc$ be a $p$-local compact group, and let $((-)^{\star}_{\FF}, (-)^{\star}_{\LL})$ be a finite retraction pair. Let also $\wL$ be the associated telescopic transporter system, and let $\Psi$ be a fine unstable Adams operation on $\LL$ (in the sense of \ref{uAo}). By a slight abuse of notation, we denote by $\Psi$ the corresponding extension of $\Psi$ to $\wL$ (see \ref{extendL2}), which is again a fine unstable Adams operation. Set $\Psi_0 = \Psi$, and for all $i \geq 0$, define
\begin{enumerate}[(a)]

\item $S_i = C_S(\Psi_i) = \{x \in S \, | \, \Psi_i(x) = x\}$; and

\item $\Psi_{i+1} = (\Psi_i)^p$.

\end{enumerate}

\end{hyp}

\begin{lmm}\label{SiS}

The following properties hold.
\begin{enumerate}[(i)]

\item $S = \bigcup_{i \geq 0} S_i$.

\item There is some $M_a \in \N$ such that $(S_i)^{\star} = S$ for all $i \geq M_a$.

\end{enumerate}

\end{lmm}

\begin{proof}

Let $T\leq S$ be the maximal torus of $S$, and set $T_i = T \cap S_i$. Notice that by definition $T_i$ is the subgroup of $T$ of elements fixed by $\Psi_i$. By hypothesis, $\Psi = \Psi_0$ has degree $\zeta \in 1 + p^m \Z^{\wedge}_p$ for some $m > 0$, and $\Psi_{i+1} = (\Psi_i)^p$. Thus, we have $T_i \lneqq T_{i+1}$ for all $i \geq 0$, and $T = \bigcup_{i \geq 0} T_i$. By \cite[Lemma 2.6]{JLL}, there is a subgroup $H\leq S_0$ such that $S = H \cdot T$. Thus, $H\leq S_i$ for all $i$, and we get
$$
S = H \cdot T = \bigcup_{i \geq 0} H \cdot T_i \subseteq \bigcup_{i \geq 0} S_i.
$$
This proves part (i). To prove part (ii), suppose otherwise that no such $M_a \in \N$ exists, that is, $(S_i)^{\star} \lneqq S$ for all $i \geq 0$. Since $\Ob(\FF^{\star})$ contains only finitely many conjugacy classes of elements, this means that there exist some $R \in \Ob(\FF^{\star})$ and some $M \in \N$ such that $R \lneqq S$ and such that $(S_i)^{\star} = R$ for all $i \geq M$. Notice that this contradicts part (i): if $R \lneqq S$, then there is some $x \in S \setminus R$. On the other hand, by part (i) we have $x \in S_i$ for $i$ big enough, and thus $x \in (S_i)^{\star} = R$, hence a contradiction.
\end{proof}

For simplicity we may assume that $M_a = 0$. In particular, $S_i \in \Ob(\wL)$ for all $i \geq 0$.

\begin{defi}\label{Li}

With the conventions above, for each $i \geq 0$ define $\LL_i$ as the category with object set $\Ob(\LL_i) = \{P\leq S_i \, | \, P \in \Ob(\wL)\}$, and with morphism sets
$$
\Mor_{\LL_i}(P,Q) = \{\varphi \in \Mor_{\wL}(P,Q) \, | \, \Psi_i(\varphi) = \varphi\}.
$$
Define also $\FF_i$ as the fusion system over $S_i$ generated by the restriction of $\4{\rho} \colon \wL \to \FF$ to $\LL_i$ (i.e. $\FF_i$ is $\Ob(\LL_i)$-generated). Finally, define functors
$$
\TT_{\Ob(\LL_i)}(S_i) \Right3{\varepsilon_i} \LL_i \qquad \mbox{and} \qquad \LL_i \Right3{\rho_i} \FF_i
$$
as the obvious restrictions of the structural functors $\4{\varepsilon} \colon \TT_{\Ob(\wL)}(S) \to \wL$ and $\4{\rho} \colon \wL \to \FF$.

\end{defi}

Despite its simplicity, the following example illustrates why it is necessary to work with a telescopic linking system rather than a centric linking system.

\begin{expl}\label{expl0}

Let $T$ be a discrete $p$-torus, i.e. $T \cong (\prufferp)^r$ for some $r \geq 1$. Let also $\g = \ploc$ be the \emph{trivial} $p$-local compact group associated to $T$. That is, $S = T$ and $\FF = \FF_T(T)$ is the fusion system over $T$ whose only morphisms are inclusions (since $T$ is abelian). This fusion system is obviously saturated, and has only one centric object, namely $T$ itself. Thus, $\LL$ has a single object, $T$, with $\Aut_{\LL}(T) = T$. On the other hand, since $T = S$, we have $P^{\bullet} = T$ for all $P\leq T$, and the telescopic linking system $\wL$ associated to $\LL$ in \ref{expl1} is the actual transporter category of the group $T$. That is, $\Ob(\wL) = \{P\leq T\}$, and $\Mor_{\wL}(P,Q) = N_T(P,Q) = T$ for all $P, Q\leq T$. Let now $\Psi$ be an unstable Adams operation as fixed in \ref{hyp1}. An easy computation reveals that $C_T(\Psi)$ must be a finite subgroup of $T$, and thus is not an object in $\LL$. In particular, without replacing $\LL$ by $\wL$, the subcategories $\LL_i$ defined above would be empty for all $i \geq 0$.

\end{expl}

\begin{prop}\label{fix1}

The following holds.
\begin{enumerate}[(i)]

\item For each $P \in \Ob(\FF^{\star})$ there exists some $M_P \in \N$ such that $(P \cap S_i)^{\star} = P$ for all $i \geq M_P$.

\item For each $\varphi \in \Mor(\wL)$ there exists some $M_{\varphi} \in \N$ such that $\Psi_i(\varphi) = \varphi$ for all $i \geq M_{\varphi}$.

\item For each $i \geq 0$ and each $x \in S$, $x^{-1} \cdot \Psi_i(x) \in T$.

\end{enumerate}

\end{prop}

\begin{proof}

To prove part (i), notice that $P^{\star} = P = \bigcup_{i \geq 0} P \cap S_i$ by Lemma \ref{SiS} (i). Suppose that $(P \cap S_i)^{\star} \lneqq P$ for all $i \geq 0$. Since $\FF^{\star}$ contains finitely many conjugacy classes of objects, this means that there is some $R \in \Ob(\FF^{\star})$ such that $(P \cap S_i)^{\star} = R \lneqq P$ for all $i$ big enough, contradicting the identity $P = \bigcup_{i \geq 0} P \cap S_i$. Part (iii) follows immediately by definition of unstable Adams operation, since the  $\Psi_i \in \Aut(S)$ induces the identity on $S/T$.

Finally, part (ii) follows by construction of the unstable Adams operations $\Psi_i$. More specifically, as stablished in \ref{hyp1}, the unstable Adams operation $\Psi$ satisfies the following property (see Remark \ref{uAo1}, or the proof of \cite[Theorem 4.1]{JLL} for a more detailed explanation): there is a (finite) set $\mathcal{M}$ of morphisms in $\LL^{\star}$ such that the following holds
\begin{enumerate}[(1)]

\item $\Psi(\varphi) = \varphi$ for all $\varphi \in \mathcal{M}$; and

\item every morphism $\psi$ in $\LL$ (and hence in $\wL$ by definition) decomposes as $\psi = \varepsilon(g) \circ \varphi$, where $\varphi$ is (the restriction of) a morphism in $\mathcal{M}$, and $g$ is an element of $S$.

\end{enumerate}
Moreover, these properties depend only on $\LL^{\star}$ containing finitely many $S$-conjugacy classes, and on $\LL$ being a centric linking system, but not on the functor $(-)^{\star}_{\LL}$. Fix $\psi \in \Mor(\wL)$, and let $\psi = \varepsilon(g) \circ \varphi$ be the corresponding decomposition, as described above. By assumption, $\Psi(\varphi) = \varphi$, and thus $\Psi_i(\varphi) = \varphi$ for all $i \geq 0$. Also, $S = \bigcup_{i \geq 0} S_i$ by Lemma \ref{SiS} (i), and thus there exists some $M_{\varphi} \in \N$ such that $g \in S_i$ for all $i \geq M_{\varphi}$. It follows that $\Psi_i(\psi) = \Psi_i(\varepsilon(g)) \circ \varphi = \varepsilon(g) \circ \varphi = \psi$, and part (ii) follows.
\end{proof}

\begin{prop}\label{fix2}

There exists some $M_b \in \N$ such that, for all $i \geq M_b$, the triple $(\LL_i, \varepsilon_i, \rho_i)$ is a transporter system associated to the fusion system $\FF_i$.

\end{prop}

\begin{proof}

We have to check the axioms in Definition \ref{defitransporter}. Notice that $\LL_i$ is a finite category, and thus we do not have to deal with axiom (III). Axioms (A1), (B) and (C) follow immediately by definition of $\LL_i$ as a subcategory of $\wL$. We deal with the remaining axioms of transporter systems in separate steps for the reader's convenience.

For each $P \in \Ob(\LL_i)$, set
$$
E_i(P) = \Ker(\Aut_{\LL_i}(P) \Right2{} \Aut_{\FF_i}(P)).
$$
Note that, via the inclusion $\Aut_{\LL_i}(P) \to \Aut_{\wL}(P^{\star})$, and by property (a) in Definition \ref{compsyst} of the functor $(-)^{\star}_{\LL}$, the subgroup $E_i(P)$ is mapped to a subgroup of $E(P^{\star}) = \Ker(\Aut_{\wL}(P^{\star}) \to \Aut_{\FF}(P^{\star}))$. Indeed, if $\varphi \in E_i(P)$, then $f = \rho_i(\varphi) = \Id$. Thus, $\4{\rho}((\varphi)^{\star}_{\LL}) = (\Id)^{\star}_{\FF} = \Id$.

\textbf{Step 1.} Axiom (A2) is satisfied: for each $P, Q \in \Ob(\LL_i)$, the group $E_i(P)$ acts freely on $\Mor_{\LL_i}(P,Q)$ by right composition, and $\rho_i \colon \Mor_{\LL_i}(P,Q) \to \Hom_{\FF_i}(P,Q)$ is the orbit map of this action. Also, $E_i(Q)$ acts freely on $\Mor_{\LL_i}(P,Q)$ by left composition.

The freeness of both the left action by $E_i(P)$ and the right action of $E_i(Q)$ follows by Lemma \ref{epimono}, which states that morphisms in $\wL$ (and in particular in $\LL_i$ by definition) are monomorphisms and epimorphisms in the categorical sense. That $\rho_i$ is the orbit map of the left conjugation action of $E_i(P)$ is now immediate.

\textbf{Step 2.} Axiom (I) is satisfied. In fact, since $\LL_i$ is a finite category, it is enough to show that there exists some $M \in \N$ such that, for all $i \geq M$, $\LL_i$ satisfies axiom (I') in Remark \ref{rmktransp}: $\varepsilon_i(S_i) \in \Syl_p(\Aut_{\LL_i}(S_i))$ or, equivalently, the group $\Out_{\FF_i}(S_i)$ has trivial Sylow $p$-subgroup.

Fix a set $\nn \subseteq \Aut_{\wL}(S)$ of representatives of the elements of $\Out_{\FF}(S) = \Aut_{\wL}(S)/S$. Then there exists some $M_b \in \N$ such that, for all $i \geq M_b$, $(S_i)^{\star} = S$ and $\nn \subseteq \Aut_{\wL}(S_i)$ (by abuse of notation we consider $\nn$ as the restriction of its elements to $S_i$). Furthermore, by Proposition \ref{fix1} (ii) we can assume that $\nn \subseteq \Aut_{\LL_i}(S_i)$ for all $i \geq M_b$. Thus, there is a commutative diagram of group extensions
$$
\xymatrix{
\4{\varepsilon}_S(S)^{\Psi_i} \ar[r] \ar[d]_{\res} & \Aut_{\wL}(S)^{\Psi_i} \ar[r] \ar[d]^{\res} & \Out_{\FF}(S) \ar[d] \\
\varepsilon_i(S_i) \ar[r] & \Aut_{\LL_i}(S_i) \ar[r] & \Out_{\FF_i}(S_i)
}
$$
where $G^{\Psi_i} = \{g \in G \, | \, \Psi_i(g) = g\}$, for $G = \varepsilon_S(S)$ or $G = \Aut_{\LL}(S)$. Furthermore, note that the restrictions $\res \colon \4{\varepsilon}_S(S)^{\Psi_i} \to \varepsilon_i(S_i)$ and $\res \colon \Aut_{\wL}(S)^{\Psi_i} \to \Aut_{\LL_i}(S_i)$ are isomorphisms by definition. Thus, for all $i \geq M_b$ we have $\Out_{\FF_i}(S_i) \cong \Out_{\FF}(S)$, and axiom (I') follows since $\{1\} \in \Syl_p(\Out_{\FF}(S))$.

\textbf{Step 3.} Axiom (II) is satisfied: let $\varphi \in \Iso_{\LL_i}(P,Q)$, $P \lhd \4{P}\leq S_i$, and $Q \lhd \4{Q}\leq S_i$ be such that $\varphi \circ \varepsilon_i(\4{P}) \circ \varphi^{-1}\leq \varepsilon_i(\4{Q})$; then there is some $\4{\varphi} \in \Mor_{\LL_i}(\4{P}, \4{Q})$ such that $\4{\varphi} \circ \varepsilon_i(1) = \varepsilon_i(1) \circ \varphi$.

Fix some $i \geq 0$, and let $\varphi \in \Iso_{\LL_i}(P,Q)$, $P \lhd \4{P}\leq S_i$, and $Q \lhd \4{Q}\leq S_i$ be as above, and notice that in this case we have $\varepsilon_i(1) = \4{\varepsilon} \colon X \to \4{X}$, where $X = P,Q$. Since $\wL$ is a transporter system, there is some $\4{\varphi} \in \Mor_{\wL}(\4{P}, \4{Q})$ such that $\4{\varphi} \circ \varepsilon_i(1) = \varepsilon_i(1) \circ \varphi$. Applying $\Psi_i$ to this equation we get
$$
\4{\varphi} \circ \varepsilon_i(1) = \varepsilon_i(1) \circ \varphi = \Psi_i(\varepsilon_i(1) \circ \varphi) = \Psi_i(\4{\varphi}) \circ \varepsilon_i(1),
$$
and thus Lemma \ref{epimono} implies that $\Psi_i(\4{\varphi}) = \4{\varphi}$.
\end{proof}

Again, we may assume that $M_b = 0$ for simplicity.

\begin{cor}\label{fix3}

For all $i \geq 0$, the fusion system $\FF_i$ is $\Ob(\LL_i)$-generated and $\Ob(\LL_i)$-saturated.

\end{cor}

\begin{proof}

The fusion system $\FF_i$ is $\Ob(\LL_i)$-generated by definition, and the $\Ob(\LL_i)$-saturation follows by \cite[Proposition 3.6]{OV}.
\end{proof}


\subsection{Existence of finite approximations}\label{Ssapprox}

We are ready to prove that every $p$-local compact group has an approximation by $p$-local finite groups. In fact, we prove something stronger: every fine unstable Adams operation (in the sense of \ref{uAo}) defines an approximation by $p$-local finite groups.

\begin{hyp}\label{hyp2}

For the rest of this section, let $\g = \ploc$ be a $p$-local compact group, and let $((-)^{\bullet}_{\FF}, (-)^{\bullet}_{\LL})$ be the finite retraction pair of \ref{expl1}. Let also $\wL$ be the associated telescopic linking system, and let $\Psi$ be a fine unstable Adams operation on $\wL$ (that is, $\Psi$ is an unstable Adams operation whose degree $\zeta \neq 1$ is congruent to $1$ modulo $p$). Let also $\{\Psi_i\}_{i \geq 0}$ and $\{S_i\}_{i \geq 0}$ be as defined in \ref{hyp1}, and let $\{(S_i, \FF_i, \LL_i)\}_{i \geq 0}$ be the associated family of finite transporter systems defined in \ref{Li}. Finally, for all $i \geq 0$ let
\begin{enumerate}[(a)]

\item $\Gamma_i =\{P\leq S_i \, | \, P^{\bullet} \notin \Ob(\LL^{\bullet}) \mbox{ and } P^{\bullet} \cap S_i = P\}$; and

\item $\Omega_i = \{R \in \Ob(\FF^{\bullet})\setminus \Ob(\LL^{\bullet}) \, | \, (R \cap S_i)^{\bullet} = R\}$.

\end{enumerate}
Note that $\Psi_i$ is a fine unstable Adams operation for all $i$. Also, for each $i \geq 0$ the sets $\Gamma_i$ and $\Omega_i$ are in one-to-one correspondence with each other for all $i \geq 0$. The bijection is given in one direction by $P \mapsto P^{\bullet}$, and in the reverse direction by $R \mapsto R \cap S_i$. Also note that $\Gamma_i \cap \Ob(\LL_i) = \emptyset$ for all $i \geq 0$.

\end{hyp}

\begin{prop}\label{fix2-1}

For all $i \geq 0$, $\LL_i$ is a linking system, and $P$ is $\FF_i$-centric for each $P \in \Ob(\LL_i)$.

\end{prop}

\begin{proof}

By definition of $\LL_i$, for all $P \in \Ob(\LL_i)$ we have $\Ker(\Aut_{\wL}(P) \to \Aut_{\FF}(P)) = \varepsilon_P(Z(P^{\bullet}))$ by property (1) in \ref{expl1}, and thus $\Ker(\Aut_{\LL_i}(P) \to \Aut_{\FF_i}(P)) = \varepsilon_P(Z(P^{\bullet}) \cap S_i) = \varepsilon_i(Z(P))$. Furthermore, we have $\varepsilon_i(C_{S_i}(P))\leq \Ker(\Aut_{\LL_i}(P) \to \Aut_{\FF_i}(P))$ by axiom (C) of transporter systems on $\LL_i$, and the statement follows immediately.
\end{proof}

\begin{prop}\label{fix4}

Let $i \geq 0$ and let $P\leq S_i$ be such that $P \lneqq P^{\bullet} \cap S_i$. Then
$$
\Out_{S_i}(P) \cap O_p(\Out_{\FF_i}(P)) \neq 1.
$$
In particular, $P$ is not $\FF_i$-centric $\FF_i$-radical.

\end{prop}

\begin{proof}

Suppose $P$ is $\FF_i$-centric. Since $P \notin \Gamma_i$, we have $P \lneqq P^{\bullet} \cap S_i \defin Q$. Notice that the functor $\bullf$ induces an inclusion $\Aut_{\FF_i}(P)\leq \Aut_{\FF_i}(Q)$. Consider the subgroup $A = \{c_x \in \Aut_{\FF_i}(P) \, | \, x \in N_Q(P)\}$. Via the above inclusion, we have $A = \Aut_{\FF_i}(P) \cap \Inn(Q)$. Since $P \lneqq Q$, it follows that $P \lneqq N_Q(P)$, and hence $\Inn(P) \lneqq A$, since $P$ is $\FF_i$-centric by hypothesis. The group $\Aut_{\FF_i}(P)$, seen as a subgroup of $\Aut_{\FF_i}(Q)$, normalizes $\Inn(Q)$, and thus $A \lhd \Aut_{\FF_i}(P)$, and
$$
\{1\} \neq A/\Inn(P)\leq \Out_{S_i}(P) \cap O_p(\Out_{\FF_i}(P)).
$$
This finishes the proof.
\end{proof}

\begin{prop}\label{fix5}

Assume Hypothesis \ref{hyp2}. Then, there exists some $M \in \N$ such that, for all $i \geq M$, the following holds: if $P \in \Gamma_i$, then either $P$ is not $\FF_i$-centric or $P$ is $\FF_i$-conjugate to some $Q\leq S_i$ such that
$$
\Out_{S_i}(Q) \cap O_p(\Out_{\FF_i}(Q)) \neq 1.
$$

\end{prop}

\begin{proof}

We start with some general observations, after which we deduce a certain condition ($\ddagger$) which will imply the statement. The rest of the proof consists of a series of steps to show that ($\ddagger$) holds.

By Proposition \ref{3.2BLO3} (i), the set $\Ob(\FF^{\bullet})$ contains finitely many $\FF$-conjugacy classes, and the same applies to $\Ob(\FF^{\bullet}) \setminus \Ob(\LL^{\bullet})$. Let $\pp = \{X_1, \ldots, X_n\}$ be a set of representatives of the $\FF$-conjugacy classes in $\Ob(\FF^{\bullet}) \setminus \Ob(\LL^{\bullet})$. By Proposition \ref{fix1} (i), there exists some $M' \in \N$ such that $\pp \subseteq \Omega_i$ for all $i \geq M'$, and we can assume that $M' = 0$ without loss of generality.

For each $X \in \pp$, fix also a set $\hh_X = \{Y_1, \ldots, Y_m\}$ of representatives of the $T$-conjugacy classes in $X^{\FF}$. Again, there exists some $M'' \in \N$ such that $\hh_X \subseteq \Omega_i$ for all $i \geq M''$, and once more we can assume $M'' = 0$ for simplicity.

Let $i \geq 0$, and let $Q\leq S_i$. Then, $T_{Q^{\bullet}} \cap Q \lhd Q$, and every automorphism in $\FF_i$ of $Q$ restricts to an automorphism of $T_{Q^{\bullet}} \cap Q$ (by the properties of $(-)^{\bullet}$, see \ref{3.2BLO3}). Set
$$
K_Q = \Ker(\Aut_{\FF_i}(Q) \Right2{} \Aut_{\FF_i}(T_{Q^{\bullet}} \cap Q) \times \Aut(Q/(T_{Q^{\bullet}} \cap Q))),
$$
which is a normal $p$-subgroup of $\Aut_{\FF_i}(Q)$ by Lemma \ref{Kpgp}. Since $\Gamma_i$ and $\Omega_i$ are in one-to-one correspondence (see \ref{hyp2}), in order to prove the statement it is enough to prove the following, slightly stronger statement for each $X \in \pp$ and each $Y \in \hh_X$.
\begin{itemize}

\item[($\ddagger$)] There exists some $M_Y \in \N$ such that the following holds for all $i \geq M_Y$: if $R \in Y^T \cap \Omega_i$, then either $R \cap S_i$ is not $\FF_i$-centric, or $R \cap S_i$ is $\FF_i$-conjugate to some $Q\leq S_i$ such that $K_Q \cap \Aut_S(Q)$ contains some element which is not in $\Inn(Q)$.

\end{itemize}

Indeed, suppose that ($\ddagger$) holds for all $X \in \pp$ and all $Y \in \hh_X$. We claim that the statement follows with $M = \max\{M_Y \, | \, X \in \pp \mbox{ and } Y \in \hh_X\}$. To prove this, let $P \in \Gamma_i$, and let $R = P^{\bullet} \in \Omega_i$. Then, there exist some $X \in \pp$ and some $Y \in \hh_X$ such that $R \in Y^T \cap \Omega_i$. Thus, ($\ddagger$) applies to $Y$, and it follows that either $P = R \cap S_i$ is not $\FF_i$-centric, or $P = R \cap S_i$ is $\FF_i$-conjugate to some $Q\leq S_i$ such that $K_Q \cap \Aut_S(Q)$ contains some element which is not in $\Inn(Q)$, in which case we have
$$
\Out_{S_i}(Q) \cap O_p(\Out_{\FF_i}(Q)) \neq 1.
$$

For the rest of the proof, fix $X \in \pp$ and $\hh_X$ as above. Since this proof is rather long, we have divided it into several steps, for the reader's convenience. We also include a brief summary of the steps in the proof. In Step 1, we give a general tool to deduce that ($\ddagger$) holds in some cases. In Step 2 we show that we may reduce to prove that ($\ddagger$) holds for all $Y \in \hh_X$ with $T_Y$ fully $\FF$-normalized. In Step 3 we justify the reduction to assume that $A = T_X$ is normal in $\FF$. In Step 4 we show some properties regarding the quotient $\g/A = (S/A, \FF/A, \LL/A)$. In Step 5 we introduce a certain subgroup $Z_Y$ for each $Y \in \hh_X$, related to $C_{T/A}(Y/A)$, and prove some of its properties. In Step 6, we show that ($\ddagger$) holds for all $Y \in \hh_X$ such that $Z_Y \not\leq Y$. In Step 7 we prove that we may reduce to prove ($\ddagger$) for all $Y \in \hh_X$ such that $Z_Y$ is fully $\FF$-normalized. In Step 8 we show that ($\ddagger$) holds for all $Y \in \hh_X$ such that $C_S(Y) \not\leq Y$ (in particular, this applies to $Y \in \hh_X$ such that $C_S(Y) \not\leq Y$ and such that $Z_Y$ is fully $\FF$-normalized). In step 9 we prove a technical property necessary for the last step of the proof. Finally, in Step 10 we show that ($\ddagger$) holds for all $Y \in \hh_X$ such that $Z_Y$ is fully $\FF$-normalized.

\textbf{Step 1.} Let $Y, Y' \in \hh_X$, and suppose that ($\ddagger$) holds for $Y$. Suppose in addition that there exists some $f \in \Hom_{\FF}(Y'T, YT)$ such that $f(Y') = Y$. Then ($\ddagger$) holds for $Y'$.

Let $f \in \Hom_{\FF}(Y'T, YT)$ as above. By Alperin's Fusion Theorem \cite[Theorem 3.6]{BLO3}, there exist subgroups $A_0 = Y'T, A_1, \ldots, A_n = YT\leq S$, objects $B_1, \ldots, B_n \in \Ob(\LL^{\bullet})$, and automorphisms $\phi_k \in \Aut_{\wL}(B_k)$, for $k = 1, \ldots, n$, such that
$$
A_{k-1}, A_k\leq B_k \qquad \mbox{and} \qquad \rho(\phi_k)(A_{k-1}) = A_k
$$
for each $k = 1, \ldots, n$, and such that $f = \rho(\phi_n) \circ \ldots \rho(\phi_1)$. By Proposition \ref{fix1} (i) and (ii), there exists some $M_1 \in \N$ such that, for all $i \geq M_1$ and all $k = 1, \ldots, n$,
$$
B_k \cap S_i \in \Ob(\LL_i) \qquad \mbox{and} \qquad (\phi_k)|_{B_k \cap S_i} \in \Aut_{\LL_i}(B_k \cap S_i).
$$
Moreover, since $T\leq A_0$, it follows that $T\leq A_k, B_k$ for all $k = 1, \ldots, n$. Note that for each $k$, there exists some $t_k \in T$ such that $t_k A_k t_k^{-1} \in \hh_X$ (in particular, $t_0 = t_n = 1$). Since $T\leq B_k$ for each $k$, we may replace $\phi_k$ by $\varepsilon(t_k) \circ \phi_k \circ \varepsilon(t_{k-1})$, and this way we may assume that $A_k \in \hh_X$ for each $k$.

Suppose now that ($\ddagger$) holds for $Y$. That is, there exists some $M_Y \in \N$ (we may choose $M_Y \geq M_1$) such that, for all $i \geq M_Y$, the following holds: if $R \in Y^T \cap \Omega_i$, then $R \cap S_i$ satisfies the conclusion of ($\ddagger$). Let $i \geq M_Y$, and let $K \in (Y')^T \cap \Omega_i$, and let $H = f(K)$. Since $f \colon Y'T \to YT$, it follows that $H \in Y^T$. Moreover, $f(K \cap S_i) = H \cap S_i$ (since each $\phi_k$ above is fixed by $\Psi_i$). Thus, $K \cap S_i$ is $\FF_i$-conjugate to $H \cap S_i$, and ($\ddagger$) holds for $Y'$.

\textbf{Step 2.} We show that if ($\ddagger$) holds for all $Y \in \hh_X$ such that $T_Y$ is fully $\FF$-normalized, then ($\ddagger$) holds for all $Y \in \hh_X$.

Indeed, let $Y' \in \hh_X$, and suppose that $T_{Y'}$ is not fully $\FF$-normalized. Then, there exists some $\gamma \in \Hom_{\FF}(N_S(T_{Y'}), S)$ such that $\gamma(T_{Y'})$ is fully $\FF$-normalized. Since $Y'\leq N_S(T_{Y'})$, we may define $Y = \gamma(Y')$. Since $\hh_X$ contains representatives of all the $T$-conjugacy classes in $X^{\FF}$, there is some $t \in T$ such that $tYt^{-1} \in \hh_X$. Thus, upon replacing $\gamma$ by $c_t \circ \gamma$, we may assume that $Y \in \hh_X$. Note that $T_{Y} = \gamma(T_{Y'})$. Notice also that $T\leq N_S(T_{Y'})$, and thus $\gamma$ restricts to some $f \in \Hom_{\FF}(Y'T, YT)$ such that $Y = \gamma(Y')$. The claim follows by Step 1.

\textbf{Step 3.} Let $\mm = \{A = T_Y\leq T \, | \, Y \in \hh_X \mbox{ and $T_Y$ is fully $\FF$-normalized}\}$. Since $\hh_X$ is finite, so is $\mm$. For each $A \in \mm$, let $\hh_{X,A} = \{U \in \hh_X \, | \, T_U = A\} \subseteq \hh_X$. To prove the statement of \ref{fix5}, it is enough to prove that ($\ddagger$) holds for all $U \in \hh_{X,A}$, for a fixed $A \in \mm$ at a time. Fix $A \in \mm$, and let $\hh_{X,A}$ be as above. The main goal of this step is to justify the reduction to the case where $A$ is normal in $\FF$.

Since $A$ is fully $\FF$-normalized, to prove that ($\ddagger$) holds for all $U \in \hh_{X,A}$,  we can reduce to work on the normalizer $p$-local compact group $N_{\g}(A) = (N_S(A), N_{\FF}(A), N_{\LL}(A))$ defined in \ref{rmknorm}, instead of the whole $\g = \ploc$. Set for short $N = N_S(A)$, $\ee = N_{\FF}(A)$, and $\TT = N_{\LL}(A)$, and note that $T\leq N$ since $A\leq T$. Let $((-)^{\bullet}_{\FF}, (-)^{\bullet}_{\LL})$ be the finite retraction pair for $\g$ fixed in \ref{hyp2}, and let $((-)^{\star}_{\ee}, (-)^{\star}_{\TT})$ be the finite retraction pair for $N_{\g}(A)$ described in \ref{expl3}. Recall that $P^{\star} = P^{\bullet}$ for all $P\leq N$. Let also $\wL$ and $\4{\TT}$ be the corresponding associated telescopic transporter systems.

We start by stating and proving several general properties.
\begin{itemize}

\item[(3-a)] Let $U \in \hh_{X,A}$ and let $R \in U^T$. Then, $T_R = A$, $N_S(T_R) = N_S(A)$, and $N_S(R)\leq N_S(A)$. Moreover, every automorphism of $R$ preserves $T_R$, and thus $\Aut_{\FF}(R) = \Aut_{\ee}(R)$.

\item[(3-b)] Since $A\leq T$ is a subtorus, it follows that $\Psi_i(A) = A$ for all $i \geq 0$. Thus, each $\Psi_i$ restricts to a fine unstable Adams operation (see \ref{uAo}) on $N_{\g}(A)$, which extends to $\4{\TT}$ by \ref{extendL2}. Let us denote by $\Psi_i$ the resulting unstable Adams operation.

\end{itemize}
Consider the family of transporter systems $\{(N_i, \ee_i, \TT_i)\}_{i \geq 0}$ associated to $(N, \ee, \4{\TT})$ in \ref{Li}. As pointed out in \ref{expl3}, $\4{\TT}$ is not a subcategory of $\wL$, and thus it is hard to compare the fusion systems $\ee_i$ and $\FF_i$. However, if we restrict to the full subcategory $\4{\TT}_{\geq A} \subseteq \4{\TT}$ of subgroups $P\leq N$ such that $A\leq P^{\star} = P^{\bullet}$, then $\4{\TT}_{\geq A}$ is a subcategory of $\wL$, by \ref{expl3} (a). In particular we have the following.
\begin{itemize}

\item[(3-c)] For all $i \geq 0$, there is an inclusion $\4{\TT}_{\geq A} \cap \TT_i \subseteq \LL_i$.

\end{itemize}
Finally, note that $A^{\bullet} = A$, since $A$ is the maximal torus of $U$, for some $U \in \hh_{X,A}$, and $U = U^{\bullet}$ by hypothesis. The following holds.
\begin{itemize}

\item[(3-d)] For each $i \geq 0$, let $A_i = A \cap S_i$. By \ref{fix1} (i), there exists some $M_2 \in \N$ such that $(A_i)^{\bullet} = A$ for all $i \geq M_2$. For simplicity we may assume that $M_2 = 0$.

\item[(3-e)] For each $U \in \hh_{X,A}$ and each $R \in U^T \cap \Omega_i$, we have $R \cap N_i = R \cap S_i$ and $\Aut_{N_i}(R \cap S_i) = \Aut_{S_i}(R \cap S_i)$. This follows since $N_S(R \cap S_i)\leq N_S(R)\leq N_S(A) = N$ by (3-a).

\item[(3-f)] Let $H,K\leq N_i$ be such that $A\leq H^{\bullet}, K^{\bullet}$. Then, for all $f \in \Hom_{\ee}(H,K)$, we have $(f)^{\star}_{\ee} = (f)^{\bullet}_{\FF}$. Similarly, if $H, K \in \Ob(\4{\TT}_{\geq A})$ and $\varphi \in \Mor_{\4{\TT}}(H,K)$, then $H, K \in \Ob(\wL)$, and $(\varphi)^{\star}_{\4{\TT}} = (\varphi)^{\bullet}_{\wL}$. The first part follows by definition of $(-)^{\star}_{\ee}$ in \ref{expl3}, and the second part follows by \ref{expl3} (c).

\end{itemize}
In fact, we deduce more. Let $W\leq N_i$ be such that $A\leq W^{\bullet}$. Then, the following holds.
\begin{itemize}

\item[(3-g)] $\Hom_{\ee_i}(W, N_i) \subseteq \Hom_{\FF_i}(W, N_i)$.

\end{itemize}
Let $f \in \Hom_{\ee_i}(W, N_i)$. Since $\ee_i$ is $\Ob(\TT_i)$-generated, there exist objects $H_k, H_k' \in \Ob(\TT_i)$ and morphisms $\gamma_k \in \Hom_{\TT_i}(H_k, H'_k)$, for $k = 1, \ldots, n$, and such that, upon taking the necessary restrictions,
$$
f = \4{\rho}_i(\gamma_n) \circ \ldots \circ \4{\rho}_i(\gamma_1).
$$
Notice that $(A_i)^{\bullet} = A\leq W^{\bullet}$. Thus, $A\leq W^{\bullet}\leq (H_1)^{\bullet}$, and $A = \4{\rho}(\gamma_1)(A)\leq (H_1')^{\bullet}$ since $\gamma_1 \in \Mor(\4{\TT})$. Inductively, it follows that $A\leq (H_k)^{\bullet}, (H'_k)^{\bullet}$ for all $k$, and thus $H_k, H_k' \in \Ob(\4{\TT}_{\geq A})$. By (3-c) above, for all $k = 1, \ldots, n$ we have $H_k, H'_k \in \LL_i$, and $\gamma_k \in \Mor(\LL_i)$. Thus $f \in \Mor(\FF_i)$, and (3-g) follows. In particular, if $U \in \hh_{X,A}$, and $P \in U^{\ee} \cap \Omega_i$, then
$$
\Hom_{\ee_i}(P \cap S_i, N_i) \subseteq \Hom_{\FF_i}(P \cap S_i, N_i).
$$

Let $U \in \hh_{X,A}$. For each $P \in U^{\ee} \cap \Omega_i$, note that $A \cap S_i \lhd P \cap S_i$, and every automorphism of $P \cap S_i$ in $\FF_i$ restricts to an automorphism of $A \cap S_i$. Let $\3{P} = P/A$. Note that $P = (P \cap S_i)A$ (since $P \in \Omega_i$), and thus we have $(P \cap S_i)/(A \cap S_i) \cong \3{P}$. Set for short $Q = P \cap S_i$, and let
$$
K_Q = \Ker(\Aut_{\FF_i}(Q) \Right2{} \Aut_{\FF_i}(A \cap S_i) \times \Aut(\3{P}));
$$
$$
K'_Q = \Ker(\Aut_{\ee_i}(Q) \Right2{} \Aut_{\ee_i}(A \cap S_i) \times \Aut(\3{P})).
$$
By Lemma \ref{Kpgp}, we have $K_Q\leq O_p(\Aut_{\FF_i}(Q))$ and $K'_Q\leq O_p(\Aut_{\ee_i}(Q))$. Moreover, by (3-g) there is an inclusion $\Aut_{\ee_i}(Q)\leq \Aut_{\FF_i}(Q)$, and it follows that $K'_Q\leq K_Q$. We claim that, in order to prove that ($\ddagger$) holds for $U \in \hh_{X,A}$, it is enough to prove that the following version of ($\ddagger$) in terms of $\ee_i$ holds.
\begin{itemize}

\item[(3-h)] There exists some $M_U \in \N$ such that, for all $i \geq M_U$ and all $R \in U^T \cap \Omega_i$, either $R \cap S_i$ is not $\ee_i$-centric, or $R \cap S_i$ is $\ee_i$-conjugate to some $Q$ such that $K'_{Q} \cap \Aut_{N_i}(Q)$ contains some element that is not in $\Inn(Q)$.

\end{itemize}
Indeed, suppose that (3-h) holds for $U$, and let $i \geq M_U$ and $R \in U^T \cap \Omega_i$. Clearly, if $R \cap S_i$ is not $\ee_i$-centric, then it is not $\FF_i$-centric by (3-g). Suppose that $R \cap S_i$ is $\ee_i$-centric. Then, $R \cap S_i$ is $\ee_i$-conjugate to some $Q$ such that $K'_Q \cap \Aut_{N_i}(Q)$ contains some element that is not in $\Inn(Q)$. Then, $R \cap S_i$ is $\FF_i$-conjugate to $Q$ by (3-g), and ($\ddagger$) holds for $U$ since $K'_Q \cap \Aut_{N_i}(Q)\leq K_Q \cap \Aut_{S_i}(Q)$.

The above just shows that, in order to prove ($\ddagger$) for $U \in \hh_{X, A}$ with respect to $\{(S_i, \FF_i, \LL_i)\}_{i \geq 0}$, it is enough to prove that the conclusion ($\ddagger$) holds for $U$ with respect to $\{(N_i, \ee_i, \TT_i)\}_{i \geq 0}$. For the rest of the proof, we assume that $A = T_X$ is normal in $\FF$, so $\hh_{X, A} = \hh_X$ (and $T_R = A$ for all $R \in X^{\FF}$). Moreover, notice that we are only concerned about subgroups $H\leq S$ such that $A\leq H^{\bullet}$, and for these subgroups there is no difference between the finite retraction pairs $((-)^{\bullet}_{\FF}, (-)^{\bullet}_{\LL})$ and $((-)^{\star}_{\ee}, (-)^{\star}_{\TT})$, by (3-f). Thus, we may assume also that we are still working with $((-)^{\bullet}_{\FF}, (-)^{\bullet}_{\LL})$.

To finish this step, we prove the following.
\begin{itemize}

\item[(3-i)] If $A = T$ then ($\ddagger$) holds for all $Y \in \hh_X$.

\end{itemize}
Indeed, in this case we have $Y^T = \{Y\}$ for each $Y \in \hh_X$. Since $\hh_X$ is finite, it is not hard to check in this case that in fact there exists some $M_X \in \N$ such that, for all $i \geq M_X$ and all $Y \in \hh_X$, the subgroup $Y \cap S_i$ is not $\FF_i$-centric.

\textbf{Step 4.} Since $A = T_X$ is normal in $\FF$, consider the quotient $(S/A, \FF/A, \wL/A)$. By Lemma \ref{quotient22}, $\FF/A$ is a saturated fusion system over $S/A$, and $\wL/A$ is a transporter system associated to $\FF/A$ which contains all the centric subgroups of $\FF/A$. In this step we prove several properties relating $(S, \FF, \wL)$ to $(S/A, \FF/A, \wL/A)$, which we will need in later steps.

Set for short $\3{S} = S/A$, $\3{\FF} = \FF/A$ and $\3{\TT} = \wL/A$. Let also $\wL_{\geq A} \subseteq \wL$ be the full subcategory of subgroups that contain $A$, and let $\tau \colon \wL_{\geq A} \to \3{\TT}$ be the projection functor. By a slight abuse of notation, we also write $\tau \colon S \to \3{S}$ for the projection homomorphism. We adopt the notation $\3{P}, \3{Q}, \ldots$ to denote subgroups of $\3{S}$, $\3{f}, \3{f}', \ldots$ to denote morphisms in $\3{\FF}$, and $\3{\varphi}, \3{\psi}, \ldots$ to denote morphisms in $\3{\TT}$. In particular, $T/A = \3{T}\leq \3{S}$ denotes the maximal torus of $\3{S}$. Also, for each $P\leq S$ that contains $A$, we will write $\3{P}$ instead of $\tau(P) = P/A$, unless there is risk of confusion.

The following is easily check to hold since $A\leq T$ is normal in $\FF$.
\begin{itemize}

\item[(4-a)] Let $P, Q\leq S$ be such that $A\leq P,Q$. Then, $Q \in P^{\FF}$ if and only if $\3{Q} \in \3{P}^{\3{\FF}}$. Similarly, $Q \in P^T$ if and only if $\3{Q} \in \3{P}^{\3{T}}$, since $A\leq T$.

\end{itemize}

By property (3-b), $\Psi_i(A) = A$ for all $i \geq 0$. By (3-i) we may assume that $A \lneqq T$, and thus, for each $i \geq 0$, the unstable Adams operation $\Psi_i$ induces a fine unstable Adams (see \ref{uAo}) operation $\3{\Psi}_i$ on $\3{\TT}$. By definition, if $\varphi \in \Mor(\wL)$ is such that $\Psi_i(\varphi) = \varphi$ for some $i$, then $\3{\Psi}_i(\3{\varphi}) = \3{\varphi}$. The following is some sort of converse of this statement.
\begin{itemize}

\item[(4-b)] Let $P, Q \in \Ob(\wL)$ be such that $A\leq P, Q$, and let $\3{\varphi} \in \Mor_{\3{\TT}}(\3{P}, \3{Q})$ be such that $\3{\Psi}_i(\3{\varphi}) = \3{\varphi}$. Then, there exists some $\varphi \in \Mor_{\wL}(P,Q)$ such that $\tau(\varphi) = \3{\varphi}$ and such that $\Psi_i(\varphi) = \varphi$. Similarly, let $\3{x} \in \3{S}$ be such that $\3{\Psi}_i(\3{x}) = \3{x}$. Then there exists some $x \in S$ such that $\Psi_i(x) = x$ (i.e. $x \in S_i$) and $\tau(x) = \3{x}$.

\end{itemize}
Let $\3{\varphi}$ be as above, and let $\varphi \in \Mor_{\wL}(P,Q)$ be such that $\tau(\varphi) = \3{\varphi}$. Since $\tau(\varphi) = \3{\varphi} = \3{\Psi}_i(\3{\varphi}) = \tau(\Psi_i(\varphi))$, it follows that $\Psi_i(\varphi) = \varphi \circ \varepsilon_P(a)$, for some $a \in A$ (see Definition \ref{quotient1}). Consider the map $A \to A$ defined by $t \mapsto t^{-1} \Psi_i(t)$. Since $A$ is abelian, this turns out to be a group homomorphism, which is in fact surjective, since $\Ker(A)$ is the subgroup of fixed points of $A$, and this is a finite subgroup of $A$ for all $i \geq 0$. In particular, there exists some $t \in A$ such that $t^{-1} \Psi_i(t) = a^{-1}$, and we get
$$
\Psi_i(\varphi \circ \varepsilon_P(t)) = \Psi_i(\varphi) \circ \varepsilon_P(\Psi_i(t)) = \varphi \circ \varepsilon_P(a) \circ \varepsilon_P(\Psi_i(t)) = \varphi \circ \varepsilon_P(t).
$$
A similar argument shows that, for each $\3{x} \in \3{S}$ such that $\3{\Psi}_i(\3{x}) = \3{x}$, there is some $x \in S_i$ such that $\tau(x) = \3{x}$.

For each $i \geq 0$, let $\3{S}_i\leq \3{S}$ be the subgroup of fixed points by $\3{\Psi}_i$. By (4-b), we deduce that $S_iA/A = \3{S}_i$, and thus $S_i/(S_i \cap A) \cong \3{S}_i$. Let $U \in \hh_X$. By (4-a), $R \in U^T$ if and only if $\3{R} \in \3{U}^{\3{T}}$.
\begin{itemize}

\item[(4-c)] For all $i \geq 0$ and all $R \in U^{T}$, $R \in \Omega_i$ if and only if $\3{R}\leq \3{S}_i$.

\end{itemize}
Suppose first that $R \in \Omega_i$, and let $P = R \cap S_i$. Then, $R = PA$ (since $A = T_R$), and $\3{R} = PA/A\leq S_iA/A = \3{S}_i$. Conversely, let $R \in U^{T}$ be such that $\3{R}\leq \3{S}_i$. Then, it follows by (4-b) that $P = R \cap S_i$ contains representatives of all the elements in $\3{R}$. Since $A_i = A \cap S_i\leq R \cap S_i = P$ and $(A_i)^{\bullet} = A$, it follows that $P^{\bullet} \geq R$. The inclusion $P^{\bullet}\leq R$ follows from $P = R \cap S_i\leq R$ and the fact that $R^{\bullet} = R$. Thus $R \in \Omega_i$.

\textbf{Step 5.} For each $R \in X^{\FF}$, recall the notation $\3{R} = R/A$. Set also
$$
Z_{\3{R}} \defin C_{\3{T}}(\3{R}) \qquad \mbox{and} \qquad Z_R \defin \{t \in N_T(R) \, | \, \tau(t) \in Z_{\3{R}}\}.
$$
Note that $A\leq Z_R$ and $Z_{\3{R}} = Z_R/A$. The main goal of this step is to show the following: if $Z_U$ is fully $\FF$-normalized for some $U \in \hh_X$, then $Z_{\3{U}}$ is fully $\3{\FF}$-normalized.

Indeed, let $U \in \hh_X$ be such that $Z_U$ is fully $\FF$-normalized. Note that $T\leq N_S(Z_U)$, and
$$
\3{T}\leq N_{\3{S}}(Z_{\3{U}}) = N_S(Z_U)/A.
$$
Suppose that $Z_{\3{U}}$ is not fully $\3{\FF}$-normalized, and let $\3{\gamma} \in \Hom_{\3{\FF}}(N_{\3{S}}(Z_{\3{U}}), \3{S})$ be such that $\3{H} = \3{\gamma}(Z_{\3{U}})$ is fully $\3{\FF}$-normalized, and $|N_{\3{S}}(Z_{\3{U}})| < |N_{\3{S}}(\3{\gamma}(Z_{\3{U}}))|$. Then $\3{\gamma}$ lifts to a map $\gamma \in \Hom_{\FF}(N_S(Z_U), S)$ such that $|N_S(Z_U)| = |\gamma(N_S(Z_U))| < |N_S(\gamma(Z_U))|$, since
$$
\3{T}\leq N_{\3{S}}(\3{\gamma}(Z_{\3{U}})) = N_S(\gamma(Z_U))/A,
$$
and this contradicts the maximality of $|N_S(Z_U)|$.

\textbf{Step 6.}  We show that ($\ddagger$) holds for all $U \in \hh_X$ such that $Z_U \not\leq U$.

Let $U \in \hh_X$ be such that $Z_U \not\leq U$, and let $a \in Z_U \setminus U$. Then, by Lemma \ref{SiS} there exists some $M_U \in \N$ such that $a \in S_i$ for all $i \geq M_U$. Fix some $i \geq M_U$, and let $R \in U^T \cap \Omega_i$. Since $R$ is $T$-conjugate to $U$, we have $Z_R = Z_U$. Let
$$
K_{R \cap S_i} = \Ker(\Aut_{\FF_i}(R \cap S_i) \Right2{} \Aut_{\FF_i}(T_R \cap S_i) \times \Aut(\3{R})).
$$
By Lemma \ref{Kpgp} we know that $K_{R \cap S_i}\leq O_p(\Aut_{\FF_i}(R \cap S_i))$.

The element $a \in Z_U \cap S_i$ fixed above satisfies $a \in Z_R \cap S_i$, since $Z_R = Z_U$, and in particular $a \in N_{S_i}(R \cap S_i)$. Moreover, $c_a \in K_{R \cap S_i} \cap \Aut_{S_i}(R \cap S_i)$, since $a \in T \cap S_i$. In particular, if $R \cap S_i$ is $\FF_i$-centric then $c_a \notin \Inn(R \cap S_i)$, and this shows that ($\ddagger$) holds for $U$.

\textbf{Step 7.} Suppose ($\ddagger$) holds for all $U \in \hh_X$ such that $Z_U$ is fully $\FF$-normalized. Then we claim that ($\ddagger$) holds for all $U \in \hh_X$.

Let $U \in \hh_X$ be such that $Z_U$ is not fully $\FF$-normalized. Then there exists some $\gamma \in \Hom_{\FF}(N_S(Z_U), S)$ such that $\gamma(Z_U)$ is fully $\FF$-normalized. Since $U\leq N_S(Z_U)$, we may set $V = \gamma(U)$. Moreover, since $\hh_X$ contains representatives of all the $T$-conjugacy classes in $X^{\FF}$, we may assume that $V \in \hh_X$. Finally, note that $T\leq N_S(Z_U)$, since $Z_U\leq T$. Thus, $Z_V = \gamma(Z_U)$.

Suppose now that ($\ddagger$) holds for all $V \in \hh_X$ such that $Z_V$ is fully $\FF$-normalized, and let $U \in \hh_X$. By the above discussion, there exists some $f \in \Hom_{\FF}(UT, S)$ such that $V = f(U) \in \hh_X$ is such that $Z_V$ is fully $\FF$-normalized, and the claim follows by Step 1.

\textbf{Step 8.} We show that ($\ddagger$) holds for all $U \in \hh_X$ such that $C_S(U) \not\leq U$. 

By Step 6, if $Z_U \not\leq U$ then ($\ddagger$) holds for $U$. Thus we may assume that $Z_U\leq U$. Also, since $C_S(U) \not\leq U$, it follows that $C_{\3{S}}(\3{U}) \not\leq \3{U}$. Notice that $C_{\3{T}}(\3{U}) = Z_{\3{U}}\leq \3{U}$ and $\3{U}$ is finite, which implies that $C_{\3{S}}(\3{U})$ is a finite group. Thus, there exists some $M_U \in \N$ such that $C_{\3{S}}(\3{U})\leq \3{S}_i$. Moreover, since $U \in \Omega_i$ by assumption, we have $\3{U}\leq \3{S}_i$ by (4-c).

Let $i \geq M_U$ and let $R \in U^T \cap \Omega_i$. Then $C_{\3{T}}(\3{R}) = Z_{\3{U}}$. By (4-a) we have $\3{R} \in \3{U}^{\3{T}}$, and by (4-c) we know that $\3{R}\leq \3{S}_i$. Thus, by Lemma \ref{invar1} (i), together with Proposition \ref{fix1} (iii), we have
$$
\3{x}^{-1} \cdot \3{\Psi}_i(\3{x}) \in C_{\3{T}}(\3{R}) = Z_{\3{R}} = Z_{\3{U}}
$$
for some $\3{x} \in N_{\3{T}}(\3{R}, \3{U})$. Fix such $\3{x} \in N_{\3{T}}(\3{R}, \3{U})$, and note that $\3{x}$ conjugates $C_{\3{S}}(\3{R})$ to $C_{\3{S}}(\3{U})$, and in particular $C_{\3{S}}(\3{R}) \not\leq \3{R}$. Moreover, since $Z_{\3{U}} = C_{\3{T}}(\3{R})\leq \3{R}$, it follows that
$$
\3{x}^{-1} \cdot \3{\Psi}_i(\3{x}) \in C_{\3{T}}(\3{R})\leq \3{R} \cap \3{T}\leq C_{\3{T}}(C_{\3{S}}(\3{R})).
$$
Thus, by \ref{invar1} (i) we deduce that $C_{\3{S}}(\3{R})\leq \3{S}_i$.

Set $X_R = C_S(R)A$. Then, by (4-b) $X_R \cap S_i$ contains representatives of all the elements in $C_S(R)A/A\leq C_{\3{S}}(\3{R})\leq \3{S}_i$. Note that $C_S(R)A \not\leq R$, and thus $X_R \cap S_i$ contains elements which are not in $R \cap S_i$. If $(X_R \cap S_i)\setminus (R \cap S_i)$ contains some element of $C_S(R)$, then clearly $R \cap S_i$ is not $\FF_i$-centric. Suppose then that $(X_R \cap S_i) \cap C_S(R)\leq R \cap S_i$, and let
$$
K_{R \cap S_i} = \Ker(\Aut_{\FF_i}(R \cap S_i) \Right2{} \Aut_{\FF_i}(T_R \cap S_i) \times \Aut(\3{R})),
$$
which is a normal $p$-subgroup of $\Aut_{\FF_i}(R \cap S_i)$ by Lemma \ref{Kpgp}. If $R \cap S_i$ is $\FF_i$-centric, then, for every $x \in (X_R \cap S_i)\setminus(R \cap S_i)$, we have $c_x \in K_{R \cap S_i} \cap \Aut_{S_i}(R \cap S_i)$ and $c_x \notin \Inn(R \cap S_i)$, and ($\ddagger$) holds for $U$.

\textbf{Step 9.} Recall the notation $\3{\TT} = \wL/A$, introduced in Step 4, which is a transporter system associated to $\3{\FF}$. Fix $U \in \hh_X$ such that $Z_U$ is fully $\FF$-normalized. Then $Z_{\3{U}}$ is fully $\3{\FF}$-normalized by Step 5, and in particular it is fully $\3{\FF}$-centralized. Thus we may consider the centralizer fusion system $\3{\ee} = C_{\3{\FF}}(Z_{\3{U}})$ over $\3{Z} = C_{\3{S}}(Z_{\3{U}})$, and note that $\3{T}\leq \3{Z}$. Since $\3{\FF}$ is saturated, it follows that $\3{\ee}$ is also saturated. The main goal of this step is to prove the following property.
\begin{itemize}

\item[(9-a)] Suppose that $Z_U\leq U$. Then, there exists $V \in \hh_X$ such that $C_S(V) \not\leq V$ and such that $\3{V} \in \3{U}^{\3{\ee}}$.

\end{itemize}

First, note that the following holds.
\begin{itemize}

\item[(9-b)] For each $\3{R} \in \3{U}^{\3{T}}$, we have $Z_{\3{R}} = C_{\3{T}}(\3{R}) = Z_{\3{U}}$.

\item[(9-c)] Since $Z_{\3{U}}\leq \3{T}$ is abelian, every $\3{\ee}$-centric subgroup of $\3{Z}$ must contain $Z_{\3{U}}$. Thus, by Lemma \ref{centricNFKA}, every $\3{\ee}$-centric subgroup of $\3{Z}$ is also an $\3{\FF}$-centric subgroup of $\3{S}$, and hence also an object in $\3{\TT}$, since $\3{\TT}$ contains all the $\3{\FF}$-centric subgroups of $\3{S}$.

\end{itemize}
We are ready now to prove (9-a). Let $K = \Aut_U(Z_U)$, and let $N_S^K(Z_U) = \{x \in S \, | \, c_x \in K\}$. Then,
$$
\Aut_S^K(Z_U) = K \cap \Aut_S(Z_U) = \Aut_U(Z_U) = K \cap \Aut_{\FF}(Z_U) = \Aut_{\FF}^K(Z_U).
$$
In particular, $Z_U$ is fully $\FF$-centralized (since it is fully $\FF$-normalized), and clearly $\Aut_S^K(Z_U) \in \Syl_p(\Aut_{\FF}^K(Z_U))$. By \cite[Lemma 2.2]{BLO6} it follows that $Z_U$ is fully $K$-normalized in $\FF$ (see Section \ref{Squotient}). Thus, the fusion system $N_{\FF}^K(Z_U)$ over $N_S^K(Z_U)$ (as defined in \ref{definorm}) is saturated.

Note that $U, T\leq N_S^K(Z_U) = U C_S(Z_U)$. Since $U$ is not $\FF$-centric and $Z_U\leq U$, it follows that $U$ is not $N_{\FF}^K(Z_U)$-centric by Lemma \ref{centricNFKA}. Let $V$ be $N_{\FF}^K(Z_U)$-conjugate to $U$ and such that $C_{N_S^K(Z_U)}(V) \not\leq V$. Since $T\leq N_S^K(Z_U)$, we may choose $V \in \hh_X$. Then,
$$
C_{N_S^K(Z_U)}(V)\leq C_S(V) \not\leq V.
$$
Also, if $f \in \Hom_{N_{\FF}^K(Z_U)}(U,V)$, then by definition $f|_{Z_U} \in \Aut_U(Z_U)$. Hence, the induced morphism $\3{f} \in \Hom_{\3{\FF}}(\3{U}, \3{V})$ satisfies $\3{f}|_{Z_{\3{U}}} = \Id$, and thus $\3{V} \in \3{U}^{\3{\ee}}$. Note that in particular $Z_U\leq Z_V$, although this may not be an equality.

\textbf{Step 10.} We show that ($\ddagger$) holds for all $U \in \hh_X$ such that $Z_U$ is fully $\FF$-normalized.

Fix $U \in \hh_X$ such that $Z_U$ is fully $\FF$-normalized. By Step 6, we may assume that $Z_U\leq U$, and thus $Z_{\3{U}}\leq \3{U}$. By Step 5, we know that $Z_{\3{U}}$ is fully $\3{\FF}$-normalized (and in particular fully $\3{\FF}$-centralized). Let $\3{Z} = C_{\3{S}}(Z_{\3{U}})$ for short, and let $\3{\ee} = C_{\3{\FF}}(Z_{\3{U}})$ be the centralizer fusion system of $Z_{\3{U}}$ in $\3{\FF}$, which is saturated by Step 9. By (9-a), there exists $V \in \hh_X$ such that $\3{V} \in \3{U}^{\3{\ee}}$ and $C_S(V) \not\leq V$.

Notice that the set $\{T_K \, | \, K \in \Ob(\LL^{\bullet})\}$ is finite since $\Ob(\LL^{\bullet})$ contains finitely many $T$-conjugacy classes. Moreover, since $K^{\bullet} = K$, it follows that $(T_K)^{\bullet} = T_K$ for all $K \in \Ob(\LL^{\bullet})$. By Proposition \ref{fix1} (i) there exists some $M_3 \in \N$ such that $(T_K \cap S_i)^{\bullet} = T_K = (T_K)^{\bullet}$ for all $i \geq M_3$ and all $T_K$ in the set above. Without loss of generality we may assume that $M_3 = 0$.

Fix $\3{f} \in \Hom_{\3{\ee}}(\3{U}, \3{V})$. By Alperin's Fusion Theorem \cite[Theorem 3.6]{BLO3}, there exist sequences of subgroups $\3{W}_0 = \3{U}, \3{W}_1, \ldots, \3{W}_n = \3{V}\leq\3{Z}$ and $\3{K}_1, \ldots, \3{K}_n \in \Ob((\3{\ee})^c) \subseteq \Ob(\3{\FF})$, and morphisms $\3{\gamma}_j \in \Aut_{\3{\ee}}(\3{K}_j)$ for each $j = 1, \ldots, n$, such that
$$
\3{W}_{j-1}, \3{W}_j\leq \3{K}_j \qquad \mbox{and} \qquad \3{\gamma}_j(\3{W}_{j-1}) = \3{W}_j
$$
for each $j = 1, \ldots, n$. For each $j = 1, \ldots, n$, let $\3{\varphi}_j \in \Aut_{\3{\TT}}(\3{K}_j)$ be such that $\3{\rho}(\3{\varphi}_j) = \3{\gamma}_j$. Let also $K_j \in \Ob(\wL)$ be the preimage of $\3{K}_j$, and let $\varphi_j \in \Aut_{\wL}(K_j)$ be such that $\tau(\varphi_j) = \3{\varphi}_j$. Note that $U\leq K_1$ and $V\leq K_n$.

We claim that we may assume that $K_j \in \Ob(\LL^{\bullet})$ for all $j$ without loss of generality. Indeed, we have $\3{K}_j \in \Ob(\3{\TT})$ by (9-c) (where $\3{\TT} = \wL/A$), and thus $K_j \in \Ob(\wL)$. This implies that $(K_j)^{\bullet} \in \Ob(\LL^{\bullet})$ by definition of $\wL$. Furthermore, we have $A\leq (K_j)^{\bullet}$, and it follows that $(K_j)^{\bullet}/A\leq \3{Z}$, since $(K_j)^{\bullet}\leq K_j T$ and $\3{T}\leq \3{Z}$. By Proposition \ref{fix1} (i) and (ii), there exists some $M_U \in \N$ such that
$$
K_j \cap S_i \in \Ob(\LL_i) \qquad \mbox{and} \qquad \Psi_i(\varphi_j) = \varphi_j
$$
for all $i \geq M_U$ and all $j = 1, \ldots, n$. In particular, this implies that $U \cap S_i$ is $\FF_i$-conjugate to $V \cap S_i$. Note that this also implies that $\3{\Psi}_i(\3{\varphi}_j) = \3{\varphi}_j$ for all $i \geq M_U$.

Let $i \geq M_U$, and let $R \in U^T \cap \Omega_i$. By assumption, $\3{U}\leq \3{S}_i$. By (4-a) we have $\3{R} \in \3{U}^{\3{T}}$, which implies that $C_{\3{T}}(\3{R}) = Z_{\3{R}} = Z_{\3{U}}$, and by (4-c) we know that $\3{R}\leq \3{S}_i$. Thus, by Lemma \ref{invar1} (i), together with Proposition \ref{fix1} (iii), we have
$$
\3{x}^{-1} \cdot \3{\Psi}_i(\3{x}) \in C_{\3{T}}(\3{R}) = Z_{\3{U}}\leq \3{R}
$$
for some $\3{x} \in N_{\3{T}}(\3{R}, \3{U})$. Fix such $\3{x} \in N_{\3{T}}(\3{R}, \3{U})$, set $\3{W}'_0 = \3{R}$, and for all $j = 1, \ldots, n$ set
$$
\3{W}'_j = \3{x}^{-1} \cdot \3{W}_j \cdot \3{x} \qquad \mbox{and} \qquad \3{K}'_j = \3{x}^{-1} \cdot \3{K}_j \cdot \3{x}.
$$
Set also $\3{\varphi}_j' = \3{\varepsilon}(\3{x})^{-1} \circ \3{\varphi}_j \circ \3{\varepsilon}(\3{x}) \in \Aut_{\3{\TT}}(\3{K}'_j)$. Let also $W_j'$ and $K_j'$ be the corresponding preimages in $S$. The following holds.
\begin{itemize}

\item[(10-a)] Since $Z_{\3{U}}$ is abelian and $\3{K}_j, \3{K}_j'$ are $\3{\ee}$-centric, it follows that $Z_{\3{U}}\leq \3{K}_j, \3{K}_j'$ for each $j$.

\item[(10-b)] Since $\3{x} \in \3{T}\leq \3{Z} = C_{\3{S}}(Z_{\3{U}})$ and $\3{K}_j\leq \3{Z}$, it follows that $\3{K}'_j\leq \3{Z}$. Moreover, since $\3{x}^{-1} \cdot \3{\Psi}_i(\3{x}) \in Z_{\3{U}}$ and $Z_{\3{U}}\leq C_{\3{T}}(\3{K}_j), C_{\3{T}}(\3{K}'_j)$, it follows from Lemma \ref{invar1} (i) that $\3{x} \in N_{\3{T}}(\3{K}'_j \cap \3{S}_i, \3{K}_j \cap \3{S}_i)$.

\item[(10-c)] Since $Z_{\3{U}}\leq \3{K}_j, \3{K}'_j$, $\3{\varphi}_j \in \Mor(\3{\ee})$, and $\3{x}^{-1} \cdot \3{\Psi}_i(\3{x}) \in Z_{\3{U}}$, it follows from axiom (C) of transporter systems for $\3{\TT}$ that $\3{\varepsilon}(\3{x}^{-1} \Psi_i(\3{x})) \circ \3{\varphi_j}' = \3{\varphi}'_j \circ \3{\varepsilon}(\3{x}^{-1}\Psi_i(\3{x}))$. Thus, $\3{\Psi}_i(\3{\varphi}'_j) = \3{\varphi}'_j$ by Lemma \ref{invar1} (ii). By (4-b), this implies that there is some $\varphi'_j \in \Aut_{\wL}(K'_j)$ such that $\Psi_i(\varphi'_j) = \varphi'_j$ and $\tau(\varphi'_j) = \3{\varphi}'_j$.

\end{itemize}
Let $\3{R'} = \3{W}'_n \in \3{V}^{\3{T}}$, and let $R'$ be its preimage in $S$, and note that $R' \in V^T$ by (4-a). Note that $R\leq K'_1$ and $R'\leq K'_n$. Thus, if $K'_j \cap S_i \in \Ob(\LL_i)$ for each $j$, then $R \cap S_i$ is $\FF_i$-conjugate to $R' \cap S_i$.

It remains to show that $K'_j \cap S_i \in \Ob(\LL_i)$ for each $j$. Recall that, by definition, $K_j' \cap S_i \in \Ob(\LL_i)$ if $(K'_j \cap S_i)^{\bullet} \in \Ob(\LL)$. For each $j = 1, \ldots, n$, let $T_j$ be the maximal torus of $K_j$. Note that $K_j'$ is $T$-conjugate to $K_j$, so in particular $K_j' \in \Ob(\LL^{\bullet})$. As noted above, we have $(T_j \cap S_i)^{\bullet} = T_j$ for all $i \geq M_U$. Also, $T_j \cap S_i\leq K'_j \cap S_i$ and
$$
(T_j \cap S_i)^{\bullet} = T_j\leq (K'_j \cap S_i)^{\bullet}.
$$
In particular, this implies that $K'_j \cap S_i, T_j\leq (K'_j \cap S_i)^{\bullet}$. Since $(K_j \cap S_i)^{\bullet} = K_j$, it follows that $K_j \cap S_i$ contains representatives of all the elements in $\3{K}_j \cap \3{S}_i$. Thus, by (10-b) we deduce that $K_j' \cap S_i$ contains representatives of all the elements in $\3{K}'_j \cap \3{S}_i$, and thus $K'_j = (K'_j \cap S_i)T_j\leq (K'_j \cap S_i)^{\bullet}$. It follows that $(K'_j \cap S_i)^{\bullet} = K_j'$, and thus $K'_j \cap S_i \in \Ob(\LL_i)$.

By the above, $R \cap S_i$ is $\FF_i$-conjugate to $R' \cap S_i$. Moreover, $R' \in V^T$, and $C_S(V) \not\leq V$ by assumption. By Step 8, ($\ddagger$) holds for $V$, and thus the conclusion of ($\ddagger$) holds for $R$ and $R'$. This shows that ($\ddagger$) holds for $U$ as well. Repeating this step, we deduce that ($\ddagger$) holds for all $U \in \hh_X$ such that $Z_U$ is fully $\FF$-normalized, and thus ($\ddagger$) holds for all $U \in \hh_X$ by Step 7. This finishes the proof.
\end{proof}

\begin{thm}\label{fix6}

Every $p$-local compact group admits an approximation by $p$-local finite groups.

\end{thm}

\begin{proof}

Let $\g = \ploc$ be a $p$-local compact group, and let $\wL$ and $\Psi$ be as fixed in \ref{hyp2}. Let also $\Psi_0 = \Psi$, and let $\{\Psi_i\}_{i \geq 0}$ be such that $\Psi_{i+1} = \Psi_i^p$. Let $\{(S_i, \FF_i, \LL_i)\}_{i \geq 0}$ be the family of transporter systems defined in \ref{Li}. For simplicity, we can assume that the degree of $\Psi$ is high enough so that Lemma \ref{SiS} and Propositions \ref{fix2} and \ref{fix5} hold (otherwise replace $\Psi$ by an appropriate power of it), and we claim that $\{(S_i, \FF_i, \LL_i)\}_{i \geq 0}$ is an approximation of $\g$ by $p$-local finite groups.

Condition (i) in \ref{defiapprox} is satisfied, by Lemma \ref{SiS} (i). Also, $S_i$ is a finite $p$-group for all $i \geq 0$, and $\LL_i$ is a linking system associated to $\FF_i$ by Proposition \ref{fix2-1}, and there are inclusions $\LL_i \subseteq \LL_{i+1}$ and $\LL_i \subseteq \wL$ for all $i \geq 0$. Thus, to show that condition (ii) in \ref{defiapprox} is satisfied, it remains to check that $\FF_i$ is saturated and $\Ob(\FF_i^{cr}) \subseteq \Ob(\LL_i)$. By Propositions \ref{fix4} and \ref{fix5}, if $P\leq S_i$ is $\FF_i$-centric but $P \notin \Ob(\LL_i)$, then $P$ is $\FF_i$-conjugate to some $Q$ such that
$$
\Out_{S_i}(Q) \cap O_p(\Out_{\FF_i}(Q)) \neq 1.
$$
Moreover, $\FF_i$ is $\Ob(\LL_i)$-generated and $\Ob(\LL_i)$-saturated by Corollary \ref{fix3}. Thus, the conditions of Theorem \ref{5A} are satisfied: $\FF_i$ is saturated, and $\Ob(\LL_i)$ contains all the centric radical subgroups of $\FF_i$. Thus condition (ii) in \ref{defiapprox} is satisfied.

Finally, we have to check condition (iii): for each $P, Q \in \Ob(\wL)$ and each $\varphi \in \Mor_{\wL}(P,Q)$, there exists some $M \in \N$ such that, for all $i \geq M$, there are objects $P_i, Q_i \in \Ob(\LL_i)$ and morphisms $\varphi_i \in \Mor_{\LL_i}(P_i, Q_i)$, such that $P = \bigcup_{i \geq M} P_i$ and $Q = \bigcup_{i \geq M} Q_i$, and $\4{\varepsilon}_{Q_i, Q}(1) \circ \varphi_i = \varphi \circ \4{\varepsilon}_{P_i, P}(1)$. By Proposition \ref{fix1} (i) and (ii), there is some $M \in \N$ such that $P \cap S_i, Q \cap S_i \in \Ob(\LL_i)$ for all $i \geq M$, and the restriction $\varphi_i = \varphi|_{P \cap S_i}$ is a morphism in $\Mor_{\LL_i}(P \cap S_i, Q \cap S_i)$. Since $S = \bigcup_{i \geq 0} S_i$ by Lemma \ref{SiS} (i), it follows that
$$
P = \bigcup_{i \geq 0} P_i \qquad \mbox{and} \qquad Q = \bigcup_{i \geq 0} Q_i.
$$
The condition $\4{\varepsilon}_{Q_i, Q}(1) \circ \varphi_i = \varphi \circ \4{\varepsilon}_{P_i, P}(1)$ is easily checked.
\end{proof}

\begin{rmk}\label{fix7}

Let $\g = \ploc$ be a $p$-local compact group, and let $\oo(\FF)$ be its \emph{orbit category}: the category $\oo(\FF)$ with object set $\Ob(\FF)$, and with morphism sets
\begin{equation}\label{orbitcat}
\Mor_{\oo(\FF_{\hh})}(P,Q) = \Inn(Q)\backslash \Hom_{\FF}(P,Q).
\end{equation}
Notice that the subcategory $\oo(\FF^{\bullet c}) \subseteq \oo(\FF)$ has a finite (full) subcategory as a skeletal subcategory, which in addition contains a representative of each $\FF$-conjugacy class of centric radical subgroups. Let $\wL$ be the associated telescopic linking system, and let $\{(S_i, \FF_i, \LL_i)\}_{i \geq 0}$ be an approximation of $\g$ by $p$-local finite groups. Denote by $\LL_i^{\bullet} \subseteq \LL^{\bullet}$ the image of $\LL_i$ through the functor $\bullf \colon \wL \to \LL^{\bullet}$ for each $i \geq 0$. Then, by Definition \ref{defiapprox} there is some $M \in \N$ such that, for all $i \geq M$, the category $\LL_i^{\bullet}$ contains representatives of all the morphisms in $\oo(\FF^{\bullet c})$ up to $S$-conjugation. In particular, this implies that
\begin{enumerate}[(a)]

\item the fusion system $\FF$ is generated by $\FF_i$ and $\Inn(S)$; and

\item the linking system $\wL$ is generated by $\LL_i^{\bullet}$ and $S$ (and thus so is $\LL$).

\end{enumerate}
More precisely, property (b) means that every object in $\LL^{\bullet}$ is $S$-conjugate to an object in $\LL_i^{\bullet}$, and every morphism in $\LL^{\bullet}$ is the composition of a morphism in $\LL_i^{\bullet}$ with (the restriction of) a morphism in $\varepsilon(S)\leq \Aut_{\LL}(S)$. Since $\LL^{\bullet}$ is a deformation retract of $\wL$, we can say that $\wL$ is generated by $\LL_i^{\bullet}$ and $S$.

\end{rmk}


\subsection{An example}\label{Ssexample}

In this subsection we analyze our constructions in detail on a specific example: the $2$-local compact group associated to $SO(3)$. In particular, this example reveals that there are approximations by $p$-local finite groups that do not appear as fixed points of any (family of) fine unstable Adams operation.

Let us first fix some notation and facts. The reader is referred to \cite[Example 3.7]{Gonza2} for further details. Let $\g = \ploc$ be the $2$-local compact group associated to $SO(3)$. Then,
$$
S = \gen{\{t_n\}_{n \geq 1}, \, x \, | \, \forall n, \, t_n^{2^n} = x^2 = 1, \, t_{n+1}^2 = t_n, \, x \cdot t_n \cdot x^{-1} = t_n^{-1}} \cong D_{2^{\infty}}.
$$
For each $n \geq 1$, set $T_n = \gen{t_n}$, so that the maximal torus of $S$ is $T = \gen{\{t_n\}_{n \geq 1}}\leq S$. Set also $V = \gen{t_1, x} \cong \Z/2 \times \Z/2$. Then $\FF$ is generated by $\Aut_{\FF}(S) = \Inn(S)$ and $\Aut_{\FF}(V) = \Aut(V) \cong \Sigma_3$. Regarding the linking system $\LL$, the only centric radical subgroups (up to $S$-conjugation) are $S$ and $V$, with
$$
\Aut_{\LL}(S) = S \qquad \mbox{and} \qquad \Aut_{\LL}(V) \cong \Sigma_4.
$$
Let $\wL$ be the associated telescopic linking system. An easy computation shows that the only new objects in $\wL$ are the subgroups $T_n$, for $n \geq 2$, with $\Aut_{\wL}(T_n) = S$ for all $n \geq 2$.

Let now $\Psi \in \Aut_{\typ}^{I}(\LL)$ be an unstable Adams operation. As usual, set $\Psi_0 = \Psi$ and $\Psi_{i+1} = (\Psi_i)^2$. Then, for each $i \geq 0$ we have $S_i = C_S(\Psi_i) \cong D_{2^{n_i}}$ for some $n_i \in \N$, and we may assume without loss of generality that $V\leq S_i$ for all $i \geq 0$. Notice that $S_i$ contains two different $S_i$-conjugacy classes of maximal elementary abelian subgroups, and $V$ is a representative of one of them. Fix a representative $W_i\leq S_i$ of the other $S_i$-conjugacy class of maximal elementary abelian subgroups. Notice that, after embedding $S_i$ into $S_{i+1}$, the subgroup $W_i$ becomes $S_{i+1}$-conjugate to $V$.

Let also $\{(S_i, \FF_i, \LL_i)\}_{i \geq 0}$ be the approximation by $2$-local finite groups associated to $\{\Psi_i\}_{i \geq 0}$. A careful inspection reveals that, in order to describe $(S_i, \FF_i, \LL_i)$ it is enough to specify the groups $\Aut_{\LL_i}(S_i)$, $\Aut_{\LL_i}(V)$ and $\Aut_{\LL_i}(W_i)$. By construction, we have
$$
\Aut_{\LL_i}(S_i) = S_i \qquad \mbox{and} \qquad \Aut_{\LL_i}(V) \cong \Sigma_4,
$$
and we have to determine the group $\Aut_{\LL_i}(W_i)$. Let $\varphi \in \Aut_{\LL_i}(V)$ be an automorphism of order $3$ that conjugates $t_1$ to $x$. An easy computation in $S$ shows that there is some $t \in N_T(W_i, V)$ such that $t^{-1} \cdot \Psi_i(t) = t_1$. By Lemma \ref{invar1} (ii), it follows that $\varepsilon(x)^{-1} \circ \varphi \circ \varepsilon(x)$ is not fixed by $\Psi_i$, and thus
$$
\Aut_{\LL_i}(W_i) = N_{S_i}(W_i) \cong D_8.
$$
These computations imply that, for all $i$, $(S_i, \FF_i, \LL_i)$ is the $2$-local finite group associated to $PGL_2(\F_q)$, where $q$ is some power of some odd prime $p$.

Let now $(S_i, \ee_i, \TT_i)$ be the $2$-local finite group associated to $PSL_2(\F_q)$. This $2$-local finite group is determined by $\Aut_{\TT_i}(S_i) = S_i$ and $\Aut_{\TT_i}(V) \cong \Sigma_4 \cong \Aut_{\TT_i}(W_i)$. Clearly, $\{(S_i, \ee_i, \TT_i)\}_{i \geq 0}$ is an approximation of $\ploc$ by $2$-local finite groups, but our computations above show that it cannot be the product of our constructions in (\ref{Li}).


\section{Stable Elements Theorem for \texorpdfstring{$p$}{p}-local compact groups}\label{Sstable}

In this section we use the approximations constructed in the previous section to prove a version of the Stable Elements Theorem for $p$-local compact groups, which computes the cohomology of the classifying space of a $p$-local compact group with coefficients in a trivial $\Z_{(p)}$-module as the \emph{stable} elements of the cohomology of its Sylow $p$-subgroup with the same module of coefficients. The finite version of this result was proved in \cite[Theorem 5.8]{BLO2} for coefficients in the trivial module $\F_p$, and then generalized to any trivial $\Z_{(p)}$-module in \cite[Lemma 6.12]{BCGLO2}.

Let $S$ be a discrete $p$-toral group, and let $\FF$ be a saturated fusion system. In this section we consider certain contravariant functors $A \colon \FF \to \cat$, where $\cat$ is either $\Z_{(p)}\mathrm{\mathbf{-Mod}}$, the category of $\Z_{(p)}$-modules, or $\mathrm{\mathbf{Gr}}-\Z_{(p)}\mathrm{\mathbf{-Mod}}$, the category of graded $\Z_{(p)}$-modules.  We start with a brief discussion of some general properties of such functors. For each $X\leq S$, let $\iota_X \colon X \to S$ denote the inclusion homomorphism. This way, given a contravariant functor $A \colon \FF \to \cat$ and a subcategory $\ee \subseteq \FF$, we define
\begin{equation}\label{Aee}
A^{\ee} \defin \{z \in A(S) \, | \, A(\iota_P)(z) = A(\iota_Q \circ f)(z), \, \forall P, Q \in \Ob(\ee) \mbox{ and } \forall f \in \Hom_{\ee}(P,Q)\}.
\end{equation}
Given $\ee_1 \subseteq \ee_2$ two subcategories of $\FF$, there is an obvious inclusion $A^{\ee_2} \subseteq A^{\ee_1}$.

\begin{lmm}\label{aux34}

Let $S$ be a discrete $p$-toral group, let $\FF$ be a fusion system over $S$, and let $\{(S_i, \FF_i)\}_{i \geq 0}$ be a family of finite fusion subsystems of $\FF$ with $\FF_i \subseteq \FF_{i+1}$ for all $i$ (in particular, $S_i\leq S_{i+1}$ are finite subgroups of $S$ for all $i$), and satisfying the following properties:
\begin{enumerate}[(i)]

\item $S = \bigcup_{i \geq 0} S_i$; and

\item for all $P\leq S$ and for all $f \in \Hom_{\FF}(P, S)$ there exists some $M_f \in \N$ such that, for all $i \geq M_f$, $f|_{P \cap S_i} \in \Hom_{\FF_i}(P \cap S_i, S_i)$.

\end{enumerate}
Let also $\cat$ be either $\Z_{(p)}\mathrm{\mathbf{-Mod}}$, the category of $\Z_{(p)}$-modules, or $\mathrm{\mathbf{Gr}}-\Z_{(p)}\mathrm{\mathbf{-Mod}}$, the category of graded-$\Z_{(p)}$-modules, and let $A \colon \FF \to \cat$ be a contravariant functor satisfying the following property:
\begin{itemize}

\item[(\textasteriskcentered)] For each $P\leq S$, the natural map $A(P) \to \varprojlim_i A(P \cap S_i)$ is an isomorphism.

\end{itemize}
Then, upon setting $\FF^{\circ} = \bigcup_{i \geq 0} \FF_i \subseteq \FF$, there are equalities
$$
A^{\FF} = A^{\FF^{\circ}} = \varprojlim_i A^{\FF_i}
$$
as subsets of $A(S)$.

\end{lmm}

\begin{proof}

We claim first that $\FF^{\circ} \subseteq \FF$ is the full subcategory of $\FF$ whose objects are the finite subgroups of $S$. Indeed, if $P\leq S$ is a finite subgroup, then there exists some $M_P \in \N$ such that $P\leq S_i$ for all $i \geq M_P$, since $S = \bigcup_{i \geq 0} S_i$. Similarly, if $P, Q\leq S$ are finite subgroups and $f \in \Hom_{\FF}(P,Q)$, then by condition (ii) there exists some $M_f \in \N$ such that $P, Q\leq S_i$, and $f = f|_{P \cap S_i} \in \Hom_{\FF_i}(P,Q)$ for all $i \geq M_f$.

By (\textasteriskcentered), we have $A(S) = \varprojlim_i A(S_i)$, which implies that $A^{\FF^{\circ}} = \varprojlim_i A^{\FF_i}$. The inclusion $\FF^{\circ} \subseteq \FF$ implies that $A^{\FF} \subseteq A^{\FF^{\circ}}$ by (\ref{Aee}). To show the reverse inclusion, let $P, Q \in \Ob(\FF)$ and $f \in \Hom_{\FF}(P,Q)$. For $X = P, Q, S$, set $X_i = X \cap S_i$, and notice that
$$
X = \bigcup_{i \geq 0} X_i \qquad \mbox{and} \qquad A(X) = \varprojlim_i A(X_i),
$$
by (i) and (\textasteriskcentered) respectively. By condition (ii), there exists some $M_f \in \N$ such that,  for all $i \geq M_f$, $f|_{P_i} \in \Hom_{\FF_i}(P_i, Q_i)$. Thus, if $z \in A^{\FF^{\circ}}$, then
$$
A(\iota_{P_i})(z) = A(\iota_{Q_i} \circ f|_{P_i})(z)
$$
for all $i \geq M_f$, and thus $A(\iota_P)(z) = A(\iota_Q \circ f)(z)$. Hence $A^{\FF^{\circ}} \subseteq A^{\FF}$.
\end{proof}

\begin{rmk}

If $\g = \ploc$ is a $p$-local compact group and $\{(S_i, \FF_i, \LL_i)\}_{i \geq 0}$ is an approximation of $\g$ by $p$-local finite groups, then conditions (i) and (ii) in Lemma \ref{aux34} follow from condition (i) in \ref{defiapprox} and \ref{finmorph}, respectively. Hence, in this case \ref{aux34} applies to any functor $A \colon \FF \to \cat$ as above that satisfies condition (\textasteriskcentered).

\end{rmk}

\begin{prop}\label{stable1}

Let $\g = \ploc$ be a $p$-local compact group, and let $M$ be a (finite) $\Z_{(p)}$-module with trivial $S$-action. Let also $\{(S_i, \FF_i, \LL_i)\}_{i \geq 0}$ be an approximation of $\g$ by $p$-local finite groups, and let $P\leq S$. Then, there are isomorphisms
$$
H^{\ast}(BP; M) \cong \varprojlim H^{\ast}(B(P \cap S_i); M) \qquad \mbox{and} \qquad H^{\ast}(B\g; M) \cong \varprojlim H^{\ast}(B\g_i; M).
$$
In particular, the functor $H^{\ast}(-; M) \colon \FF \to \mathrm{\mathbf{Gr}}-\Z_{(p)}\mathrm{\mathbf{-Mod}}$ satisfies condition (\textasteriskcentered) in Lemma \ref{aux34}.

\end{prop}

\begin{proof}

Fix some $P\leq S$, and set $P_i = P \cap S_i$ for all $i \geq 0$. Let $X$ be either $B\g$ or $BP$, and similarly let $X_i$ be either $B\g_i$ or $BP_i$, depending on which case we want to prove. Note that the following holds.
\begin{enumerate}[(i)]

\item If $X = BP$, then $X = \hocolim X_i$, since $P = \bigcup_{i \geq 0} P_i$ by hypothesis.

\item If $X = B\g$, then $X \simeq (\hocolim X_i)^{\wedge}_p$ by Lemma \ref{approx0} and Remark \ref{approx-1}. In particular, $H^{\ast}(X; M) \cong H^{\ast}(\hocolim X_i; M)$.

\end{enumerate}
Consider the homotopy colimit spectral sequence for cohomology \cite[XII.5.7]{BK}:
$$
E^{r,s}_2 = \varprojlim \!\! \phantom{i}^rH^s(X_i;M) \Longrightarrow H^{r+s}(X;M).
$$
We will see that, for $r \geq 1$, $E_2^{r,s} = \{0\}$, which, in particular, will imply the statement.

For each $s$, let $H^s_i = H^s(X_i;M)$, and let $F_i$ be the induced morphism in cohomology (in degree $s$) by the map $|\Theta_i| \colon |\LL_i| \to |\LL_{i+1}|$ induced by the inclusion $\LL_i \subseteq \LL_{i+1}$. The cohomology ring $H^{\ast}(X_i;M)$ is noetherian by \cite[Theorem 5.8]{BLO2}, and in particular $H^s_i$ is finite for all $s$ and all $i$. Thus, the inverse system $\{H^s_i;F_i\}$ satisfies the Mittag-Leffler condition \cite[3.5.6]{Weibel}, and hence the higher limits $\varprojlim^rH^s_i$ vanish for all $r \geq 1$. This in turn implies that the differentials in the above spectral sequence are all trivial, and thus it collapses.
\end{proof}

We are ready to prove Theorem \ref{thmB}, the Stable Elements Theorem for $p$-local compact groups, which we restate below.

\begin{thm}\label{stable2}

Let $\g = \ploc$ be a $p$-local compact group, and let $M$ be a finite $\Z_{(p)}$-module $M$ with trivial $S$-action. Then, the natural map
$$
H^{\ast}(B\g; M) \Right3{\cong} H^{\ast}(\FF;M) \defin \varprojlim_{\FF} H^{\ast}(-; M) \subseteq H^{\ast}(BS; M)
$$
is an isomorphism.

\end{thm}

\begin{proof}

Let $\{\g_i = (S_i, \FF_i, \LL_i)\}_{i \geq 0}$ be an approximation of $\g$ by $p$-local finite groups, with respect to some telescopic transporter system $\wL$ satisfying the conditions in \ref{defiapprox}. By Remark \ref{approx-1}, the space $B\g_i \defin |\LL_i|^{\wedge}_p$ is the classifying space of a $p$-local finite group, and we can apply the Stable Elements Theorem for $p$-local finite groups: there is a natural isomorphism
$$
H^{\ast}(B\g_i; M) \stackrel{\cong} \longrightarrow H^{\ast}(\FF_i;M).
$$
By Proposition \ref{stable1} there are natural isomorphisms
$$
H^{\ast}(B\g;M) \cong \varprojlim_i H^{\ast}(B\g_i; M) \cong \varprojlim_i H^{\ast}(\FF_i;M) \subseteq \varprojlim_i H^{\ast}(BS_i; M) \cong H^{\ast}(BS;M),
$$
and we have to show that $\varprojlim_i H^{\ast}(\FF_i;M) \cong H^{\ast}(\FF;M)$. If we set $\FF^{\circ} = \bigcup_{i \geq 0} \FF_i$, then Lemma \ref{aux34} and Proposition \ref{stable1} combined imply that there are isomorphisms
$$
\varprojlim_i H^{\ast}(\FF_i;M) \cong H^{\ast}(\FF^{\circ};M) \cong H^{\ast}(\FF;M),
$$
and this finishes the proof.
\end{proof}

\begin{rmk}\label{stable21}

The reader may have noticed the difference between the original statement of the Stable Elements Theorem for $p$-local finite groups, \cite[Theorem 5.8]{BLO3}, in terms of the orbit category of $\FF$ defined in (\ref{orbitcat}), and our statement, in terms of $\FF$. In fact, a formulation of the Stable Elements Theorem in terms of $\FF$, rather than $\oo(\FF)$, is already found in \cite[6.12]{BCGLO2}, as well as in other papers, and it is rather straightforward to justify the equivalence of statements. Consider the projection functor $\tau \colon \FF \to \oo(\FF)$, which is the identity on objects. Since conjugation by elements of $S$ induces the identity on cohomology, we have a commutative triangle
$$
\xymatrix{
\FF \ar[rrr]^{H^{\ast}(-)} \ar[d]_{\tau} & & & \mathrm{\mathbf{Gr}}-\Z_{(p)}\mathrm{\mathbf{-Mod}}\\
\oo(\FF) \ar[rrru]_{H^{\ast}(-)} & & &
}
$$
and an induced morphism between the corresponding inverse limits
$$
\varprojlim_{\oo(\FF)} H^{\ast}(-;\F_p) \Right2{} \varprojlim_{\FF} H^{\ast}(-;\F_p).
$$
This morphism is easily checked to be an isomorphism upon considering both groups as subgroups of stable elements in $H^{\ast}(S; \F_p)$, since every element of $H^{\ast}(S; \F_p)$ is stable by any morphism in $\FF_S(S)$.

\end{rmk}

We finish this section proving Theorem \ref{thmD}, restated as Theorem \ref{thmd} below, which states the existence of a certain spectral sequence associated to a strongly closed subgroup of a given saturated fusion system. We first fix some notation.

Let $\FF$ be a saturated fusion system over a discrete $p$-toral group $S$, let $R\leq S$ be a strongly $\FF$-closed subgroup, and let $M$ be an $\Z_{(p)}$-module with trivial $S$-action. For each $X\leq S$, set $\overline{X} = X/(X \cap R) \cong XR/R\leq S/R$. Then, each $f \in \Hom_{\FF}(P,Q)$ induces a morphism of extensions
$$
\xymatrix{
P \cap R \ar[r] \ar[d]_{f_0} & P \ar[d]^f \ar[r] & \overline{P} \ar[d]^{\overline{f}} \\
Q \cap R \ar[r] & Q \ar[r] & \overline{Q},
}
$$
and thus also homomorphisms
$$
\gamma(f) \colon H^n(\overline{Q}; H^m(Q \cap R;M)) \Right2{f_0^{\ast}} H^n(\overline{Q}; H^m(P \cap R; M))  \Right2{\overline{f}^{\,\ast}} H^n(\overline{P}; H^m(P \cap R;M))
$$
for all $n, m \geq 0$. This defines a contravariant functor
\begin{equation}\label{frakx}
\frakx^{n,m} \colon \FF \Right3{} \Z_{(p)}\mathrm{\mathbf{-Mod}},
\end{equation}
which satisfies condition (\textasteriskcentered) in Lemma \ref{aux34}, since, for each $P\leq S$ and each $i \geq 0$, the group $H^n(\overline{P \cap S_i}; H^m(P \cap S_i \cap R; M))$ is finite and the Mittag-Leffler Condition \cite[3.5.6]{Weibel} applies. Set
$$
H^n(S/R; H^m(R; M))^{\FF} = (\frakx^{n,m})^{\FF}
$$
in order to match the notation of \cite{Diaz}

\begin{thm}\label{thmd}

Let $\g = \ploc$ be a $p$-local compact group, let $R\leq S$ be a strongly $\FF$-closed subgroup, and let $M$ be a finite $\Z_{(p)}$-module with trivial $S$-action. Then there is a first quadrant cohomological spectral sequence with second page
$$
E^{n,m}_2 = H^n(S/R; H^m(R;M))^{\FF}
$$
and converging to $H^{n+m}(B\g; M)$.

\end{thm}

\begin{proof}

The spectral sequence of the statement is constructed as an inverse limit of spectral sequences. For the reader's convenience, the proof is divided into smaller steps.

\textbf{Step 1.} Construction of inverse systems of spectral sequences. Let $\{(S_i, \FF_i, \LL_i)\}_{i \geq 0}$ be an approximation of $\g$ by $p$-local finite groups. Notice that $R_i \defin R \cap S_i$ is strongly $\FF_i$-closed for all $i \geq 0$, since $R$ is strongly $\FF$-closed.

For each $i \geq 0$, consider the group extension
$$
R_i \Right2{} S_i \Right2{} S_i/R_i,
$$
and its associated Lyndon-Hochschild-Serre spectral sequence $\{E(i)_k^{\ast, \ast}, d_k^i\}_{k \geq 2}$, with second page $E(i)_2^{n,m} = H^n(S_i/R_i; H^m(R_i;M))$, and converging to $H^{n+m}(S;M)$. For all $k \geq 2$ and all $i, n, m \geq 0$, the inclusion homomorphism $\iota_{i,i+1} \colon S_i \to S_{i+1}$ induces homomorphisms of spectral sequences for each $i \geq 0$,
$$
\gamma(\iota_{i, i+1}) \colon \{E(i+1)_k^{\ast,\ast}, d_k^{i+1}\}_{k \geq 2} \Right3{} \{E(i)_k^{\ast,\ast}, d_k^i\}_{k \geq 2}
$$
and thus an inverse system of spectral sequences $\{\{E(i)_k^{\ast,\ast}, d_k^i\}_{k \geq 2}, \gamma(\iota_{i, i+1})\}_{i \geq 0}$.

Note that $R_i$ is strongly $\FF_i$-closed, and thus for each $n, m \geq 0$ there is a functor
$$
\frakx_i^{n,m} \colon \FF_i \Right3{} \Z_{(p)}\mathrm{\mathbf{-Mod}},
$$
defined by similar arguments as those used in the definition the functor $\frakx^{n,m}$ in (\ref{frakx}). Furthermore, we may consider the spectral sequence $\{\4{E}(i)_k^{\ast, \ast}, \4{d}_k^i\}_{k \geq 2}$ of \cite[Theorem 1.1]{Diaz}, whose second page is $\4{E}(i)_2^{n,m} = H^n(S_i/R_i; H^m(R_i;M))^{\FF_i} = (\frakx_i^{n,m})^{\FF_i}$, and which converges to $H^{n+m}(B\g_i;M)$. We may see this spectral sequence as the restriction of the spectral sequence $\{E(i)_k^{\ast, \ast}, d_k^i\}_{k \geq 2}$. We claim that $\gamma(\iota_{i,i+1})$ restricts to a morphism of spectral sequences
$$
\gamma(\alpha_{i, i+1}) \colon \{\4{E}(i+1)_k^{\ast,\ast}, \4{d}_k^{i+1}\}_{k \geq 2} \Right3{} \{\4{E}(i)_k^{\ast,\ast}, \4{d}_k^i\}_{k \geq 2}.
$$

For each $i,n, m \geq 0$ and for $X\leq S_i$, let $A^{n,m}_i(X) = \Hom_X(\bb^n_{\overline{X}} \otimes \bb^m_X, M)$ be the double complex defined in \cite[Section 3]{Diaz} (although we do not give an explicit description here, notice that the complex itself does not actually depend on $i$). With this notation, for each $k \geq 2$, the group $\4{E}(i)_k^{n,m}$ can be seen as $(\xi(i)^{n,m}_k)^{\FF_i}$, for a certain functor
$$
\xi(i)_k^{n,m} \colon \FF_i \Right3{} \Z_{(p)}\mathrm{\mathbf{-Mod}},
$$
defined in terms of the double complexes $A^{n,m}_i(X)$ above. In particular, the functor $\frakx_i^{n,m}$ defined above corresponding to the case $k = 2$ (for the sake of brevity, the reader is referred to \cite[Section 4]{Diaz} for details). Furthermore, we have a commutative triangle of functors
$$
\xymatrix{
\FF_i \ar[rrr]^{\xi(i)_k^{n,m}} \ar[d]_{\incl} & & & \Z_{(p)}\mathrm{\mathbf{-Mod}} \\
\FF_{i+1} \ar[rrru]_{\xi(i+1)_k^{n,m}} & & &
}
$$

Fix some $k \geq 2$, and some $i, n, m \geq 0$. For a subgroup $X\leq S_i\leq S_{i+1}$, let $
\iota_X \colon X \to S_i$ and $\4{\iota}_X \colon X \to S_{i+1}$ be the corresponding inclusion monomorphisms, and note that $\4{\iota}_X = \iota_{i, i+1} \circ \iota_X$. Also, recall that $\4{E}(i+1)_k^{n,m}$ is the subgroup of elements $z \in E(i+1)_k^{n,m}$ such that
$$
\xi(i+1)_k^{n,m}(\4{\iota}_P)(z) = \xi(i+1)_k^{n,m}(\4{\iota}_Q \circ f)(z)
$$
for all $P,Q\leq S_{i+1}$ and all $f \in \Hom_{\FF_{i+1}}(P,Q)$, and a similar, equation with $\xi(i)_k^{n,m}$ replacing $\xi(i+1)_k^{n,m}$, describes the elements of $\4{E}(i)_k^{n,m}$.

Fix some $z \in E(i+1)_k^{n,m}$, and set $w = \gamma(\iota_{i, i+1})(z) \in E(i)_k^{n,m}$. Then, for all $P, Q\leq S_i$ and all $f \in \Hom_{\FF_i}(P,Q) \subseteq \Hom_{\FF_{i+1}}(P,Q)$, the commutativity of the above triangle implies that
$$
\begin{aligned}
\xi(i)_k^{n,m}(\iota_P)(w) & = \xi(i+1)_k^{n,m}(\iota_{i, i+1} \circ \iota_P)(z) = \xi(i+1)_k^{n,m}(\4{\iota}_P)(z) = \\
 & = \xi(i+1)_k^{n,m}(\4{\iota}_Q \circ f)(z) = \xi(i+1)_k^{n,m}(\iota_{i,i+1} \circ \iota_Q \circ f)(z) = \\
 & = \xi(i)_k^{n,m}(\iota_Q \circ f)(w).
\end{aligned}
$$
Thus $w \in \4{E}(i)_k^{n,m}$ and the claim follows.

\textbf{Step 2.} The inverse limit spectral sequences and their convergence. Consider the inverse systems of spectral sequences defined in Step 1, $\{\{E(i)_k^{\ast,\ast}, d_k^i\}_{k \geq 2}, \gamma(\iota_{i, i+1})\}_{i \geq 0}$ and $\{\{\4{E}(i)_k^{\ast,\ast}, \4{d}_k^i\}_{k \geq 2}, \gamma(\alpha_{i, i+1})\}_{i \geq 0}$. For each $k \geq 2$, and for each $n, m \geq 0$, define
$$
E_k^{n,m} \defin \varprojlim_i E_k^{n,m}(i) \qquad \qquad \4{E}_k^{n,m} \defin \varprojlim_i \4{E}_k^{n,m}(i).
$$

For each $k \geq 2$ and each $i, n,m  \geq 0$, the group $E(i)_k^{n,m}$ is finite by definition, and thus the higher limits of $\{E(i)_k^{n,m}\}_{i \geq 0}$ all vanish, by the Mittag-Leffler Condition \cite[3.5.6]{Weibel}. A similar conclusion applies to $\{\4{E}(i)_k^{n,m}\}_{i \geq 0}$, since it is a restriction of $\{E(i)_k^{\ast, \ast}\}_{i \geq 0}$. Furthermore, the differentials $\{d_k^i\}_{i \geq 0}$ and $\{\4{d}_k^i\}_{i \geq 0}$ induce differentials $d_k$ (on $E_k^{\ast, \ast}$) and $\4{d}_k$ (on $\4{E}_k^{\ast, \ast}$) respectively, and an easy computation shows that
$$
\begin{array}{c}
\Ker(d_k) \cong \varprojlim_i \Ker(d_k^i) \qquad \qquad \Im(d_k) \cong \varprojlim_i \Im(d_k^i)\\[4pt]
\Ker(\4{d}_k) \cong \varprojlim_i \Ker(\4{d}_k^i) \qquad \qquad \Im(\4{d}_k) \cong \varprojlim_i \Im(\4{d}_k^i).
\end{array}
$$
Hence, $\{E_k^{\ast, \ast}, d_k\}_{k \geq 2}$ and $\{\4{E}_k^{\ast, \ast}, \4{d}_k\}_{k \geq 2}$ are well defined spectral sequences. Moreover, their corresponding second and infinity pages are, respectively,
$$
\begin{array}{c}
E_2^{n,m} \cong \varprojlim_i H^n(S_i/R_i; H^m(R_i;M)) \qquad \4{E}_2^{n,m} \cong \varprojlim_i H^n(S_i/R_i; H^m(R_i;M))^{\FF_i}\\[6pt]
E_{\infty}^{n,m} \cong \varprojlim_i H^{k}(S_i;M) \cong H^{k}(S; M) \qquad \4{E}_{\infty}^{n,m} \cong \varprojlim_i H^{k}(B\g_i; M) \cong H^{k}(B\g;M)
\end{array}
$$
where $k = n+m$, and where the last isomorphism on each position of the bottom row holds by Proposition \ref{stable1}.

It remains to show that the corresponding second pages of the limit spectral sequences satisfy, respectively,
$$
\begin{array}{l}
E_2^{n,m} = \varprojlim_i H^n(S_i/R_i; H^m(R_i;M)) \cong H^n(S/R; H^m(R;M)) \\[2pt]
\4{E}_2^{n,m} = \varprojlim_i H^n(S_i/R_i; H^m(R_i;M))^{\FF_i} \cong H^n(S/R; H^m(R;M))^{\FF}
\end{array}
$$
By Proposition \ref{stable1}, there are isomorphisms 
$$
H^{\ast}(R; M) \cong \varprojlim H^{\ast}(R_i;M) \quad \mbox{and} \quad H^{\ast}(S/R; H^{\ast}(R;M))\cong \varprojlim H^{\ast}(S_i/R_i; H^{\ast}(R;M)).
$$
Furthermore, \cite[Proposition B.2.3]{Rubin} implies that there are isomorphisms
$$
\begin{aligned}
H^n(S/R; H^m(R;M)) & \cong \varprojlim_i H^n(S_i/R_i; H^m(R;M)) \cong \\
 & \cong \varprojlim_i \varprojlim_j H^n(S_i/R_i; H^m(R_j;M)) \cong \varprojlim_i H^n(S_i/R_i; H^m(R_i;M)),
\end{aligned}
$$
since $H^n(S_i/R_i; H^m(R_j;M))$ is finite for all $i, n, m \geq 0$. To finish the proof, we have to show that the $H^n(S/R; H^m(R;M))^{\FF}$ is isomorphic to $\varprojlim_i H^n(S_iR_i; H^m(R_i;M))^{\FF_i}$, and this follows from Lemma \ref{aux34}, since the functor $\frakx^{n,m}$ in (\ref{frakx}) satisfies the required condition (\textasteriskcentered).
\end{proof}


\section{Mapping spaces}\label{Smap}

In this section we describe the mapping space $\map(BP, B\g)$, where $P$ is a discrete $p$-toral group and $\g$ is a $p$-local compact group, in terms of centralizers in $\g$ of subgroups of $S$. Such mapping spaces where described in \cite{BLO3} when $P\leq S$ is centric, and in full generality in \cite{BLO2} when $P$ is a finite $p$-group and $\g$ is a $p$-local finite group. Our proof follows the same lines as the proof in \cite[Theorem 6.3]{BLO2}.

Let $\g = \ploc$ be a $p$-local compact group, and let $H^{\ast}(\FF) \subseteq H^{\ast}(BS)$ be the subring of stable elements for $\FF$. Let also $E\leq S$ be an elementary abelian subgroup which is fully $\FF$-centralized, and let $j_E \colon H^{\ast}(\FF) \to H^{\ast}(BE)$ be the map induced by inclusion. Let $T_E$ be Lannes' $T$-functor (see \cite{Lannes}), and let $T_E(H^{\ast}(\FF); j_E)$ be the component in $T_E(H^{\ast}(\FF))$ of $j_E \in T_E^0(H^{\ast}(\FF)) \cong \Hom_{\kk}(H^{\ast}(\FF), H^{\ast}(BE))$, where $\kk$ is the category of unstable algebras over the mod $p$ Steenrod algebra.

\begin{lmm}\label{mapping1}

There is an isomorphism
$$
T_E(H^{\ast}(\FF); j_E) \Right3{\cong} H^{\ast}(C_{\FF}(E)) \defin \varprojlim_{C_{\FF}(E)} H^{\ast}(-),
$$
which is the restriction of the homomorphism $T_E(H^{\ast}(BS);j_E) \to H^{\ast}(C_S(E))$ induced by the natural homomorphism $C_S(E) \times E \to S$.

\end{lmm}

\begin{proof}

Let $\{(S_i, \FF_i, \LL_i)\}_{i \geq 0}$ be an approximation of $\g$ by $p$-local finite groups. For the reader's convenience, we have divided the proof into several shorter steps.

\textbf{Step 1.} Given a discrete $p$-toral group $X$ and a morphism $\sigma \colon E \to X$, we claim that
$$
T_E(H^{\ast}(BX);\sigma^{\ast}) \cong H^{\ast}(\Map(BE, BX)_{B\sigma}) \cong H^{\ast}(BC_X(\sigma E)).
$$
Indeed, since $X$ is discrete $p$-toral, it follows that $BX$ is $p$-good and $H^{\ast}(BX)$ is of finite type (i.e. finite in every dimension). Moreover, since both $E$ and $X$ are discrete groups, it follows that $\Map(BE, BX)_{B\sigma} \simeq BC_X(\sigma E)$. The claim follows immediately by \cite[Proposition 3.4.4]{Lannes}.

\textbf{Step 2.} Next we claim that there is some $M \in \N$ such that $E$ is fully centralized in $\FF_i$, for all $i \geq M$.

To prove this, fix representatives $V_0 = E, V_1, \ldots, V_n\leq S$ of the different $S$-conjugacy classes in $E^{\FF}$. Note that, since $E$ is abelian, we have $VC_S(V) = C_S(V)$ for each $V \in E^{\FF}$. We will just write $C_S(V)$ instead of $VC_S(V)$ for simplicity. There is some $M \in \N$ such that, for all $i \geq M$, we have $V_0, \ldots, V_n\leq S_i$, and the following holds for each $j = 0, \ldots, n$:
\begin{enumerate}[(a)]

\item $|C_{S_i}(V_j)|\leq |C_{S_i}(E)|$; and

\item $|C_{S_i}(V_j)/C_{T_i}(V_j)| = |C_S(V_j)/C_T(V_j)|$.

\end{enumerate}
Indeed, since $V_j$ is $\FF$-conjugate to $E$, there exists some $f_j \in \Hom_{\FF}(C_S(V_j), C_S(E))$ such that $f_j(V_j) = E$. Moreover, by Lemma \ref{finmorph}, there is some $M \in \N$ such that, for all $i \geq M$ and all $j = 0, \ldots, n$, the restriction of $f_j$ to $C_S(V_j) \cap S_i$ is a morphism in $\FF_i$. Property (a) follows immediately. Property (b) is easily checked, since $S/T \cong S_i/T_i$ for all $i \geq 0$ by assumption.

For each $i \geq M$, set $T_i = T \cap S_i$ for short. If $V\leq S_i$ is $\FF_i$-conjugate to $E$, then $V$ is $S$-conjugate to some $V_j$, for some $j \in \{1, \ldots, n\}$. Fix $x \in N_S(V, V_j)$. By Lemma \ref{invar1} (i), together with Proposition \ref{fix1} (iv), it follows that $x^{-1} \Psi_i(x) \in C_T(V)$. Since $T$ is abelian and $T_i\leq T$, it follows that
$$
x^{-1} \Psi_i(x) \in C_T(C_{T_i}(V)) = T.
$$
Thus, $C_{T_i}(V)$ is $S$-conjugate (by $x$) to a subgroup of $C_{T_i}(V_j)$, again by Lemma \ref{invar1} (i), and a similar argument with $x^{-1}$ instead of $x$ shows that in fact the element $x$ conjugates $C_{T_i}(V)$ onto $C_{T_i}(V_j)$. Moreover, via the inclusion $C_{S_i}(V)\leq C_S(V)$, the quotient $C_{S_i}(V)/C_{T_i}(V)$ can be identified with a subgroup of $C_S(V)/C_T(V)$, and
$$
|C_{S_i}(V)/C_{T_i}(V)|\leq |C_S(V)/C_T(V)| = |C_S(V_j)/C_T(V_j)| = |C_{S_i}(V_j)/C_{T_i}(V_j)|.
$$
Since $C_{S_i}(V)$ and $C_{S_i}(V_j)$ are finite groups, it follows that
$$
\begin{aligned}
|C_{S_i}(V)| & = |C_{T_i}(V)| \cdot |C_{S_i}(V)/C_{T_i}(V)|\leq \\
 &\leq |C_{T_i}(V_j)| \cdot  |C_{S_i}(V_j)/C_{T_i}(V_j)| = |C_{S_i}(V_j)|\leq |C_{S_i}(E)|,
\end{aligned}
$$
and $E$ is fully centralized in $\FF_i$.

\textbf{Step 3.} For each $P\leq S$, set $\TT_P = \Rep_{\FF}(E,P) = \Inn(P)\backslash \Hom_{\FF}(E,P)$. Notice that $\TT_P$ is finite by \cite[Lemma 2.5]{BLO3}. Consider the functor $\4{T}_E \colon \oo(\FF) \Right2{} \mathrm{\mathbf{Gr}}-\Z_{(p)}\mathrm{\mathbf{-Mod}}$, defined on objects by
$$
\4{T}_E(P) = \bigoplus_{\rho \in \TT_P} T_E(H^{\ast}(BP);\rho^{\ast}).
$$
We claim that there is an isomorphism
$$
T_E(H^{\ast}(\FF);j_E) \cong \varprojlim_{P \in \oo(\FF)} \big(\bigoplus_{\rho \in \TT_P} T_E(H^{\ast}(BP);\rho^{\ast}) \big).
$$

To prove this, consider the orbit category $\oo(\FF)$, defined in (\ref{orbitcat}). By \cite[Lemma 2.5]{BLO3}, all morphism sets in $\oo(\FF)$ are finite. Recall also that the full subcategory $\FF^{\bullet} \subseteq \FF$ contains only finitely many $\FF$-conjugacy classes by \cite[Lemma 3.2 (a)]{BLO3}, and thus the full subcategory $\oo(\FF^{\bullet}) \subseteq \oo(\FF)$ contains a finite skeletal subcategory. Furthermore, $\FF^{\bullet}$ contains all the $\FF$-centric $\FF$-radical subgroups of $S$ by \cite[Corollary 3.5]{BLO3}.

Fix a finite skeletal subcategory $\oo_{sk}$ of $\oo(\FF^{\bullet})$. The functor $T_E$ is exact and commutes with direct limits. As a consequence, it also commutes with inverse limits over finite categories, and we have
$$
\begin{aligned}
T_E(H^{\ast}(\FF)) & = T_E(\varprojlim_{\oo(\FF)} H^{\ast}(-)) = \\
 & =  T_E(\varprojlim_{\oo_{sk}} H^{\ast}(-)) \cong \varprojlim_{\oo_{sk}} T_E(H^{\ast}(-)) = \varprojlim_{\oo(\FF)} T_E(H^{\ast}(-)),
\end{aligned}
$$
Restricting the above to $T_E(H^{\ast}(\FF);j_E)$, we obtain
$$
\begin{aligned}
T_E(H^{\ast}(\FF);j_E) & \cong \varprojlim_{P \in \oo_{sk}} \big(\bigoplus_{\rho \in \TT_P} T_E(H^{\ast}(BP);\rho^{\ast}) \big) \cong \\
 & \cong \varprojlim_{P \in \oo(\FF)} \big(\bigoplus_{\rho \in \TT_P} T_E(H^{\ast}(BP);\rho^{\ast}) \big) = \varprojlim_{\oo(\FF)} \4{T}_E(-).
\end{aligned}
$$

\textbf{Step 4.} Consider the functor $\4{T}_E \colon \oo(\FF) \Right2{} \mathrm{\mathbf{Gr}}-\Z_{(p)}\mathrm{\mathbf{-Mod}}$ defined in Step 3. By precomposing with the projection functor $\tau \colon \FF \to \oo(\FF)$, we obtain a functor
$$
\FF \Right2{} \mathrm{\mathbf{Gr}}-\Z_{(p)}\mathrm{\mathbf{-Mod}},
$$
which by abuse of notation we also denote by $\4{T}_E$. This will not lead to confusion since both functors take the same values on objects (as $\tau$ is the identity on objects), as well as the same value on any two morphisms in $\FF$ representing a given morphism in $\oo(\FF)$. Note also that, for each $P \in \Ob(\FF) = \Ob(\oo(\FF))$, the set $\TT_P$ does not depend on the choice of category between $\FF$ and $\oo(\FF)$, and an argument similar to that used in Remark \ref{stable21} implies that
$$
\begin{aligned}
\varprojlim_{\oo(\FF)} \4{T}_E(-) & = \varprojlim_{P \in \oo(\FF)} \big(\bigoplus_{\rho \in \TT_P} T_E(H^{\ast}(BP);\rho^{\ast}) \big) \cong\\
 & \cong \varprojlim_{P \in \FF} \big(\bigoplus_{\rho \in \TT_P} T_E(H^{\ast}(BP);\rho^{\ast}) \big) = \varprojlim_{\FF} \4{T}_E(-).
\end{aligned}
$$
From now on we consider $\4{T}_E$ as a functor on $\FF$. In this step we prove that the functor $\4{T}_E$ satisfies condition (\textasteriskcentered) in Lemma \ref{aux34}: for each $P\leq S$, the natural map $\4{T}_E(P) \to \varprojlim_i \4{T}_E(P \cap S_i)$ is an isomorphism.

Fix $P\leq S$, and let $\Omega \subseteq \Hom_{\FF}(E,P)$ be a set of representatives of the classes in $\TT_P$. Note that, by definition, there is an equality
\begin{equation}\label{TEP}
\4{T}_E(P) \defin \bigoplus_{\rho \in \TT_P} T_E(H^{\ast}(BP);\rho^{\ast}) = \bigoplus_{f \in \Omega} T_E(H^{\ast}(BP); f^{\ast}).
\end{equation}
For each $i \geq 0$, set $P_i = P \cap S_i$. Since $\Omega$ is a finite set, there exists some $M_P \in \N$ such that $f(E)\leq P_i$ for all $i \geq M_P$ and all $f \in \Omega$. For simplicity we may assume that $M_P = 0$.

For each $0\leq i\leq j$, consider the maps $\TT_{P_i} \Right1{\alpha_i} \TT_P$ and $\TT_{P_i} \Right1{\beta_{i,j}} \TT_{P_j}$, defined by $\alpha_i(\rho) = \incl_{P_i}^P \circ \rho$ and $\beta_{i,j}(\rho) = \incl_{P_i}^{P_j} \circ \rho$, respectively. Notice that $\alpha_i = \alpha_j \circ \beta_{i,j}$ for all $0\leq i\leq j$. Moreover, we claim that $\TT_P$ is the colimit of the system $\{\TT_{P_i}, \beta_{i,j}\}$. Indeed, surjectivity of $\colim \TT_{P_i} \to \TT_P$ follows from the discussion above. To prove injectivity, fix some $i \geq 0$, and let $\rho, \rho' \in \TT_{P_i}$ be such that $\alpha_i(\rho) = \alpha_i(\rho')$. Let also $f, f' \in \Hom_{\FF}(E, P_i)$ be representatives of $\rho$ and $\rho'$ respectively. Then, by definition of $\TT_P$ there exists some $x \in P$ such that
$$
\incl_{P_i}^P \circ f' = c_x \circ \incl_{P_i}^P \circ f.
$$
Since $P = \bigcup_{i \geq 0} P_i$, it follows that $x \in P_j$ for some $j \geq i$, and thus $\beta_{i,j}(\rho) = \beta_{i,j}(\rho') \in \TT_{P_j}$.

For each $i \geq 0$ and each $\rho \in \TT_P$, set $\TT_{P_i}^{\rho} = \alpha_i^{-1}(\rho)$. We also fix the following.
\begin{enumerate}[(a)]

\item For each $f \in \Omega$, let $f_i \in \Hom_{\FF}(E,P_i)$ be the restriction of $f$.

\item For each $\rho \in \TT_P$, fix a set $\Omega_i^{\rho} \subseteq \Hom_{\FF}(E, P_i)$ of representatives of the classes in $\TT_{P_i}^{\rho}$. In particular, if $f \in \Omega$ represents the class $\rho \in \TT_P$, then we choose $f_i$ as representative of its own class in $\TT_{P_i}^{\rho}$. Let also $\Omega_i = \coprod_{\rho \in \TT_P} \Omega_i^{\rho}$. By definition there is an equality
\begin{equation}\label{TEPi}
\4{T}_E(P_i) \defin \bigoplus_{\gamma \in \TT_{P_i}} T_E(H^{\ast}(BP_i);\gamma^{\ast}) = \bigoplus_{\omega \in \Omega_i} T_E(H^{\ast}(BP_i); \omega^{\ast}).
\end{equation}

\item For each $\rho \in \TT_P$, each $i \geq 0$, and each $\omega \in \Omega_i^{\rho}$, fix an element $x_{\omega} \in P_{i+1}$ such that $c_{x_{\omega}} \circ \incl_{P_i}^{P_{i+1}} \circ \omega \in \Omega_{i+1}^{\rho}$ (such an element must exist since $\beta_{i,i+1}[\omega] \in \TT_{P_{i+1}}^{\rho}$). In the particular case where $\omega = f_i$ (see (a) above), we choose $x_{\omega} = 1$. Although the element $x_{\omega}$ clearly depends on $i$, we omit this dependence from the notation since it will be clear at all times which $i$ is involved.

\end{enumerate}
Note that, since $x_{\omega} \in P_{i+1}$, we have
\begin{equation}\label{samemap}
\Map(BE, BP_{i+1})_{\4{\omega}} = \Map(BE, BP_{i+1})_{\sigma},
\end{equation}
where $\4{\omega} = \incl_{P_i}^{P_{i+1}} \circ \omega$ and $\sigma = c_{x_{\omega}} \circ \incl_{P_i}^{P_{i+1}} \circ \omega$.

Fix $\rho \in \TT_P$. For each $i \geq 0$ and each $\omega \in \Omega_i^{\rho}$, consider the homomorphism
$$
\Gamma_{\omega} \colon C_{P_i}(\omega E) \times E \Right2{} P_i,
$$
defined by $\Gamma_{\omega}(a,y) = a \cdot \omega(y) = \omega(y) \cdot a$. Fix $i \geq 0$ and $\omega \in \Omega_i^{\rho}$, and let $x_{\omega}$ be as fixed in (c) above. Let also $\4{\omega} = \incl_{P_i}^{P_{i+1}} \circ \omega$ and $\sigma = c_{x_{\omega}} \circ \incl_{P_i}^{P_{i+1}} \circ \omega$. Then there is a commutative diagram
$$
\xymatrix@C=2cm{
C_{P_i}(\omega E) \times E \ar[d]_{(\incl, \Id)} \ar[r]^{\Gamma_{\omega}} & P_i \ar[d]^{\incl} \\
C_{P_{i+1}}(\4{\omega} E) \times E \ar[d]^{\cong}_{(c_{x_{\omega}}, \Id)} \ar[r]^{\Gamma_{\4{\omega}}} & P_{i+1} \ar[d]_{\cong}^{c_{x_{\omega}}}\\
C_{P_{i+1}}(\sigma E) \times E \ar[r]_{\Gamma_{\sigma}} & P_{i+1}
}
$$
which in turn induces the commutative diagram below by first passing to classifying spaces and then applying adjunction.
\begin{equation}\label{diagrams}
\vcenter{
\xymatrix@C=2cm{
BC_{P_i}(\omega E) \ar[r]^{\simeq} \ar[d]_{B\incl} & \Map(BE, BP_i)_{\omega} \ar[d]^{B\incl_{\ast}}\\
BC_{P_{i+1}}(\4{\omega} E) \ar[r]^{\simeq} \ar[d]_{Bc_{x_{\omega}}}^{\cong} & \Map(BE, BP_{i+1})_{\4{\omega}} \ar[d]_{\cong}^{(Bc_{x_{\omega}})_{\ast}}\\
BC_{P_{i+1}}(\sigma E) \ar[r]_{\simeq} & \Map(BE, BP_{i+1})_{\sigma}
}
}
\end{equation}
Note that the horizontal maps in the diagram above are homotopy equivalences by Step 1. Moreover, recall from (\ref{samemap}) that $\Map(BE, BP_{i+1})_{\4{\omega}} = \Map(BE, BP_{i+1})_{\sigma}$.

Let $\rho \in \TT_P$, and let $f \in \Omega$ be its representative. We claim that there is an isomorphism
\begin{equation}\label{isoTP}
H^{\ast}(BC_P(f E)) \cong \varprojlim_i \big( \bigoplus_{\omega \in \Omega_i^{\rho}} H^{\ast}(BC_{P_i}(\omega E)) \big),
\end{equation}
where the limit is defined by the homomorphisms
$$
H_{\omega} \colon BC_{P_i}(\omega E) \Right2{B\incl} BC_{P_{i+1}}(\4{\omega} E) \Right2{Bc_{x_{\omega}}} BC_{P_{i+1}}(\sigma E).
$$
Note that the above limit does not depend on the choice of the elements $x_{\omega}$ in (c), since a different choice would differ from $x_{\omega}$ by an element in $C_{P_{i+1}}(\sigma E)$.

For each $i \geq 0$, let $f_i \in \Hom_{\FF}(E, P_i)$ be the restriction of $f$ as fixed in (a) above, which is an element of $\Omega_i^{\rho}$ by (b), and let $\3{\Omega}_i^{\rho} = \Omega_i^{\rho} \setminus \{f_i\}$. For each $\omega \in \Omega_i^{\rho}$, set for short $A_{\omega} = H^{\ast}(BC_{P_i}(\omega E))$. Then, there are short exact sequences
$$
0 \to \bigoplus_{\omega \in \3{\Omega}_i^{\rho}} A_{\omega} \Right2{\iota_i} \bigoplus_{\omega \in \Omega_i^{\rho}} A_{\omega} \Right2{\pi_i} A_{f_i} \to 0
$$
for all $i \geq 0$. Let also $\kappa_{i+1} \colon \bigoplus_{\sigma \in \Omega_{i+1}^{\rho}} A_{\sigma} \to \bigoplus_{\omega \in \Omega_i^{\rho}} A_{\omega}$ be the morphism induced by the maps $H_{\omega}$ described above, and let $\3{\kappa}_{i+1}$ be the restriction of $\kappa_{i+1}$ to $\bigoplus_{\sigma \in \3{\Omega}_{i+1}^{\rho}} A_{\sigma}$. In order to prove (\ref{isoTP}) we use the above exact sequences to construct an exact sequence of inverse limits.

Recall that $\TT_P$ is the colimit of the sets $\TT_{P_i}$. Thus, given $i \geq 0$, there exists some $M \in \N$ such that $[\incl_{P_i}^{P_{i+M}} \circ \omega] = [f_{i + M}] \in \TT_{P_{i+M}}^{\rho}$ for all $\omega \in \Omega_i^{\rho}$, where $f_{i+M}$ is the restriction of $f$ to $P_{i+M}$ as fixed in (a) above. For simplicity let us assume that $M = 1$. In terms of mapping spaces, this means that, composing with the inclusion map $BP_i \to BP_{i+1}$, we have
$$
\coprod_{\omega \in \Omega_i^{\rho}} \Map(BE, BP_i)_{\omega} \Right2{} \Map(BE, P_{i+1})_{f_{i+1}},
$$
and it follows that the morphism $\3{\kappa}_{i+1}$, defined in the previous paragraph, is simply the trivial homomorphism. The morphism $\kappa_{i+1}$ also induces a morphism $A_{f_{i+1}} \to A_{f_i}$, which is easily seen to coincide with the restriction homomorphism induced by the inclusion $C_{P_i}(f_iE)\leq C_{P_{i+1}}(f_{i+1}E)$. Summarizing, we obtain commutative diagrams
$$
\xymatrix@C=1.5cm{
0 \ar[r] & \bigoplus_{\sigma \in \3{\Omega}_{i+1}^{\rho}} A_{\sigma} \ar[d]_{\3{\kappa}_{i+1}= 0} \ar[r]^{\iota_{i+1}} & \bigoplus_{\sigma \in \Omega_{i+1}^{\rho}} A_{\sigma} \ar[r]^{\pi_{i+1}} \ar[d]^{\kappa_{i+1}} & A_{f_{i+1}} \ar[r] \ar[d]^{\res_{i+1}} & 0 \\
0 \ar[r] & \bigoplus_{\omega \in \3{\Omega}_i^{\rho}} A_{\omega} \ar[r]_{\iota_i} & \bigoplus_{\omega \in \Omega_i^{\rho}} A_{\omega} \ar[r]_{\pi_i} & A_{f_i} \ar[r] & 0
}
$$
This produces a morphism of inverse systems, and hence, since the inverse limit functor is left exact, there is an exact sequence
\begin{equation}\label{exactlim}
0 \to \varprojlim_i \big(\bigoplus_{\omega \in \3{\Omega}_i^{\rho}} A_{\omega} \big) \Right1{} \varprojlim_i \big( \bigoplus_{\omega \in \Omega_i^{\rho}} A_{\omega} \big) \Right2{\pi} \varprojlim_i A_{f_i},
\end{equation}
where $\varprojlim_i \big(\bigoplus_{\omega \in \3{\Omega}_i^{\rho}} H^{\ast}(BC_{P_i}(\omega E))\big) = 0$ since the morphisms in the corresponding inverse system, namely $\3{\kappa}_{i+1}$, are trivial. Thus, $\pi$ is a monomorphism, and it remains to prove the it is also surjective. Let $(x_i)_{i \in \N}$ be an element in $\varprojlim_i H^{\ast}(BC_{P_i}(f_i E))$. That is, $\res_{i+1}(x_{i+1}) = x_i$ for all $i \geq 0$. Let also
$$
X_i = \pi_i^{-1}(\{x_i\}) \subseteq \bigoplus_{\omega \in \Omega_i^{\rho}} H^{\ast}(BC_{P_i}(\omega E)).
$$
Restricting the maps $\kappa_{i+1}$ to the sets $X_{i+1}$ produces an inverse system of sets, and the resulting inverse limit $X = \varprojlim_i X_i$ is nonempty by \cite[Proposition 1.1.4]{RZ}. Moreover, by definition of the maps $\kappa_i$, any element in $X$ is a preimage of $(x_i)_{i \in \N}$, and $\pi$ is surjective. Thus, we have
$$
H^{\ast}(BC_P(f E)) \cong \varprojlim_i H^{\ast}(BC_{P_i}(f_i E)) \cong \varprojlim_i \big(\bigoplus_{\omega \in \Omega_i^{\rho}} H^{\ast}(BC_{P_i}(\omega E))\big),
$$
where the leftmost isomorphism follows from Proposition \ref{stable1}, and the rightmost isomorphism corresponds to $\pi$ in (\ref{exactlim}). This proves (\ref{isoTP}).

Combining Step 1 and the above discussion, we deduce the following.
$$
\begin{aligned}
\varprojlim_i \4{T}_E(P_i) & = \varprojlim_i \big( \bigoplus_{\omega \in \Omega_i} T_E(H^{\ast}(BP_i); \omega^{\ast}) \big) \cong \varprojlim_i \big( \bigoplus_{\omega \in \Omega_i} H^{\ast}(\Map(BE, BP_i)_{\omega}) \big) \cong \\
 & \cong \varprojlim_i \big( \bigoplus_{\omega \in \Omega_i} H^{\ast}(BC_{P_i}(\omega E)) \big) \cong \bigoplus_{f \in \Omega} H^{\ast}(BC_P(f E)) \cong \\
 & \cong \bigoplus_{f \in \Omega} H^{\ast}(\Map(BE, BP)_{f}) \cong \bigoplus_{f \in \Omega} T_E(H^{\ast}(BP); f^{\ast}) = \4{T}_E(P).
\end{aligned}
$$
More precisely, the first equality holds by (\ref{TEPi}), while the last equality corresponds to (\ref{TEP}), the first isomorphism in the first row follows from Step 1, the isomorphism between the last term in the first row and the first term in the second row follows from Step 1 together with the leftmost diagram in (\ref{diagrams}), the middle isomorphism in the second row holds by (\ref{isoTP}), and the rest follows again from Step 1.

\textbf{Step 5.} The isomorphism $T_E(H^{\ast}(\FF); j_E) \Right3{\cong} H^{\ast}(C_{\FF}(E))$. Consider the fusion system $C_{\FF}(E)$ over $C_S(E)$, and consider also the set of fusion subsystems $\{C_{\FF_i}(E)\}_{i \geq 0}$. This setup satisfies the conditions of Lemma \ref{aux34}: clearly, $C_S(E) = \bigcup_{i \geq 0} C_{S_i}(E)$, since $S = \bigcup_{i \geq 0} S_i$, and every morphism in $C_{\FF}(E)$ eventually restricts to a morphism in $C_{\FF_i}(E)$ for $i$ big enough, since the fusion systems $\FF_i$ are part of an approximation of $\ploc$ by $p$-local finite groups. It follows that there is a sequence of isomorphisms
$$
T_E(H^{\ast}(\FF); j_E) \cong \varprojlim_{\FF} \4{T}_E(-) \cong \varprojlim_i \big(\varprojlim_{\FF_i} \4{T}_E(-) \big) \cong \varprojlim_i H^{\ast}(C_{\FF_i}(E)) \cong H^{\ast}(C_{\FF}(E))
$$
where the first isomorphism follows from Steps 3 and 4 combined, the second isomorphism is a consequence of Lemma \ref{aux34}, since we have checked in Step 4 that $\4{T}_E$ satisfies condition (\textasteriskcentered), the third isomorphism follows from \cite[Lemma 5.7]{BLO2}, since we have shown in Step 2 that $E$ can be assumed to be fully $\FF_i$-centralized for all $i$, and the last isomorphism is a consequence of Proposition \ref{stable1}.
\end{proof}

\begin{rmk}

Let $\g = \ploc$ be a $p$-local compact group, and let $\{(S_i, \FF_i, \LL_i)\}_{i \geq 0}$ be an approximation of $\g$ by $p$-local finite groups. Let also $E\leq S$ be an elementary abelian subgroup which is fully centralized, and let $j_E$ be as above. Without loss of generality we may assume that $E\leq S_i$ for all $i \geq 0$, and that $E$ is fully $\FF_i$-centralized for all $i \geq 0$. Let $j_{E,i} \colon H^{\ast}(\FF_i) \to H^{\ast}(BE)$ be the map induced by the inclusion $E\leq S_i$, for all $i \geq 0$. In this situation, we have just proved that
$$
T_E(H^{\ast}(\FF); j_E) \cong T_E(\varprojlim_i H^{\ast}(\FF_i); j_E) \cong \varprojlim_i T_E(H^{\ast}(\FF_i);j_{E,i}).
$$
In other words, the functor $T_E$ commutes with the inverse limit $\varprojlim_i H^{\ast}(\FF_i)$. We do not know of any general result about the functor $T_E$ commuting with infinite inverse limits.

\end{rmk}

The proof of Theorem \ref{mapping} below requires the following result involving the quotient of a linking system by a normal discrete $p$-toral subgroup. The following is a generalization of \cite[Lemma 5.6]{BLO2} to the compact case.

\begin{lmm}\label{quotient2}

Let $\g = \ploc$ be a $p$-local compact group, and let $\wL$ be the telescopic linking system associated to $\LL$ in \ref{expl1}. Let also $V\leq S$ be a central subgroup in $\FF$ of order $p$, and let $(\wL/V, \3{\varepsilon}, \3{\rho})$ be the quotient transporter system, associated to the saturated fusion system $\FF/V$. Finally, let $(\wL/V)^c \subseteq \wL/V$ and $\wL_0 \subseteq \wL$ be the full subcategories with object sets
$$
\begin{array}{l}
\Ob((\wL/V)^c) = \{P/V \in \Ob(\wL/V) \, | \, P/V \mbox{ is $\FF/V$-centric}\} \\[2pt]
\Ob(\wL_0) = \{P \in \Ob(\wL) \, | \, P/V \in \Ob((\wL/V)^c)\}
\end{array}
$$
respectively. Then the following holds.
\begin{enumerate}[(i)]

\item $(\wL/V)^c$ is a linking system associated to $\FF/V$.

\item $BV \to |\wL_0|^{\wedge}_p \to |(\wL/V)^c|^{\wedge}_p$ is a fibration sequence.

\item The inclusion $\wL_0 \subseteq \wL$ induces a homotopy equivalence $|\wL_0|^{\wedge}_p \simeq B\g$.

\end{enumerate}

\end{lmm}

\begin{proof}

Clearly, all the $\FF$-centric subgroups of $S$ must contain $V$, since $V$ is abelian and $\FF$-central. Thus, for all $P, Q \in \Ob(\wL)$, it follows by axiom (C) of transporter systems that the left and right actions of $V$ on $\Mor_{\wL}(P,Q)$ (via composition with $\varepsilon_Q(V)$ and $\varepsilon_P(V)$ respectively) are the same. Since the proof is rather long, it is divided into steps for the reader's convenience.

\textbf{Step 1.} For each $P\leq S$ which contains $V$, we claim that
$$
P \mbox{ fully $\FF$-normalized } \Longrightarrow P/V \mbox{ fully $\FF/V$-centralized } \Longrightarrow \Gamma_P\leq \Aut_S(P).
$$

Let $P\leq S$ be such that $V\leq P$, and note that $\Aut_{\FF}(V) = \{\Id\}$, since $V$ is central in $\FF$. Thus, the subgroup
$$
\Gamma_P \defin \Ker(\Aut_{\FF}(P) \Right2{} \Aut_{\FF/V}(P/V))
$$
is a discrete $p$-toral normal subgroup of $\Aut_{\FF}(P)$ by Lemma \ref{Kpgp}. Since $\FF$ is saturated, every subgroup of $S$ is $\FF$-conjugate to a fully $\FF$-normalized subgroup. Similarly, since $\FF/V$ is saturated, every subgroup of $S/V$ is $\FF/V$-conjugate to a fully $\FF/V$-centralized subgroup. Thus it is enough to show the following: if $P, Q\leq S$ are $\FF$-conjugate subgroups such that $P/V$ is fully $\FF/V$-centralized and $Q$ is fully $\FF$-normalized, then $Q/V$ is fully $\FF/V$-centralized and $\Gamma_P\leq \Aut_S(P)$.

As show above, the group $\Gamma_Q$ is a normal discrete $p$-toral subgroup of $\Aut_{\FF}(Q)$. Furthermore, since $Q$ is fully $\FF$-normalized we have $\Aut_S(Q) \in \Syl_p(\Aut_{\FF}(Q))$, and thus $\Gamma_Q\leq \Aut_S(Q)$. By axiom (II) of saturated fusion systems, every isomorphism $f \in \Iso_{\FF}(P,Q)$ extends to some $\gamma \in \Hom_{\FF}(N_f, N_S(Q))$, where
$$
N_f = \{g \in N_S(P) \, | \, f \circ c_g \circ f^{-1} \in \Aut_S(Q)\}.
$$
Set $N_S^0(P) = \{g \in N_S(P) \, | \, c_g \in \Gamma_P\}$, and notice that $N_S^0(P)/V = C_{S/V}(P/V)$. We claim that $N_S^0(P)\leq N_f$. To prove that, fix $g \in N_S^0(P)$ and $a \in Q$, and set $b = f^{-1}(a) \in P$. By definition, $c_g(b) = bv$ for some $v \in V$, and we have
$$
(f \circ g \circ f^{-1})(a) = (f \circ c_g)(b) = f(bv) = f(b) v = a v,
$$
where $f(bv) = f(b) v$ since $V$ is central in $\FF$ and $V,\leq P,Q$ (and thus the morphism $f$ restricts to the identity on $V$).

The above implies that $\gamma$ restricts to $\gamma \in \Hom_{\FF}(N_S^0(P), N_S^0(Q))$, which in turn factors through a homomorphism
$$
\3{\gamma} \in \Hom_{\FF/V}(C_{S/V}(P/V), C_{S/V}(Q/V)).
$$
Since $P/V$ is fully $\FF/V$-centralized, it follows that $\3{\gamma}$ is an isomorphism and $Q/V$ is also fully $\FF/V$-centralized. Furthermore, $\gamma$ must be an isomorphism too, and thus $\Gamma_P\leq \Aut_S(P)$.

\textbf{Step 2.} We show now that the category $(\wL/V)^c$ is a centric linking system associated to $\FF/V$. First notice that $(\wL/V)^c$ is a transporter system associated to $\FF/V$ by the above remarks, and thus we only have to check that, for each object $P/V$ of $(\wL/V)^c$,
$$
E(P/V) \defin \Ker(\Aut_{(\wL/V)^c}(P/V) \to \Aut_{\FF/V}(P/V)) = \3{\varepsilon}_{P/V}(Z(P/V)).
$$

Fix $P/V \in \Ob((\wL/V)^c)$, and consider the following commutative diagram
$$
\xymatrix{
\Aut_{\wL}(P) \ar[rr]^{\tau_P} \ar[d]_{\rho_P} & & \Aut_{\wL/V}(P/V) \ar[d]^{\3{\rho}_{P/V}} \\
\Aut_{\FF}(P) \ar[rr]_{\omega_P} & & \Aut_{\FF/V}(P/V),
}
$$
where $\rho_{P}$ and $\3{\rho}_{P/V}$ denote the corresponding structural functors in the transporter systems $\wL$ and $\wL/V$, respectively. By definition, we have
$$
\Ker(\tau_P) = \varepsilon_P(V) \qquad \qquad \Ker(\rho_P) = \varepsilon_P(C_S(P)) \qquad \qquad \Ker(\omega_P) = \Gamma_P,
$$
and an easy computation shows then that $E(P/V) = \tau_P(\gen{\varepsilon_P(C_S(P)), \varepsilon_P(N_S^0(P))})$. In particular it follows that $E(P/V)$ is a discrete $p$-toral group. To finish the proof, recall that $P/V$ is $\FF/V$-centric, and in particular it is fully $\FF/V$-centralized. Thus $Z(P/V) \in \Syl_p(E(P/V))$, which implies that $E(P/V) = Z(P/V)$.

\textbf{Step 3.} $BV \to |\wL_0|^{\wedge}_p \to |(\wL/V)^c|^{\wedge}_p$ is a fibration sequence. Using \cite[Lemma 4.3 (a)]{BLO3} it is easy to check that each undercategory for the projection of $\wL_0$ onto $(\wL/V)^c$ contains a category equivalent to $\bb(V)$ as a deformation retract. Thus, by Quillen's Theorem B, the map $|\wL_0| \to |(\wL/V)^c|$ has homotopy fiber $BV$. By \cite[II.5.1]{BK}, the fibration sequence $BV \to |\wL_0| \to |(\wL/V)^c|$ is still a fibration sequence after $p$-completion.

\textbf{Step 4.} The inclusion $\wL_0 \subseteq \wL$ induces a homotopy equivalence $|\wL_0|^{\wedge}_p \simeq B\g$. Notice that the functor $\bullf$ restricts to a functor $\bullf$ on $\wL_0$. Indeed, let $P \in \Ob(\wL^{\bullet})\setminus \Ob(\wL_0)$, and let $Q \in \Ob(\wL)$ be such that $V\leq Q$ and $Q^{\bullet} = P$. Then $Q/V\leq P/V$, and since $P/V$ is not $\FF/V$-centric by assumption, neither is $Q/V$. Thus $Q \notin \Ob(\wL_0)$. Conversely, if $Q \in \Ob(\wL_0)$, then clearly $Q^{\bullet}\in \Ob(\wL_0)$.

Also, note that if $\wL_0$ contains $\LL^{\bullet}$ then the claim follows easily. Thus, fix some $P\leq S$ in $\LL^{\bullet}$ but not in $\wL_0$. In particular, $P$ is $\FF$-centric, and $V\leq P$, but $P/V$ is not $\FF/V$-centric. By replacing $P$ by a conjugate if necessary, we may assume that $C_{S/V}(P/V) \neq Z(P/V)$, and thus there is some $gV \in S/V$ such that $gV \notin P/V$ and $[gV, P/V] = 1$. Equivalently, there is some $g \in S$ such that $g \notin P$ and $[g, P] = V$. Furthermore, $c_g \in \Aut_{\FF}(P)$ cannot be an inner automorphism, because if this was the case then $c_g = c_x$ for some $x \in P$, and $gx^{-1} \in C_S(P) \setminus P = \emptyset$, since $P$ is $\FF$-centric. Thus $c_g$ is a nontrivial element of $\Ker(\Out_{\FF}(P) \to \Out_{\FF/V}(P/V))$, and thus $P$ is not $\FF$-radical. In other words, we have shown that $\wL_0$ contains all $\FF$-centric $\FF$-radical subgroups. By \cite[Corollary A.10]{BLO6} the inclusion $\wL_0 \subseteq \wL$ induces a homotopy equivalence $|\wL_0|^{\wedge}_p \simeq B\g$.
\end{proof}

Let $\g = \ploc$ be a $p$-local compact group, and let $P\leq S$ be a fully $\FF$-centralized subgroup. Let also $C_{\g}(P) = (C_S(P), C_{\FF}(P), C_{\LL}(P))$ be the centralizer $p$-local compact group of $P$ defined in \ref{rmknorm}, with classifying space $BC_{\g}(P)$, and let $\bb(P)$ be the category with a single object $\circ_P$ and $P$ as automorphism group. By Lemma \ref{centricNFKA}, if $Q\leq C_S(P)$ is $C_{\FF}(P)$-centric then $QP$ is $\FF$-centric. Thus we can define a functor
$$
\Gamma_{\LL, P} \colon C_{\LL}(P) \times \bb(P) \Right3{} \LL
$$
by setting $\Gamma_{\LL,P}(Q, \circ_P) = QP$ for each $C_{\FF}(P)$-centric subgroup of $C_S(P)$. Given a morphism $(\varphi, g) \in \Mor_{C_{\LL}(P) \times \bb(P)}((Q,\circ), (R,\circ))$, the functor $\Gamma$ is defined as follows
$$
\Gamma_{\LL,P}(\varphi, g) = \varphi \circ \varepsilon_{QP}(g) = \varepsilon_{RP}(g) \circ \varphi,
$$
where the last equality follows from condition (C) of transporter systems, since the underlying homomorphism of $\varphi \in \Mor_{\LL}(QP, RP)$ restricts to the identity on $P$ by definition of $C_{\LL}(P)$. By first realizing nerves and them $p$-completing, we obtain a map
$$
BC_{\g}(P) \times (BP)^{\wedge}_p \Right3{} B\g
$$
(notice that $(BP)^{\wedge}_p$ is not necessarily equivalent to $BP$ since $P$ is a discrete $p$-toral group). By taking adjoint first and then precomposing with the natural map $BP \to (BP)^{\wedge}_p$, we obtain a map
$$
\Gamma'_{\LL,P} \colon |C_{\LL}(P)|^{\wedge}_p \Right3{} \Map((BP)^{\wedge}_p, B\g)_{\incl} \Right3{} \Map(BP, B\g)_{\incl}.
$$

\begin{thm}\label{mapping}

Let $\g = \ploc$ be a $p$-local compact group, let $P$ be a discrete $p$-toral group, and $\gamma: P \to S$ be a group homomorphism such that $\gamma(P)$ is fully $\FF$-centralized in $\FF$. Then,
$$
\Gamma'_{\LL, \gamma(P)} \colon BC_{\g}(\gamma(P)) \Right2{\simeq} \Map(BP, B\g)_{B\gamma}
$$
is a homotopy equivalence.

\end{thm}

\begin{proof}

Our proof follows the same strategy as the proof of \cite[Theorem 6.3]{BLO2}. The referee suggested an alternative to the cases 1-3 below, which we briefly discuss after the proof. By \cite[Proposition 6.2]{BLO3}, for each $\gamma \in \Hom(P, S)$,
\begin{equation}\label{zero}
\Map(BP, B\g)_{B\gamma} \simeq \Map(B\gamma(P), B\g)_{\incl}.
\end{equation}
Thus, it suffices to prove the statement when $P\leq S$ is fully $\FF$-centralized and $\gamma$ is the inclusion. The proof is divided into several cases for the reader's convenience.

\textbf{Case 1.} Suppose that $P$ is elementary abelian. In this case, by Theorem 0.5 \cite{Lannes}, it is enough to prove that $\Gamma_{\LL,P}$ induces an isomorphism
$$
T_P(H^{\ast}(B\g); \incl^{\ast}) \Right2{\cong} H^{\ast}(BC_{\g}(P)).
$$
By Theorem \ref{stable2}, $H^{\ast}(B\g) \cong H^{\ast}(\FF)$, and $H^{\ast}(BC_{\g}(P)) \cong H^{\ast}(C_{\FF}(P))$, and the above isomorphism follows from Lemma \ref{mapping1}.

\textbf{Case 2.} Suppose that $P$ is a normal subgroup of $\FF$, that is, $N_{\FF}(P) = \FF$. Let $\LL_0 \subseteq \LL$ be the full subcategory whose objects are the subgroups $Q\leq S$ such that $C_Q(P)$ is $C_{\FF}(P)$-centric.

In order to prove the statement in this case, we first need to show that the inclusion of $\LL_0$ into $\LL$ induces an equivalence $|\LL_0|^{\wedge}_p \simeq |\LL|^{\wedge}_p$. By \cite[Corollary A.10]{BLO6}, it is enough to check that $\Ob(\FF^{cr}) \subseteq \Ob(\LL_0)$.

Let $Q \in \Ob(\LL)$ be an $\FF$-centric subgroup of $S$ which is not an object in $\LL_0$. We claim that $Q$ is not $\FF$-radical. Set $Q_0 = C_Q(P)$, and note that every element of $\Aut_{\FF}(Q)$ restricts to an automorphism of $Q_0$ since $P$ is normal in $\FF$. Furthermore $Q_0 \lhd Q$, and we can define
$$
K = \Ker(\Aut_{\FF}(Q) \Right2{} \Aut_{\FF}(Q_0) \times \Aut(Q/Q_0)) \lhd \Aut_{\FF}(Q).
$$
Note that $K$ is a discrete $p$-toral subgroup of $\Aut_{\FF}(Q)$ by Lemma \ref{Kpgp}. In order to prove that $Q$ is not $\FF$-radical, it is enough to check that $1 \neq K \not\leq \Inn(Q)$.

By assumption $Q_0 \notin \Ob(\LL_0)$, and thus $Q_0$ is not $C_{\FF}(P)$-centric. Thus, we may assume that $C_{C_S(P)}(Q_0) \not\leq Q_0$, since otherwise $Q$ can be replaced by an $\FF$-conjugate $R$ such that the corresponding subgroup $R_0 = C_R(P)$ is $C_{\FF}(P)$-conjugate to $Q_0$ and satisfies the desired condition. Set
$$
Q_1 \defin C_S(Q_0P) = C_{C_S(P)}(Q_0).
$$
As discussed above we have $Q_1 \not\leq Q_0 = C_Q(P)$, and thus, $Q_1 \cap Q = C_Q(Q_0P)\leq Q_0$, and $Q_1 \not\leq Q$. Also, $Q\leq N_S(Q_1)$ by definition, and thus $Q_1Q\leq S$ is a subgroup. Note that $Q \lneqq Q_1Q$, and thus $Q \lneqq N_{Q_1Q}(Q)$. Choose some $x \in N_{Q_1Q}(Q)$ such that $x \notin Q$. Then,
$$
[x, Q] = \{xax^{-1}a^{-1} \,\, | \,\, a \in Q\}\leq Q_1 \cap Q\leq Q_0,
$$
and hence $c_x \in K$. Since $x \notin Q$ and $Q$ is $\FF$-centric, it follows that $c_x \notin \Inn(Q)$, and thus $Q$ is not $\FF$-radical.

Let $(\LL_0)_{P, Id}$ be the category with object set the pairs $(Q, \alpha)$, for $Q$ in $\LL_0$ and $\alpha \in \Hom_{\FF}(P, Q)$ and such that
$$
\Mor_{(\LL_0)_{P, Id}}((Q,\alpha), (R, \alpha')) = \{\varphi \in \Mor_{\LL}(Q,R) \mbox{ } | \mbox{ } \alpha' = \rho(\varphi) \circ \alpha\}.
$$
This is equivalent to the component of the object $(P,Id)$ in the category $\LL_0^P$ of \cite[Proposition 6.2]{BLO3}, and hence there is a homotopy equivalence
$$
\Map (BP, |\LL_0|^{\wedge}_p)_{\incl} \simeq |(\LL)_{P, Id}|^{\wedge}_p.
$$

At this point, one can define functors
$$
\xymatrix{
(\LL_0)_{P,Id} \ar @<  2pt> [r]^{\sigma} & C_{\LL}(P) \ar @< 2pt> [l]^{\tau}\\
}
$$
in the same way as they are constructed in Step 2 of the proof \cite[Theorem 6.3]{BLO2}, and which are inverse to each other up to natural transformation. This implies that the composite
$$
BC_{\g}(P) \Right2{\tau} |(\LL_0)_{P, Id}|^{\wedge}_p \Right2{\simeq} \Map(BP, |\LL_0|^{\wedge}_p)_{\incl} \Right2{\simeq} \Map(BP, B\g)_{\incl}
$$
is a homotopy equivalence, and by construction it is equal to $\Gamma'_{\LL,P}$.

\textbf{Case 3.} Suppose that $P$ is a finite subgroup of $S$. As an induction hypothesis, we can assume that the statement holds for all maps with source $BP'$, with $|P'| < |P|$, and all $p$-local compact groups.

Fix a subgroup $V\leq  P \cap Z(N_S(P))$ of order $p$, and note that in particular $N_S(P)\leq C_S(V)$. By \cite[Lemma 2.2 (b)]{BLO6}, there exists some $\omega \in \Hom_{\FF}(N_S(V), S)$ such that $\omega(V)$ is fully $\FF$-centralized. Hence, there is an inequality $|N_S(\omega(P))| \geq |N_S(P)|$ which is in fact an equality since we are assuming $P$ to be fully $\FF$-normalized.

We may replace $P$ and $V$ by $\omega(P)$ and $\omega(V)$, and assume that $V$ is fully $\FF$-centralized and $P$ is fully $\FF$-normalized. Furthermore, $P$ is fully normalized in the saturated fusion system $C_{\FF}(V)$, since $N_S(P) = N_{C_S(V)}(P)$.

By Case 1, the map $BC_{\g}(V) \to B\g$, induced by the inclusion $C_{\LL}(V) \subseteq \LL$, induces a homotopy equivalence $\Map(BV, BC_{\LL}(V))_{\incl} \simeq \Map(BV, B\g)_{\incl}$, and hence also a homotopy equivalence
$$
\Map((EP)/V, BC_{\g}(V))_{\incl} \Right2{\simeq} \Map((EP)/V, B\g)_{\incl}
$$
which is $P/V$-equivariant (where the action of $P/V$ is the action induced by the original action of $P$ on $EP$). This is still a homotopy equivalence after considering homotopy fixed point sets by \cite[Remark 10.2]{DW0}, and thus we obtain another homotopy equivalence
$$
\big[\Map((EP)/V, BC_{\g}(V))_{\incl}\big]^{h(P/V)} \Right2{\simeq} \big[\Map((EP)/V, B\g)_{\incl}\big]^{h(P/V)}.
$$
Notice that $E(P/V) \times_{P/V} (EP)/V \simeq BP$ with the given actions (here $P/V$ is acting diagonally on $E(P/V) \times (EP)/V$). Let $X = BC_{\g}(V)$ or $B\g$. By definition of homotopy fixed point sets, we have
$$
\begin{aligned}
\big[\Map((EP)/V, X)\big]^{h(P/V)} & = \Map_{P/V}(E(P/V), \Map((EP)/V,X)) \simeq \\
 & \simeq \Map_{P/V}(E(P/V) \times (EP)/V, X),
\end{aligned}
$$
where the rightmost equivalence follows by adjunction. Furthermore, $P/V$ acts trivially on $X$, and thus it follows that $\Map_{P/V}(E(P/V) \times (EP)/V, X) \simeq \Map(E(P/V) \times_{P/V} (EP)/V, X) \simeq \Map(BP,X)$. Thus, the equivalence of homotopy fixed point sets above induces the following equivalence
$$
\Map(BP, BC_{\g}(V))_{\incl} \Right2{\simeq} \Map(BP, B\g)_{\incl}.
$$
We can suppose that $\LL = C_{\LL}(V)$, and hence that $V$ is central in $\LL$. Let $\g/V$ be the quotient of $\g$ by $V$, as described in Definition \ref{quotient1}. In particular, $\FF/V$ is a saturated fusion system on $S/V$, and $\LL/V$ is a transporter system.

Consider the full subcategories $\LL_0 \subseteq \LL$ and $(\LL/V)^c \subseteq \LL/V$ whose objects are the subgroups $Q\leq S$, respectively $Q/V\leq S/V$, such that $Q/V$ is $\FF/V$-centric. In particular, $(\LL/V)^c$ determines a centric linking system associated to $\FF/V$, and there are homotopy equivalences $|\LL_0|^{\wedge}_p \simeq B\g$ and $|(\LL/V)^c|^{\wedge}_p \simeq |\LL/V|^{\wedge}_p$.

Finally, let also $\FF' = N_{\FF}(P)$, $\LL' = N_{\LL}(P)$, and define $\LL'/V$, $\LL_0' \subseteq \LL'$ in a similar way as done above. It follows by Lemma \ref{quotient2} that there are fibration sequences
$$
\xymatrix@R=2mm{
BV \ar[rr] & & |\LL_0|^{\wedge}_p \ar[rr]^{\Phi} & & |(\LL/V)^c|^{\wedge}_p \\
BV \ar[rr] & & |\LL_0'|^{\wedge}_p \ar[rr]^{\Phi'} & & |(\LL'/V)^c|^{\wedge}_p,\\
}
$$
and hence also a homotopy pull-back square
$$
\xymatrix{
\Map(BP, |\LL_0'|^{\wedge}_p)_{\iota} \ar[rr]^{I_1} \ar[d]_{\Phi' \circ -} & & \Map(BP, |\LL_0|^{\wedge}_p)_{\incl} \ar[d]^{\Phi \circ -} \\
\Map(BP, |(\LL'/V)^c|^{\wedge}_p)_{\Phi \circ \incl} \ar[rr]_{I_2} & & \Map(BP, |(\LL/V)^c|^{\wedge}_p)_{\Phi \circ \incl},
}
$$
where $\Map(BP, |\LL_0'|^{\wedge}_p)_{\iota}$ is the union of the connected components which map to the inclusion in $|\LL_0|^{\wedge}_p$ and to $\Phi \circ \incl$ in $|(\LL'/V)^c|^{\wedge}_p$, and $I_1$, $I_2$ are inclusions.

By (\ref{zero}), together with the induction hypothesis, and since $P/V$ has strictly smaller order than $P$ (recall that $P$ is finite by hypothesis), there are homotopy equivalences
$$
\begin{aligned}
|C_{(\LL/V)^c}(P/V)|^{\wedge}_p & \RIGHT6{\Gamma'_{\LL/V, P/V}}{\simeq} \Map(B(P/V), |(\LL/V)^c|^{\wedge}_p)_{\incl} \Right1{} \\
 & \RIGHT6{- \circ \proj}{\simeq} \Map(BP, |(\LL/V)^c|^{\wedge}_p)_{f \circ \incl}
\end{aligned}
$$
and similarly for maps to $|(\LL'/V)^c|^{\wedge}_p$. Since $\Gamma'_{\LL/V, P/V}$ is the composite of $\Gamma'_{\LL'/V, P/V}$ with the inclusion by definition , this shows that the map $I_2$ in the diagram above is a homotopy equivalence, and hence so is $I_1$. In particular, $\Map(BP, |\LL_0'|^{\wedge}_p)_{\iota}$ is connected and contains the component of the inclusion. Thus, Case 3 follows from Case 2 applied to the mapping space $\Map(BP, |\LL_0'|^{\wedge}_p)_{\incl}$.

\textbf{Case 4.} Suppose that $P$ is an infinite discrete $p$-toral group. By Lemma \ref{central2}, there is a sequence of subgroups $P_0\leq P_1\leq \ldots$ such that $P = \bigcup_{n \geq 0} P_n$, and such that $P_n$ is fully $\FF$-centralized with $C_{\g}(P_n) = C_{\g}(P)$ for all $n \geq 0$. We have a sequence of homotopy equivalences
$$
\begin{aligned}
\Map(BP, B\g)_{\incl} & = \Map(\hocolim \mbox{ } BP_n, B\g)_{\incl} \simeq \\
& \simeq \holim \mbox{ } \Map(BP_n, B\g)_{\incl} \simeq \holim \mbox{ } BC_{\g}(P_n) = BC_{\g}(P),
\end{aligned}
$$
where the equivalence $\Map(\hocolim \mbox{ } BP_n, B\g)_{\incl} \simeq \holim \mbox{ } \Map(BP_n, B\g)_{\incl}$ follows from \cite[Proposition 2, page 187]{D-F}. This finishes the proof.
\end{proof}

The reader may think of replacing Cases 1 to 3 in the proof above by the following argument. Given a $p$-local compact group $\g = \ploc$ and a fully $\FF$-centralized finite subgroup $P\leq S$, consider the centralizer $p$-local compact group $C_{\g}(P) = (C_S(P), C_{\FF}(P), C_{\LL}(P))$ of $P$ in $\g$. Set for short $Z = C_S(P)$, $\ee = C_{\FF}(P)$, and $\TT = C_{\LL}(P)$. Given a fine unstable Adams operation (see \ref{uAo}), we may assume that $\Psi(P) = P$, since $P$ is a finite subgroup of $S$, and thus $\Psi$ restricts to a fine unstable Adams operation on $C_{\g}(P)$.

This way, $\Psi$ defines approximations of $\g$ and $C_{\g}(P)$ by $p$-local finite groups, namely $\{(S_i, \FF_i, \LL_i)\}_{i \geq 0}$ and $\{(Z_i, \ee_i, \TT_i)\}_{i \geq 0}$ respectively. Moreover, we may assume that $P\leq S_i$ for all $i \geq 0$. It is not hard to see that $Z_i = C_{S_i}(P)$ for all $i \geq 0$, since $Z_i = Z \cap S_i = C_S(P) \cap S_i$. The main difficulty of this argument is that it is not clear whether $\ee_i$ corresponds to the centralizer fusion system of $P$ in $\FF_i$ for any $i$. Essentially, the main problem is that the finite retraction pairs for $\g$ and $C_{\g}(P)$ given in \ref{expl1} do not agree with each other in general, and \ref{expl3} does not apply to this situation, since in general $C_S(P)$ does not have finite index in $S$. One could drop this last condition, but at the price of dealing with a more complicated situation.

\begin{rmk}

With the above description of the homotopy type of the mapping spaces $\map(BP, B\g)$, one could now generalize the results of C. Broto, N. Castellana, J. Grodal, R. Levi and B. Oliver \cite{BCGLO1, BCGLO2}. More precisely, \cite[Theorem A]{BCGLO1} has already been proved for $p$-local compact groups as \cite[Theorem 4.2]{BLO6}, and \cite[Theorem B]{BCGLO1} would follow easily now from our result above. Regarding \cite{BCGLO2}, some results have already been extended to $p$-local compact groups in \cite[Appendix B]{Gonza2}, and the rest would follow by the same arguments (with some minor modifications). We omit this for the sake of brevity of this paper.

\end{rmk}


\bibliographystyle{gtart}
\bibliography{Main}

\end{document}